\documentclass[a4paper,onecolumn,11pt]{article}
\usepackage{geometry}
\usepackage{amsmath}
\usepackage{indentfirst}
\usepackage{ntheorem}
\usepackage{authblk}
\usepackage{bbm}
\usepackage{graphicx}
\usepackage{multirow}
\numberwithin{equation}{section}
\newtheorem{thm}{Theorem}[section]
\newtheorem{prop}[thm]{Proposition}
\newtheorem{lem}[thm]{Lemma}
\newtheorem{cor}[thm]{Corollary}

\usepackage{bm}
\usepackage{amssymb,amsmath,mathrsfs,txfonts}
\newtheorem{remark}{\indent Remark}[section]
\usepackage[colorlinks=true]{hyperref}
\hypersetup{urlcolor=blue, citecolor=red}

\title{\bf Asymptotic stability of solutions to a hyperbolic-elliptic coupled system of the radiating gas on the half line}

\author[1]{Shanming Ji}
\author[2]{Minyi Zhang}
\author[3]{Changjiang Zhu\thanks{Corresponding author. Email: machjzhu@scut.edu.cn}}
\affil[1,2,3]{School of Mathematics, South China University of Technology, Guangzhou 510641, P.R. China}

\date{} 

\begin{document}
\maketitle

\begin{abstract}
  \indent This paper is concerned with the asymptotic stability of the solution to an initial-boundary value problem on the half line for a hyperbolic-elliptic coupled system of the radiating gas, where the data on the boundary and at the far field state are defined as $u_-$ and $u_+$ satisfying $u_-<u_+$.
  For the scalar viscous conservation law case, it is known by the work of Liu, Matsumura, and Nishihara (SIAM J. Math. Anal. {\bf 29} (1998) 293-308) that the solution tends toward rarefaction wave or stationary solution or superposition of these two kind of waves depending on the distribution of $u_\pm$. Motivated by their work, we prove the stability of the above three types of wave patterns for the hyperbolic-elliptic coupled system of the radiating gas with small perturbation.
  A singular phase plane analysis method is introduced to show the existence and the precise asymptotic behavior of the stationary solution, especially for the degenerate case: $u_-<u_+=0$ such that the system has inevitable singularities.
  The stability of rarefaction wave, stationary solution, and their superposition, is proved by applying the standard $L^2$-energy method.
\end{abstract}

\vspace{3mm}


\vspace{3mm}


\hypersetup{%
  linkcolor  = black
}

\tableofcontents

\hypersetup{%
  linkcolor  = red
}
\section{Introduction}
We consider the asymptotic behavior of solutions to a model of hyperbolic-elliptic coupled system of radiating gas on the half line.
The initial-boundary value problem (IBVP) of the above model on the half-line $\mathbbm{R}_+=(0,+\infty)$ is
\begin{equation}\label{eq1.1}
\begin{cases}
\displaystyle
u_t+f(u)_x+q_x=0, \ \ \ \ x\in \mathbbm{R}_+, \ \ \ t>0, \\
-q_{xx}+q+u_x=0, \ \ \ \  x\in \mathbbm{R}_+, \ \ \ t>0,
\end{cases}
\end{equation}
with boundary condition
\begin{equation}\label{eq1.2}
 u(0,t)=u_-, \quad t>0,
\end{equation}
and initial data
\begin{equation}\label{eq1.3}
u(x,0)=u_0(x)\
\begin{cases}
=u_-, \ \ \ \ x=0,\\
\rightarrow u_+, \ \ \ \ x\rightarrow +\infty,
\end{cases}
\end{equation}
where the flux $f(u)$ is a given smooth function of $u$, $u_\pm$ are given constants, $u$ and $q$ are unknown functions of the spacial variable $x\in \mathbbm{R}_+$ and the time variable $t$. Generally, $u$ and $q$ represent the velocity and the heat flux of the gas respectively. Throughout this paper, we impose the following condition:
\begin{equation*}
f(0)=f'(0)=0,~f''(u)>0 \text{ for } u\in \mathbbm{R}, \quad \text{and } u_-<u_+.
\end{equation*}
\indent Such a hyperbolic-elliptic coupled system appears typically in radiation hydrodynamics, cf. \cite{Vincenti, Zel66}.
The simplified model \eqref{eq1.1} was first recovered by Hamer in \cite{Hamer}, and for the derivation of system \eqref{eq1.1},
we refer to \cite{Gao2, Hamer, Vincenti}.
In the \emph{in-flow} case of $u_->0$, the boundary condition \eqref{eq1.2} is necessary for the single hyperbolic equation
\eqref{eq1.1}$_1$, we additionally impose $q(0,t)=0$ for the well-posedness of the coupled elliptic equation \eqref{eq1.1}$_2$.
The problem with additional condition $q(0,t)=0$ for \emph{in-flow} case will still be denoted by \eqref{eq1.1}-\eqref{eq1.3}
for the sake of simplicity.
On the contrary, in the \emph{out-flow} case of $u_-<0$, the boundary condition \eqref{eq1.2} is enough for the
hyperbolic-elliptic coupled system \eqref{eq1.1},
where the single hyperbolic equation \eqref{eq1.1}$_1$ is over-determined with the boundary condition \eqref{eq1.2}
and the single elliptic equation \eqref{eq1.1}$_2$ is under-determined since no information of $q$
is given at the boundary.

The system \eqref{eq1.1} has been extensively studied by several authors in different contexts recently,
but most of them are in the case of the whole space.
For the one-dimensional whole space, concerning the large-time behavior of solutions to the Cauchy problem
\begin{equation}\label{Cauchyproblem}
	\begin{cases}
		u_t+f(u)_x+q_x=0, \quad x\in \mathbbm{R},\ \ t>0,\\
		-q_{xx}+q_x+u_x=0,\quad x\in \mathbbm{R},\ \ t>0,\\
		u(x,0)=u_{0}(x)\to u_{\pm} ,\quad x\to\pm \infty,
	\end{cases}
\end{equation}
Tanaka in \cite{Tanaka} proved the stability of the diffusion wave as a self-similar solution to the viscous Burgers equation
for the special case: $u_-=u_+=0$.
Kawashima and Nishihara in \cite{Kawashima1} discussed the case of $u_->u_+$ and showed that the solution to the Cauchy problem \eqref{Cauchyproblem} approaches the travelling wave of shock profile.
Kawashima and Tanaka in \cite{Kawashima} investigated the remaining case of $u_-<u_+$ and
proved the stability of the rarefaction wave for the Cauchy problem \eqref{Cauchyproblem}
with the asymptotic convergence rate.
Ruan and Zhang in \cite{Ruan1} further studied the case: $u_-<u_+$ for general flux $f(u)$.
The radiating gas system \eqref{Cauchyproblem} in the following scalar equation form with convolution
\begin{equation} \label{eq-scalar}
u_t+f(u)_x+u-\psi*u=0, \quad x\in \mathbbm{R},\ \ t>0,
\end{equation}
where $\psi(x)$ is the fundamental solution to the elliptic operator $-\Delta +I$ in $\mathbbm{R}$,
was studied in \cite{DiFrancescoAML, Duan12, LiuH, LiuKawashima, Nguyen, Serre, Wang}.
It was Schochet and Tadmor in \cite{Schochet} who first proved the $W^{1,\infty}$ regularity of the solution to \eqref{eq-scalar}.
Then Lattanzio and Marcati in \cite{Lattanzio} studied the well-posedness and relaxation limits of weak entropy solutions.
Yang and Zhao in \cite{YangZhao} constructed the Lax-Friedrichs' scheme and obtained the BV estimates.
For the multi-dimensional case,
Gao, Ruan and Zhu investigated the asymptotic rate towards the planar rarefaction waves to the Cauchy problem for a hyperbolic-elliptic coupled system (see \cite{Gao1, Gao2, Ruan10}).
Di Francesco in \cite{DiFrancesco} studied the global well-posedness and the relaxation limits
of the multi-dimensional radiating gas system.
Duan, Fellner and Zhu in \cite{Duan10} studied the stability and optimal time decay rates of planar rarefaction waves
for a radiating gas model based on Fourier energy method.
The structure of shock waves in the radiating gas dynamics was also investigated by many authors,
see \cite{Heaslet, Lowrie08, Lowrie07, Ohnawa, Zel57}.
Moreover, the large time behaviors of the solutions for viscous conservation laws, and other system were studied by many authors, see \cite{Hsiao, HuangFM, Ito, LuoT, MatsumuraNishihara, Nishikawa, Xin, ZhaoHJ, ZhuCJ1}.

We are concerned with the asymptotic behavior of the solution to \eqref{eq1.1}-\eqref{eq1.3}.
For the case of $u_-<u_+$, the corresponding Riemann problem for the inviscid Burgers equation
\begin{equation*}
	\begin{cases}
	     u_t+f(u)_x=0, \quad  x\in \mathbbm{R},\ \ t>0,\\
	     u(x,0)=u_0^R(x):=
	     \begin{cases}
	     	u_-,\ \ x<0,\\
	     	u_+,\ \ x>0,
	     \end{cases}
    \end{cases}
\end{equation*}
admits a simple rarefaction wave solution
\begin{equation*}
    u^R(x/t)=
        \begin{cases}
	         u_-, \ \ \ \ x\le f'(u_-)t,\\[2mm]
                \displaystyle
	        (f')^{-1}\big(\frac{x}{t}\big), \ \ \ \ f'(u_-)t\le x\le f'(u_+)t, \\[2mm]
	         u_+, \ \ \ \ x\ge f'(u_+)t.
        \end{cases}	
\end{equation*}
Additionally, the rarefaction wave $(u,q)$ of the system \eqref{eq1.1} is defined as $(u,q)=(u^R(x/t),-\partial_xu^R(x/t))$.
According to the convex function $f(u)$ and the arguments used by Liu, Matsumura, and Nishihara in \cite{Liu},
we have the following five cases (taking the typical form of $f(u)=\frac{1}{2}u^2$ for example)
due to the signs of the characteristic speeds $u_\pm$:
\begin{equation*}
\begin{aligned}
(1)\ \  u_-<u_+<0 \quad (\text{equivalent to } f'(u_-)<f'(u_+)<0);\\
(2)\ \  u_-<u_+=0 \quad (\text{equivalent to } f'(u_-)<f'(u_+)=0);\\
(3)\ \  u_-<0<u_+ \quad (\text{equivalent to } f'(u_-)<0<f'(u_+));\\
(4)\ \  0=u_-<u_+ \quad (\text{equivalent to } 0=f'(u_-)<f'(u_+));\\
(5)\ \  0<u_-<u_+ \quad (\text{equivalent to } 0<f'(u_-)<f'(u_+)).
\end{aligned}
\end{equation*}
By using energy method, for the cases: $u_-<u_+\le 0$, $ 0\le u_-<u_+$ and $u_-<0<u_+$,
Liu, Matsumura, and Nishihara in \cite{Liu} proved that the initial-boundary value problem on the half line for scalar viscous generalized Burgers equation admits a unique global solution and it converges to the stationary solution, the rarefaction wave and the superposition of the nonlinear waves, respectively, as $t\to \infty$.
Since then, the initial boundary value problem on the half line $\mathbbm{R_+}$ for different models have been studied by many authors, cf. \cite{Hashimoto, Liu1, Liu2, MatsumuraMei, YangT}, and references therein.
In the case $(4)$, the convergence of the initial-boundary value problem to a rarefaction wave
has been investigated by Ruan and Zhu in \cite{Ruan2}.
However, there remain four cases to be considered in the previous studies.
Motivated by the classification by Liu, Matsumura, and Nishihara in \cite{Liu},
here we prove the asymptotic behavior of the solutions to \eqref{eq1.1}-\eqref{eq1.3} for all remaining cases:
$u_-<u_+\le 0$, $ 0<u_-<u_+$, and $u_-<0<u_+$.

The main features of the hyperbolic-elliptic coupled system \eqref{eq1.1}-\eqref{eq1.3} on the half line
are different from the previous study on scalar viscous Burgers equation on the half line in \cite{Liu},
or the Cauchy problem \eqref{Cauchyproblem} (equivalent to the scalar form \eqref{eq-scalar} with convolution),
due to the following reasons:

\textbullet \ \
The hyperbolic-elliptic coupled system \eqref{eq1.1}-\eqref{eq1.3} on the half line cannot be converted to
a scalar equation with convolution,
since the information of the solution $q$ on the boundary ($q(0,t)$ or $q_x(0,t)$) is unknown for the \emph{out-flow} cases (1) and (2).
In fact, the ``boundary condition'' of the elliptic problem \eqref{eq1.1}$_2$ is determined in an inverse problem fashion
such that the hyperbolic equation \eqref{eq1.1}$_1$ satisfies the boundary condition $u(0,t)=u_-$,
which is unnecessary for a single hyperbolic equation with characteristic curves running out of the region.

\textbullet \ \
For the degenerate case (2): $u_-<u_+=0$, there arise inevitable singularities
in the analysis of the existence and spatial decay rates of the stationary solution.
It should be noted that the spatial decay estimates of the stationary solution are essential to
the energy estimates of the perturbation problem.
We employ a singular phase plane analysis method with a series of approximated solutions to show the existence
(see Lemma \ref{le-vk}), and
then we utilize the finite series expansion to derive the precise decay estimates of higher order derivatives (see Lemma \ref{le-expansion}).

\textbullet \ \
The estimates on the boundary terms are more subtle, especially for $w_{xx}(0,t)$
of the perturbation $w:=u-\tilde u-\hat u$ when $0<u_-<u_+$,
since this case corresponds to the \emph{in-flow} problem.
To overcome it, we find out the relation \eqref{xsbwxtzxtbj} at boundary $x=0$ to estimate $w_{xt}(0,t)$,
which plays an important role in estimating $w_{xx}(0,t)$.

Additionally, in order to avoid too much tedious estimations in Sections \ref{sec-4} and \ref{sec-5},
we introduce Lemma \ref{le-Bessel} to simplify the proof of the asymptotic behavior of perturbation $z$ after we get the asymptotic behavior of $w$.

This paper is organized as follows. In Section \ref{sec-2},
we firstly prepare the basic properties of the rarefaction wave and stationary solution.
Secondly, we give some inequalities for the maximum norm to elliptic problem,
which is essential in estimating the asymptotic behavior of stationary solution.
Finally, we present our main theorems.
In Section \ref{sec-3}, we show the asymptotic behavior for the case (5), which correspond to the rarefaction wave.
The cases (1) and (2) corresponding to the stationary solutions are investigated in Section \ref{sec-4}.
In the final Section \ref{sec-5}, referring to the results from \cite{Ruan2},
the combination of the cases (2) and (4) can help us to consider the case (3) of superposition waves.

\vspace{3mm}

\textbf{Notations.}
Hereafter, we denote generic positive constants by $C$ and $c$ unless they need to be distinguished.
For function spaces, $L^p=L^p(\mathbbm{R}_+)$ with $1\le p \le \infty$ denotes the usual Lebesgue space on $\mathbbm{R}_+$ with the norm $|\cdot|_p$.
For a non-negative integer $l$, $H^l=H^l(\mathbbm{R}_+)$ denotes the $l$-th order Sobolev space in the $L^2$-sense, equipped with the norm $\|\cdot\|_l$. We note that $H^0=L^2$ and $\|\cdot\|_0=\|\cdot\|$. For simplicity, $\|f(\cdot,t)\|$ is denoted by $\|f(t)\|$, and $\|f(\cdot,t)\|_l$ is denoted by $\|f(t)\|_l$.

\section{Preliminaries and Main Theorems} \label{sec-2}
Without loss of generality, we may take the typical case of
$f(u)=\frac{1}{2} u^2$ to simplify the calculations in the following,
since the proof for general convex function $f(u)$ can be slightly modified.
\subsection{Construction and Properties of the Smooth Rarefaction Wave}

In this subsection we consider the cases (4): $0=u_-<u_+$ and (5): $0<u_-<u_+$.
Since the rarefaction wave $u^R(x/t)$ is not smooth enough, we construct the smooth approximation $\tilde{u}=\tilde{u}_i(x,t) ~ (i=4,5)$ by employing the ideal of Hattori and Nishihara in \cite{Hattori}.
We define $\tilde{u}_5(x,t)$ as a solution of the Cauchy problem
\begin{equation}\label{xsbghbj} 
  \begin{cases}
    \tilde{u}_t+\tilde{u}\tilde{u}_x=\tilde{u}_{xx}, \ \ \ \ x\in \mathbbm{R}, \ \ \ t>0\\
    \tilde{u}(x,0)=\tilde{u}_0^R(x), \ \ \ \ x\in \mathbbm{R},
  \end{cases}
\end{equation}
with $\tilde{q}=-\tilde{u}_x$ and the initial data $\tilde{u}_0^R(x)$ is defined by
\begin{equation}\label{xsbcsz}
\tilde{u}_0^R(x)=
  \begin{cases}
    u_-, \ \ \ \ x<0,\\
    u_+, \ \ \ \ x>0,
  \end{cases}
\end{equation}
for the case $u_->0$.
When $u_-=0$, $\tilde{u}_4(x,t)$ defined as above does not converge to the corresponding rarefaction wave fast enough
near the boundary $x=0$. Therefore, we need to modify $\tilde{u}_0^R(x)$ as
\begin{equation*}
  \tilde{u}_0^R(x)=
  \begin{cases}
    -u_+, \ \ \ \ x<0,\\
    u_+, \ \ \ \ x>0,
  \end{cases}
\end{equation*}
such that the solution $\tilde{u}_4(x,t)$ of \eqref{xsbghbj} satisfies $\tilde{u}_4(0,t)=0$.
Using the Hopf-Cole transformation, the explicit formula of $\tilde{u}$ can be obtained.
Here, we give some properties of smooth approximation solutions $\tilde{u}$ in Lemma \ref{lemxsb}.
The proof of Lemma \ref{lemxsb} can be found in \cite{Hattori,Kawashima}.

\begin{lem}\label{lemxsb} 
For $1\le p\le \infty$ and $t>0$, $\tilde{u}_i$ $(i=4,5)$ satisfies the following estimates:
\\[3mm] \indent
$\mathrm{(i)}$ $0\le \tilde{u}_5(0,t)-u_-\le C \delta \mathrm{e}^{-c(1+t)}$ for $u_->0$ and $\tilde{u}_4(0,t)=0$ for $u_-=0$;
\\[3mm] \indent
$\mathrm{(ii)}$ $|\partial_x^k\tilde{u}_5(0,t)|\le C \delta \mathrm{e}^{-c(1+t)}$, \ \ \ \ $|\partial_x^k\tilde{u}_4(0,t)|\le C (1+t)^{-\frac{1}{2} (k+1)}$, \ \ \ \ $k=1,2,3,4$;
\\[3mm] \indent
$\mathrm{(iii)}$ $|\tilde{u}_i(t)-\tilde{u}^R(t)|_p\le C (1+t)^{-\frac{1}{2}+\frac{1}{2p}  }$, \ \ \ \ $i=4,5 $;
\\[3mm] \indent
$\mathrm{(iv)}$ $|\tilde{u}_{ix}(t)|_p\le C \delta^{\frac{1}{p} }(1+t)^{-1+\frac{1}{p} }$, \ \ \ \ $|\tilde{u}_{it}(t)|_p\le C \delta^{\frac{1}{p} }(1+t)^{-1+\frac{1}{p} }$, \ \ \ \ $i=4,5 $;
\\[3mm] \indent
$\mathrm{(v)}$ $|\partial_x^k \partial_t^l \tilde{u}_i(t)|_p\le C \delta(1+t)^{-\frac{1}{2}(k+l -\frac{1}{p}) }$, \ \ \ \ $k+ l=1,2,3,4, \ \ \ \ k,l\in \mathbbm{N}, \ \ \ \ i=4,5$;
\\[3mm] \indent
$\mathrm{(vi)}$ $\tilde{u}_{ix}>0, \ \ \ \ u_-\le \tilde{u}_i\le u_+, \ \ \ \ x \in \mathbbm{R}, \ \ \ \ i=4,5$.
\end{lem}

We note that the boundary value $\tilde{u}_5(0,t)\not=u_-$ if $u_->0$.
In this case, the perturbation $u(x,t)-\tilde{u}(x,t)$ has a ``boundary layer'' $u_--\tilde{u}(0,t)$ at $x=0$.
To solve this problem, we need to modify $\tilde{u}$ near the boundary $x=0$.
Referring to the method of Nakamura in \cite{Nakamura}, our modified smooth approximation $(\varphi(x,t),\psi(x,t))$ is defined as
\begin{equation}\label{ghxsbdy} 
  \begin{cases}
    \varphi(x,t):=\tilde{u}(x,t)-\hat{u}(x,t),\\
    \psi(x,t):=\tilde{q}(x,t)-\hat{q}(x,t),
  \end{cases}
\end{equation}
where
\begin{equation}\label{xsbjzhs} 
  \begin{cases}
    \hat{u}(x,t):=(\tilde{u}(0,t)-u_-)\mathrm{e} ^{-x},\\
    \hat{q}(x,t):=-\tilde{u}_x(0,t)\mathrm{e} ^{-x}.
  \end{cases}
\end{equation}
Note that $\hat{u}(0,t)=0$ if $u_-=0$. Substituting \eqref{ghxsbdy} into \eqref{xsbghbj}, we get the equation of $(\varphi(x,t),\psi(x,t))$:
\begin{equation*}
  \begin{cases}
    \varphi_t+\varphi \varphi_x=\varphi_{xx}-\hat{u}_t- \varphi_x \hat{u}- \varphi \hat{u}_x- \hat{u}\hat{u}_x+\hat{u}_{xx}, \\
    \varphi(0,t)=u_-, \ \ \ \ \varphi(+\infty,t)=u_+,\\
    \varphi(x,0)=\varphi_0(x)=\tilde{u}_0(x)-\hat{u}(x,0), \ \ \ \ x\in \mathbbm{R}_+,
  \end{cases}
\end{equation*}
and $\psi(x,t)=-\tilde{u}_x(x,t)+\tilde{u}_x(0,t)\mathrm{e}^{-x}$ which satisfies $\psi(0,t)=0$.

According to Lemma \ref{lemxsb}, by simple calculations, we can conclude the following estimates of $\varphi(x,t)$.

\begin{lem}\label{lemjzxsb} 
For $1\le p\le\infty$ and $t\ge0$, $\varphi(x,t)$ satisfies:
\\[3mm] \indent
$\mathrm{(i)}$ $\varphi_x>0$, \ \ $|\varphi(x,t)|<u_+$, \ \ for $x \in \mathbbm{R}$;
\\[3mm] \indent
$\mathrm{(ii)}$ $|\varphi(t)-u^R(t)|_p\le C(1+t)^{-\frac{1}{2} +\frac{1}{2p} }$;
\\[3mm] \indent
$\mathrm{(iii)}$ $|\varphi_x(t)|_p\le C \delta^{\frac{1}{p} }(1+t)^{-1+\frac{1}{p} }$, \ \ \ \ $|\varphi_t(t)|_p\le C \delta^{\frac{1}{p} }(1+t)^{-1+\frac{1}{p} }$;
\\[3mm] \indent
$\mathrm{(iv)}$ $|\partial_x^k \partial_t^l  \varphi(t)|_p\le C \delta(1+t)^{-\frac{1}{2}(k+ l -\frac{1}{p}) }$, \ \ \ \ $k+l =1,2,3,4$, \ \ \ $k,l\in\mathbbm{N}$;
\\[3mm] \indent
$\mathrm{(v)}$ $|\partial_x^k \partial_t^l R_1(t)|_p\le C \delta \mathrm{e}^{-c(1+t)}$, \ \ \ \ $k,l=0,1,2$;
\\[3mm] \indent
$\mathrm{(vi)}$ $|\partial_x^k \partial_t^l R_2(t)|_p\le C \delta (1+t)^{-\frac{1}{2}(k+l+3-\frac{1}{p}) }$, \ \ \ \ $k,l=0,1$.
\\[2mm]
Here,
$$R_1(x,t)=\hat{u}_t+ \varphi_x \hat{u}+ \varphi \hat{u}_x+ \hat{u}\hat{u}_x+\hat{q}_x,$$
and
$$R_2(x,t)=-\varphi_{xxx}-\hat{u}_{xxx}-\hat{q}_{xx}+\hat{u}_x+\hat{q}.$$
\end{lem}

\subsection{Construction and Properties of the Stationary Solution}
In this subsection, we consider the cases (1): $u_-<u_+<0$ and (2): $u_-<u_+=0$,
and we will show that the IBVP \eqref{eq1.1}-\eqref{eq1.3} has a stationary solution $(\bar{u}(x),\bar{q}(x))$, which satisfies
\begin{equation}\label{eq2.0}
\begin{cases}
 (\frac{1}{2}\bar{u}^2)_x=-\bar{q}_x,\\
  -\bar{q}_{xx}+\bar{q}+\bar{u}_x=0,\\
 \bar{u}(0)=u_-,\\
 \bar{u}(+\infty)=u_+, \ \ \ \ \bar{q}(+\infty)=0.
 \end{cases}
\end{equation}
Integrating the system \eqref{eq2.0}$_1$ with respect to $x$ over $[x,\infty)$, we conclude that
$\frac{1}{2}\bar{u}^2=\frac{1}{2}u_+^2-\bar{q}\ge0.$
Thus, $\bar{u}$ can be expressed as
$$\bar{u}=\pm\sqrt{u_+^2-2\bar{q}}.$$
Under the assumption of $u_-<u_+\le0$, the function $\bar{u}$ connecting $u_-$ and $u_+$ should be chosen as
$$\bar{u}=-\sqrt{u_+^2-2\bar{q}},\quad \text{such~that}\quad \bar{u}_x=\frac{\bar{q}_x}{\sqrt{u_+^2-2\bar{q}}}.$$
Substituting $\bar{u}_x$ into \eqref{eq2.0}$_2$, the system \eqref{eq2.0} is converted to
\begin{equation}\label{eq2.1}
\begin{cases}
 -\bar{q}_{xx}+\bar{q}+\frac{\bar{q}_x}{\sqrt{u_+^2-2\bar{q}}}=0,\\
 \bar{q}(0)=\bar{q}_-<0, \ \ \ \ \bar{q}(+\infty)=0.
 \end{cases}
\end{equation}
Here, $\bar{q}_-$ satisfies $\bar{q}_-=\frac{1}{2}u_+^2-\frac{1}{2}u_-^2<0$.

In the cases (1): $u_-<u_+<0$ and (2): $u_-<u_+=0$, we prove that the initial-boundary value problem \eqref{eq2.0} has a stationary solution $(\bar{u},\bar{q})=(\bar{u}_i(x),\bar{q}_i(x)),i=1,2$, respectively.

\begin{lem}\label{lem2.1.1}
Suppose $u_-<u_+\le 0$ and let $\delta:=|u_--u_+|$. Then there exists a solution $\bar{u}_i,(i=1,2)$ to the stationary problem \eqref{eq2.0}, such that the following estimates hold for some positive constants C and $\lambda$:
\\[3mm] \indent
$\mathrm{(i)}$ $\bar{u}_{ix}(x)>0, \quad i=1,2$;
\\[3mm] \indent
$\mathrm{(ii)}$ $\left|\partial_x^k\left\{\bar{u}_1(x)-u_+\right\}\right|\le C\delta  \mathrm{e}^{-\lambda x},  \quad  k=0,1,2,3$;
\\[3mm] \indent
$\mathrm{(iii)}$ $\displaystyle |\partial_x^k\bar{u}_2(x)|  \le C \frac{\delta^{k+1}}{(1+\delta x)^{k+1}},  \quad  k=0,1,2,3,4$;
\\[3mm] \indent
$\mathrm{(iv)}$ $\displaystyle \left|\frac{\bar{u}_{ixx}^2}{\bar{u}_{ix}}\right|\le C\delta,  \quad  \left|\frac{\bar{u}_{ixxx}^2}{\bar{u}_{ix}}\right|\le C\delta, \ \ \ \ i=1,2$;
\\[3mm] \indent
$\mathrm{(v)}$ $\displaystyle \left|\frac{\bar{u}_{2xxx}}{\bar{u}_{2x}}\right|\le C\delta,  \quad  \left|\frac{\partial_x^4\bar{u}_{2}}{\bar{u}_{2x}}\right|\le C\delta$.
\end{lem}

In the degenerate case (2): $u_-<u_+=0$, the elliptic problem \eqref{eq2.1} has singularity near $\bar q=0$,
and $\bar q(+\infty)=0$,
which means the singularity is inevitable.
The singularity causes significant difficulty in the analysis of the existence and asymptotic behavior of the
stationary solution.
We will present the detailed proof of Lemma \ref{lem2.1.1} for the degenerate case (2)
by applying a generalized singular phase plane analysis method.
Then we sketch the main lines of the proof for case (1).

For the sake of convenience, in this subsection we set
\begin{equation} \label{eq-sq}
s(x)=-\bar q(x), \quad v(x)=\bar q_x(x), \quad x>0,
\end{equation}
then $0<s(x)<s_0:=\frac{1}{2}u_-^2-\frac{1}{2}u_+^2>0$,
$v(x)>0$, are solutions to the following problem
\begin{equation} \label{eq-sv}
\begin{cases}
s_x=-v, \\
\displaystyle
v_x=-s+\frac{v}{\sqrt{u_+^2+2s}}, \\
0<s(x)<s_0, \quad v(x)>0, \quad x\in(0,+\infty).
\end{cases}
\end{equation}
We first focus on the degenerate case of $u_+=0$, and
the problem has singularity at where $s=0$
\begin{equation} \label{eq-sv-0}
\begin{cases}
s_x=-v, \\
\displaystyle
v_x=-s+\frac{v}{\sqrt{2s}}, \\
0<s(x)<s_0, \quad v(x)>0, \quad x\in(0,+\infty).
\end{cases}
\end{equation}

\begin{lem} \label{le-phase}
For any given $s_0>0$, if $\tilde v(s)$ solves the following singular equation
\begin{equation} \label{eq-v}
\begin{cases}
\displaystyle
\frac{\mathrm{d}\tilde v}{\mathrm{d}s}=\frac{s}{\tilde v}-\frac{1}{\sqrt{2s}}, \quad s\in(0,s_0),\\
\displaystyle
\tilde v(s)>0 \mathrm{~for~}s\in(0,s_0), \quad \lim_{s\rightarrow0^+}\tilde v(s)=0,
\quad \int_0^{s_0}\frac{1}{\tilde v(s)}\mathrm{d}s=+\infty,
\end{cases}
\end{equation}
then the function $s(x)$ defined by
\begin{equation} \label{eq-sx}
x=-\int_{s_0}^{s(x)}\frac{1}{\tilde v(\tau)}\mathrm{d}\tau, \qquad x\in(0,+\infty)
\end{equation}
is a solution of \eqref{eq-sv-0}
with $v(x):=\tilde v(s(x))>0$.
\end{lem}
{\it\bfseries Proof.}
The positivity of $\tilde v(s)$ on $s\in(0,s_0)$ and the singularity at $s=0$
under the conditions in \eqref{eq-v}
imply that the function $s(x)$ is well-defined on $x\in(0,+\infty)$ by \eqref{eq-sx}
such that $s(0)=s_0>0$ and $s(+\infty)=0$, $s(x)$ is a strictly decreasing function
for $x\in(0,+\infty)$.
Differentiating the identity \eqref{eq-sx} with respect to $x$ shows that
$$
1=-\frac{1}{\tilde v(s(x))}s_x(x)=-\frac{s_x(x)}{v(x)}.
$$
Further,
$$
v_x(x)=\frac{\mathrm{d}}{\mathrm{d}x}\tilde v(s(x))=\frac{\mathrm{d}\tilde v}{\mathrm{d}s}\cdot\frac{\mathrm{d}s(x)}{\mathrm{d}x}
=\Big(\frac{s(x)}{\tilde v(s(x))}-\frac{1}{\sqrt{2s(x)}}\Big)\cdot(-v(x))
=-s(x)+\frac{v(x)}{\sqrt{2s(x)}}.
$$
The proof is completed.
$\hfill\Box$

The main feature of the problem \eqref{eq-v} is that
it has two singularities: (i) for $s$ near zero, the function $-\frac{1}{\sqrt{2s}}$
is un-bounded and not Lipschitz continuous;
(ii) for $\tilde v=0$ as the condition $\lim_{s\rightarrow0^+}\tilde v(s)=0$ required,
the term $\frac{s}{\tilde v}$ is also un-bounded and not Lipschitz continuous.
In order to handle these two singularities, we consider an approximated problem
for $k\in \mathbbm N^+$
(without loss of generality we may assume that $\frac{1}{k^2}<s_0$,
otherwise we only consider large $k\in \mathbbm N^+$ with $k\ge k_0$,
where $k_0:=[\frac{1}{\sqrt{s_0}}]+1$)
\begin{equation} \label{eq-vk}
\begin{cases}
\displaystyle
\frac{\mathrm{d}\tilde v_k}{\mathrm{d}s}=\frac{s}{\tilde v_k}-\frac{1}{\sqrt{2s}},
\quad s\in\Big(\frac{1}{k^2},s_0\Big),\\[4mm]
\displaystyle
\tilde v_k(s)>0 \mathrm{~for~}s\in\Big(\frac{1}{k^2},s_0\Big),
\quad \tilde v_k\Big(\frac{1}{k^2}\Big)=\frac{1}{k}.
\end{cases}
\end{equation}

\begin{figure}[htb]
\centering
\includegraphics[width=0.6\textwidth]{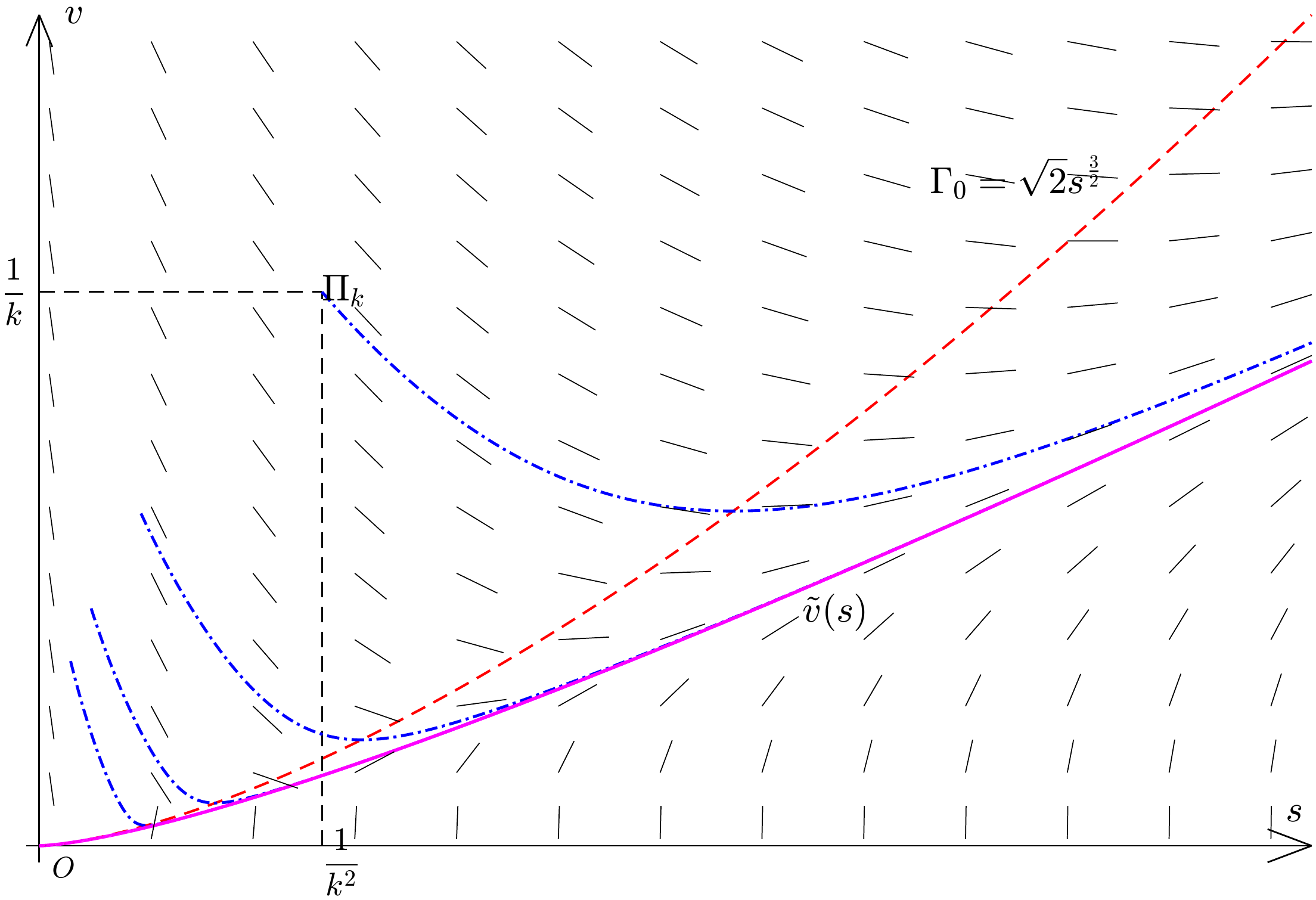}
\caption{The singular phase plane corresponding to system \eqref{eq-sv-0}:
the dash-dot lines correspond to solutions $\tilde v_k(s)$ for different $k$ in the proof of Lemma \ref{le-vk},
the solid line corresponds to the limit function $\tilde v(s)$.}
\label{fig-phase}
\end{figure}

The above approximated problem is solved by the generalized phase plane analysis method.
See the illustration Figure \ref{fig-phase}.

\begin{lem} \label{le-vk}
For any $s_0>0$ and any $\frac{1}{k^2}<s_0$ (i.e. $k\ge k_0$),
the problem \eqref{eq-vk} admits
a solution $\tilde v_k(s)>0$ on $s\in(\frac{1}{k^2},s_0)$
such that
\\[3mm] \indent
$\mathrm{(i)}$ $\tilde v_k(s)$ is monotone decreasing with respect to $k$, i.e.,
$\tilde v_{k_1}(s)>\tilde v_{k_2}(s)$ for any $k_2>k_1\ge k_0$
on their joint interval $(\frac{1}{k_1^2},s_0)\cap(\frac{1}{k_2^2},s_0)=(\frac{1}{k_1^2},s_0)$;
\\[3mm] \indent
$\mathrm{(ii)}$ $\tilde v_k(s)$ has the following upper bound estimate
$$\tilde v_k(s)\le \overline\Gamma_k(s):=
\begin{cases}
\displaystyle
\frac{1}{k}, ~ & s\in\Big(\frac{1}{k^2},(\frac{1}{\sqrt2k})^\frac{2}{3}\Big] \\[3mm]
\displaystyle
\sqrt2s^\frac{3}{2}, ~ & s\in\Big((\frac{1}{\sqrt2k})^\frac{2}{3},s_0\Big)
\end{cases}
=\max\Big\{\frac{1}{k},\sqrt2s^\frac{3}{2}\Big\}, \quad s\in \Big(\frac{1}{k^2},s_0\Big);
$$

$\mathrm{(iii)}$ $\tilde v_k(s)$ has the following uniformly
lower bound estimate for any $\mu\in(0,\sqrt2)$
$$\tilde v_k(s)\ge \Gamma_\mu(s):=
\begin{cases}
\displaystyle
(\sqrt2-\mu)s^\frac{3}{2}, ~ & s\in(0,\delta_\mu] \\[3mm]
\displaystyle
(\sqrt2-\mu)\delta_\mu^\frac{3}{2}, ~ & s\in(\delta_\mu,s_0)
\end{cases}
=\min\left\{(\sqrt2-\mu)s^\frac{3}{2},(\sqrt2-\mu)\delta_\mu^\frac{3}{2}\right\},
$$
where $\delta_\mu:=\frac{\mu}{3\sqrt2}$.
\end{lem}
{\it\bfseries Proof.}
In the phase plane of \eqref{eq-sv-0}, define an auxiliary function
$$
\Gamma_0(s):=\sqrt2s^\frac{3}{2}, \quad s>0.
$$
The system \eqref{eq-sv-0} is locally uniquely solvable at any point
$(s,v)\in\mathbbm R^+\times\mathbbm R^+$ and the vector field is denoted by
$$ \Phi(s,v):=-v, \qquad \Psi(s,v):=-s+\frac{v}{\sqrt{2s}}.$$
There is no stationary point and no closed periodic orbit
(since $\Phi(s,v)<0$ or according to Bendixson's criterion such that
$\mathrm{div}_{(s,v)}(\Phi,\Psi)=\frac{1}{\sqrt{2s}}>0$)
within the first quadrant,
then the Poinc\'are-Bendixson Theorem
implies that any trajectory must runs to the boundary of the first quadrant
for $x\rightarrow\pm\infty$.
According to the smoothness and the absence of stationary point
of the vector field $(\Phi,\Psi)$ within the first quadrant,
we know that any two trajectories can not intersect with each other at
any point $(s,v)\in\mathbbm R^+\times\mathbbm R^+$.
The graph of $\Gamma_0$ (also denoted by $\Gamma_0$)
divides the first quadrant into two parts
$$
G_1:=\{(s,v)\in\mathbbm R^+\times\mathbbm R^+;v>\Gamma_0(s)\},
\quad
G_2:=\{(s,v)\in\mathbbm R^+\times\mathbbm R^+;v<\Gamma_0(s)\}.
$$
For any point $(s,v)\in G_1$, locally $\Phi<0$ and $\Psi>0$,
which means the trajectory runs through $(s,v)$ to the left-up direction
if the autonomous independent variable $x$ grows;
similarly, for any point $(s,v)\in G_2$, locally $\Phi<0$ and $\Psi<0$.

For any $s_0>0$ and any $\frac{1}{k^2}<s_0$ (i.e. $k\ge k_0$),
consider the dynamic system \eqref{eq-sv-0}
with initial condition $s_k(x_k)=\frac{1}{k^2}$ and $v_k(x_k)=\frac{1}{k}$,
where $x_k\in(0,+\infty)$ is a constant to be determined.
The trajectory corresponding to this local solution is denoted by
$\Pi_k:=\{(s_k(x),v_k(x))\}$.
Since the system \eqref{eq-sv-0} is autonomous,
we will shift $x_k$ to a suitable position in the following proof.
Within the first quadrant, $s_x(x)=\Phi(s,v)=-v<0$, which means that
$s(x)$ is strictly decreasing with respect to $x$.
We can take $x$ as an inverse function of $s$ in the range of $s(x)$ and
then regard $v(x)$ as a function of $s$, this is the local solution
$\tilde v_k(s)$ of \eqref{eq-vk}.
The choice of $x_k$ has no influence on the function $\tilde v_k(s)$.

We only consider the trajectory $\Pi_k$ in the negative $x$ direction,
that is, we consider the solution $(s_k(x),v_k(x))$ such that $x<x_k$  and $x$ is decreasing.
Noticing that $(\frac{1}{k^2},\frac{1}{k})\in G_1$ (for $k\ge2$ without loss of generality),
we see that $\Pi_k$ runs through $(\frac{1}{k^2},\frac{1}{k})$ in the right-down direction
until it reaches some point at $\Gamma_0$.
This must happen since $\Gamma_0(s)$ is increasing and
$\Gamma_0((\frac{1}{\sqrt2k})^\frac{2}{3})=\frac{1}{k}$.
Therefore, there exists a $\hat x_k<x_k$ such that
$s_k(\hat x_k)=:\hat s_k\in(\frac{1}{k^2},(\frac{1}{\sqrt2k})^\frac{2}{3})$
and $v_k(\hat x_k)=:\hat v_k$, satisfying
\begin{equation} \label{eq-lower-vk}
\hat v_k=\sqrt2\hat s_k^\frac{3}{2}\in\Big(
\Gamma_0(\frac{1}{k^2}),\Gamma_0((\frac{1}{\sqrt2k})^\frac{2}{3})\Big)
=\Big(\frac{\sqrt2}{k^3},\frac{1}{k}\Big).
\end{equation}
Locally at the point $(\hat s_k,\hat v_k)\in \Gamma_0$, $\Phi<0$ and $\Psi=0$.
Noticing that $\Gamma_0'(\hat s_k)=\frac{3\sqrt2}{2}\sqrt{\hat s_k}>0$,
we find that $\Pi_k$ runs into $G_2$ region as $x$ decreasing from $\hat x_k$.

We assert that $\Pi_k$ is under $\Gamma_0$ for $x<\hat x_k$
(i.e., $\tilde v_k(s)<\Gamma_0(s)$ for $s>\hat s_k$)
and $\tilde v_k(s)$ is increasing for $s>\hat s_k$.
We prove by contradiction and assume that there exists a $s^*>\hat s_k$
such that $\tilde v_k(s^*)=\Gamma_0(s^*)$ and
$\tilde v_k(s)<\Gamma_0(s)$ for $s\in(\hat s_k,s^*)$,
which means there exists $x^*<\hat x_k$ such that $s_k(x^*)=s^*$ and
$v_k(x^*)=\tilde v_k(s^*)=\Gamma_0(s^*)=:v^*$.
Then $\tilde v_k'(s^*)\ge\Gamma_0'(s^*)$,
and $\Gamma_0'(s^*)=\frac{3\sqrt2}{2}\sqrt{s^*}$,
but
$$
\tilde v_k'(s^*)=\left[\frac{\frac{\mathrm{d}v_k(x)}{\mathrm{d}x}}{\frac{\mathrm{d}s_k(x)}{\mathrm{d}x}}\right]_{x=x^*}
=\frac{-s_k+\frac{v_k}{\sqrt{2s_k}}}{-v_k}\Big|_{x=x^*}
=\frac{s^*}{v^*}-\frac{1}{\sqrt{2s^*}}
=0
$$
at this point, which is a contradiction.
Now that we have proved $\tilde v_k(s)<\Gamma_0(s)$ for $s>\hat s_k$,
or equivalently, $(s_k(x),v_k(x))\in G_2$ for $x<\hat x_k$.
Therefore,
\begin{equation} \label{eq-increasing-vk}
\tilde v_k'(s)=\frac{s}{\tilde v_k}-\frac{1}{\sqrt{2s}}>0, \qquad (s,\tilde v_k(s))\in G_2,
\end{equation}
which shows that $\tilde v_k(s)$ is increasing for $s>\hat s_k$.
To summarize, we proved that
$\tilde v_k(s)$ is decreasing for $s\in(\frac{1}{k^2},\hat s_k)$
and increasing for $s>\hat s_k$,
which means $\tilde v_k(s)\ge \hat v_k$
with $\hat v_k>\frac{\sqrt2}{k^3}$ satisfying \eqref{eq-lower-vk}.
Furthermore, we see that
$$
s_k'(x)=\Phi(s_k,v_k)=-v_k(x)
=-\tilde v_k(s_k(x))
\le-\hat v_k<-\frac{\sqrt2}{k^3},
$$
which implies that
$$
-\frac{k^3}{\sqrt 2}s_0\le-\int_{\frac{1}{k^2}}^{s_0}\frac{1}{\tilde v_k(s)}\mathrm{d}s<0
$$
is finite.
We would shift $x_k=\int_{\frac{1}{k^2}}^{s_0}\frac{1}{\tilde v_k(s)}\mathrm{d}s$
such that $s_k(0)=s_0$.

The above arguments imply that $\tilde v_k(s)\le \overline\Gamma_k(s)$
for $s\in(\frac{1}{k^2},s_0)$.
Next, we show the monotone dependence of $\tilde v_k(s)$ with respect to $k$.
For any $k_2>k_1\ge k_0$,
the trajectory $\Pi_{k_2}$ runs through $(\frac{1}{k_2^2},\frac{1}{k_2})$
in the right-down direction (for $x$ decreasing) until
it reaches some point $(\frac{1}{k_1^2},v_*)\in G_1$ with $v_*<\frac{1}{k_2}<\frac{1}{k_1}$
or some point $(s_*,v_*)$ on $\Gamma_0$ with $s_*\in(\frac{1}{k_2^2},\frac{1}{k_1^2}]$.
In the latter case, $\Pi_{k_2}$ runs across $\Gamma_0$ into $G_2$,
and as shown by the above arguments
(i.e., \eqref{eq-lower-vk} and \eqref{eq-increasing-vk}),
$\tilde v_{k_2}(s)<\Gamma_0(s)$ and $\tilde v_{k_2}(s)$ is increasing for $s>s_*$.
Therefore, $\tilde v_{k_2}(\frac{1}{k_1^2})\in(v_*,\Gamma_0(\frac{1}{k_1^2}))%
\subset(v_*,\frac{1}{k_1})$.
In all cases, $\tilde v_{k_2}(\frac{1}{k_1^2})<\frac{1}{k_1}=\tilde v_{k_1}(\frac{1}{k_1^2})$.
The comparison $\tilde v_{k_2}(s)<\tilde v_{k_1}(s)$ for $s>\frac{1}{k_1^2}$
follows from the fact that any two trajectories cannot intersect with each other
in the first quadrant.

Lastly we show the uniformly lower bound of $\tilde v_k(s)$.
For any $\mu\in(0,\sqrt2)$ and $\delta_\mu:=\frac{\mu}{3\sqrt2}$,
we consider the special curve $\Gamma_\mu$ defined by $v=\Gamma_\mu(s)$.
For any point $(\check s,\check v)\in \Gamma_\mu$ with $\check s\in(0,\delta_\mu)$
and $\check v=\Gamma_\mu(\check s)$,
the direction of vector field
\begin{equation} \label{eq-PsiPhi}
\frac{\Psi(\check s,\check v)}{\Phi(\check s,\check v)}
=\frac{\check s}{\check v}-\frac{1}{\sqrt{2\check s}}
=\frac{\check s}{(\sqrt2-\mu)\check s^\frac{3}{2}}-\frac{1}{\sqrt{2\check s}}
=\frac{\mu}{(\sqrt2-\mu)\sqrt2\sqrt{\check s}}
>\frac{\mu}{2\sqrt{\check s}},
\end{equation}
and the derivative of the curve $\Gamma_\mu$
$$
\Gamma_\mu'(\check s)=\frac{3}{2}(\sqrt2-\mu)\sqrt{\check s}
<\frac{3}{2}\sqrt2\sqrt{\check s}
<\frac{\mu}{2\sqrt{\check s}}<
\frac{\Psi(\check s,\check v)}{\Phi(\check s,\check v)}
$$
for $\check s<\delta_\mu=\frac{\mu}{3\sqrt2}$.
For any point $(\check s,\check v)\in \Gamma_\mu$ with $\check s\ge\delta_\mu$
and $\check v=\Gamma_\mu(\check s)=(\sqrt2-\mu)\delta_\mu^\frac{3}{2}$,
the direction of vector field
\begin{equation} \label{eq-PsiPhi2}
\frac{\Psi(\check s,\check v)}{\Phi(\check s,\check v)}
=\frac{\check s}{\check v}-\frac{1}{\sqrt{2\check s}}
>0,
\end{equation}
since $(\check s,\check v)\in G_2$ as $\Gamma_\mu(s)<\Gamma_0(s)$,
and the curve $\Gamma_\mu$ is horizontal.
Noticing that $\Phi(s,v)=-v<0$, we see that any trajectory $(s(x),v(x))$
runs rightwards as $x$ decreasing.
It follows from \eqref{eq-PsiPhi} and \eqref{eq-PsiPhi2}
that any trajectory starting from a point $(s,v)$ above $\Gamma_\mu$
cannot run through $\Gamma_\mu$ as the independent variable $x$ decreasing.
Therefore, $\tilde v_k(s)>\Gamma_\mu(s)$ for $s\in(\frac{1}{k^2},s_0)$
since $(\frac{1}{k^2},\frac{1}{k})\in G_1$ is above $\Gamma_\mu$ for $k\ge2$.
The proof is completed.
$\hfill\Box$

The solutions $\{\tilde v_k(s)\}$
to the above approximated problem \eqref{eq-vk} are not defined for all
$s>0$, only on $(\frac{1}{k^2},s_0)$.
We define
$$
\tilde v_k^*(s):=
\begin{cases}
\tilde v_k(s), \ \ \ \ s\in\left.\left(\frac{1}{k^2},s_0\right.\right],\\[4mm]
\frac{1}{k}, \ \ \ \ s\in\left.\left(0,\frac{1}{k^2}\right.\right],
\end{cases}
$$
and
\begin{equation} \label{eq-tildev}
\tilde v(s)=\lim_{k\rightarrow\infty}\tilde v_k^*(s), \ \ \ \  s\in(0,s_0).
\end{equation}

\begin{lem} \label{le-tildev}
The function $\tilde v(s)$ is well-defined in \eqref{eq-tildev} for $s\in(0,s_0)$
and $\tilde v(s)$ is a solution to the singular problem \eqref{eq-v}.
Moreover, for any $\mu\in(0,\sqrt2)$, let $\Gamma_\mu$ be the function defined
in Lemma \ref{le-vk}, then
$\Gamma_\mu(s)\le \tilde v(s)\le \sqrt2s^\frac{3}{2}$
for all $s\in(0,s_0)$,
which means $\tilde v(s)=\sqrt2s^\frac{3}{2}+o(s^\frac{3}{2})$
as $s\rightarrow0^+$.
\end{lem}
{\it\bfseries Proof.}
For any fixed $s\in(0,s_0)$, $\tilde v_k(s)$ is defined on $(\frac{1}{k^2},s_0)$,
which contains $s$ if $\frac{1}{k^2}<s$,
i.e., $k\ge [\frac{1}{\sqrt s}]+1$.
According to Lemma \ref{le-vk}, $\{\tilde v_k(s)\}$ is monotone decreasing
with respect to $k$ and is bounded.
Meanwhile, $\{\tilde v_k^*(s)\}$ is monotone decreasing
with respect to $k$ and is bounded on $(0,s_0]$.
Therefore, the limit $\tilde v(s)=\lim_{k\rightarrow\infty}\tilde v_k^*(s)$ exists,
and satisfies $\Gamma_\mu(s)\le \tilde v(s)\le \sqrt2s^\frac{3}{2}$.
It is easy to check that
$\tilde v(s)>0$ for $s\in(0,s_0)$, $\lim_{s\rightarrow0^+}\tilde v(s)=0$, and
$\int_0^{s_0}\frac{1}{\tilde v(s)}\mathrm{d}s=+\infty$.

We show that $\tilde v(s)$ satisfies the differential equation \eqref{eq-v}.
Locally in a neighbourhood of any $s_1\in(0,s_0)$, say $I:=(\frac{s_1}{2},s_0)$,
we rewrite the differential equation \eqref{eq-vk}
(for large $k$ such that $\frac{1}{k^2}<\frac{s_1}{2}$)
as
$$
\frac{1}{2}\frac{\mathrm{d}}{\mathrm{d}s}(\tilde v_k)^2=s-\frac{\tilde v_k}{\sqrt{2s}}, \quad s\in I.
$$
Integrating from $s_1$ shows
\begin{equation} \label{eq-tildevk}
\frac{1}{2}(\tilde v_k(s))^2-\frac{1}{2}(\tilde v_k(s_1))^2
=\frac{1}{2}s^2-\frac{1}{2}s_1^2-\int_{s_1}^s\frac{\tilde v_k(\tau)}{\sqrt{2\tau}}\mathrm{d}\tau,
\quad s\in I.
\end{equation}
Since $\tilde v_k(s)$ is monotone decreasing with respect to $k$ and bounded,
Lebesgue's Dominated Convergence Theorem (or Levi's Theorem) implies that
\begin{equation} \label{eq-tildev0}
\frac{1}{2}(\tilde v(s))^2-\frac{1}{2}(\tilde v(s_1))^2
=\frac{1}{2}s^2-\frac{1}{2}s_1^2-\int_{s_1}^s\frac{\tilde v(\tau)}{\sqrt{2\tau}}\mathrm{d}\tau,
\quad s\in I.
\end{equation}
Differentiating \eqref{eq-tildev0} with respect to $s$ near $s_1$, we have
$$
\tilde v(s)\cdot\frac{\mathrm{d}\tilde v(s)}{\mathrm{d}s}=s-\frac{\tilde v(s)}{\sqrt{2s}},
\quad s\in I.
$$
Therefore, $\tilde v(s)$ is a solution to the problem \eqref{eq-v}.
$\hfill\Box$

\begin{lem}\label{le-sx}
For $s_0=\frac{1}{2}u_-^2-\frac{1}{2}u_+^2\in(0,\frac{1}{6})$,
let $\tilde v(s)$ for $s\in(0,s_0)$ be the solution to the problem \eqref{eq-v}
as defined in \eqref{eq-tildev},
and let $s(x)$ and $v(x)$ be the functions defined by \eqref{eq-sx}
in Lemma \ref{le-phase}.
Then $\bar q(x)=-s(x)$ is a stationary solution to the problem \eqref{eq2.1}
for the degenerate case of $u_+=0$, and has the following decay estimates
\begin{equation} \label{eq-qdecay-1}
-\frac{1}{\Big(\frac{1}{\sqrt{s_0}}+\frac{x}{2\sqrt2}\Big)^2}
\le \bar q(x)\le
-\frac{1}{\Big(\frac{1}{\sqrt{s_0}}+\frac{x}{\sqrt2}\Big)^2},
\quad x\in(0,+\infty),
\end{equation}
and
\begin{equation} \label{eq-qdecay-2}
\frac{\frac{\sqrt2}{2}}{\Big(\frac{1}{\sqrt{s_0}}+\frac{x}{\sqrt2}\Big)^3}
\le \bar q_x(x)\le
\frac{\sqrt2}{\Big(\frac{1}{\sqrt{s_0}}+\frac{x}{2\sqrt2}\Big)^3},
\quad x\in(0,+\infty).
\end{equation}
\end{lem}
{\it\bfseries Proof.}
According to Lemma \ref{le-tildev},
$\tilde v(s)=\sqrt2s^\frac{3}{2}+o(s^\frac{3}{2})$
as $s\rightarrow0^+$,
and $\Gamma_\mu(s)\le \tilde v(s)\le \sqrt2s^\frac{3}{2}$
for all $s\in(0,s_0)$ and $\mu=\frac{\sqrt2}{2}$,
then the function $s(x)$ defined by \eqref{eq-sx} satisfies
(note that $0<s(x)<s_0$)
\begin{equation} \label{eq-zsx1}
x=-\int_{s_0}^{s(x)}\frac{1}{\tilde v(\tau)}\mathrm{d}\tau
\ge -\int_{s_0}^{s(x)}\frac{1}{\sqrt2\tau^\frac{3}{2}}\mathrm{d}\tau
=\sqrt2\Big(\frac{1}{\sqrt{s(x)}}-\frac{1}{\sqrt{s_0}}\Big), \ \ \ \ x\in(0,+\infty),
\end{equation}
and on the other hand,
\begin{align} \nonumber
x=&-\int_{s_0}^{s(x)}\frac{1}{\tilde v(\tau)}\mathrm{d}\tau
\le -\int_{s_0}^{s(x)}\frac{1}{\Gamma_{\frac{\sqrt2}{2}}(\tau)}\mathrm{d}\tau
=-\int_{s_0}^{s(x)}\frac{2}{\sqrt2\tau^\frac{3}{2}}\mathrm{d}\tau
\\[2mm] \label{eq-zsx2}
=&2\sqrt2\Big(\frac{1}{\sqrt{s(x)}}-\frac{1}{\sqrt{s_0}}\Big),
\qquad x\in(0,+\infty),
\end{align}
since $\Gamma_{\frac{\sqrt2}{2}}(s)=(\sqrt2-\frac{\sqrt2}{2})s^\frac{3}{2}$
for $s\in(0,\delta_{\frac{\sqrt2}{2}})$ and
$\delta_{\frac{\sqrt2}{2}}=\frac{\sqrt2/2}{3\sqrt2}=\frac{1}{6}>s_0$.
Combining the above estimates \eqref{eq-zsx1} and \eqref{eq-zsx2} implies
\begin{equation} \label{eq-sxdecay}
\frac{1}{\Big(\frac{1}{\sqrt{s_0}}+\frac{x}{\sqrt2}\Big)^2}
\le s(x)\le \frac{1}{\Big(\frac{1}{\sqrt{s_0}}+\frac{x}{2\sqrt2}\Big)^2},
\quad x\in(0,+\infty).
\end{equation}
Since $s(x)$ is a solution to the problem \eqref{eq-sv-0},
we have
$$
s_x(x)=-v(x)=-\tilde v(s(x))
\in [-\sqrt2s^\frac{3}{2}(x),-\Gamma_\frac{\sqrt2}{2}(s(x))], \quad x\in(0,+\infty).
$$
That is,
$$
s_x(x)\le -\frac{\sqrt2}{2}s^\frac{3}{2}(x)
\le -\frac{\frac{\sqrt2}{2}}{\Big(\frac{1}{\sqrt{s_0}}+\frac{x}{\sqrt2}\Big)^3},
$$
and
$$
s_x(x)\ge -\sqrt2s^\frac{3}{2}(x)\ge -
\frac{\sqrt2}{\Big(\frac{1}{\sqrt{s_0}}+\frac{x}{2\sqrt2}\Big)^3},
$$
for $x\in(0,+\infty)$.
$\hfill\Box$

\begin{remark}
The restriction of $s_0<\frac{1}{6}$ is not essential for the existence and
the decay estimates of the stationary solution in Lemma \ref{le-sx}.
For general $s_0\ge\frac{1}{6}$, we can modify the estimate \eqref{eq-zsx2}
such that
$$
-\int_{s_0}^{s(x)}\frac{1}{\Gamma_{\frac{\sqrt2}{2}}(\tau)}\mathrm{d}\tau
=\begin{cases}
\displaystyle
-\int_{\frac{1}{6}}^{s(x)}\frac{2}{\sqrt2\tau^\frac{3}{2}}\mathrm{d}\tau
-\int_{s_0}^{\frac{1}{6}}\frac{1}{\Gamma_{\frac{\sqrt2}{2}}(\frac{1}{6})}\mathrm{d}\tau,
~& \text{if~} s(x)<\frac{1}{6},
\\[6mm]
\displaystyle
-\int_{s_0}^{s(x)}\frac{1}{\Gamma_{\frac{\sqrt2}{2}}(\frac{1}{6})}\mathrm{d}\tau,
~& \text{if~} s(x)>\frac{1}{6}.
\end{cases}
$$
The decay estimates follow similarly.
Here we take $s_0<\frac{1}{6}$ for the sake of the simplicity of the expressions.
This remark is valid for all the estimates in this subsection,
hence we only present the precise estimates
for small $s_0$.
\end{remark}

The asymptotic behavior $\tilde v(s)=\sqrt2s^\frac{3}{2}+o(s^\frac{3}{2})$
as $s\rightarrow0^+$ implies the decay estimates of $\bar q(x)$ and $\bar q_x(x)$.
In order to derive decay estimates of higher order derivatives,
we expand $\tilde v(s)$ to higher order.
Define sequences $\{c_k\}$ and $\{a_k\}$ as following
\begin{equation} \label{eq-ck}
c_1=1, \quad c_k=\sum_{i+j=k}(2j+1)c_i\cdot c_j, \quad \forall k\ge2,
\end{equation}
and
\begin{equation} \label{eq-ak}
a_k=(-1)^{k+1}\sqrt2c_k, \quad \forall k\ge1.
\end{equation}
For example, $a_1=\sqrt2$, $a_2=-3\sqrt2$, $a_3=24\sqrt2$, $a_4=-285\sqrt2$,
$a_5=4284\sqrt2$.
The formal series
$$
\sum_{i=1}^{\infty} a_is^\frac{2i+1}{2}
$$
is not convergent at any point $s>0$,
and then the infinite series expansion method cannot be applied.
However, the finite series expansion still gives the local behavior of the solution $\tilde v(s)$,
which leads to the precise decay estimates of the higher order derivatives of $\bar q(x)$ and $\bar u(x)$.

\begin{lem} \label{le-expansion}
For any $k\in\mathbbm N^+$, there holds
$$
\tilde v(s)=\sum_{i=1}^ka_is^\frac{2i+1}{2}+o(s^\frac{2k+1}{2}),
\quad s\rightarrow0^+.
$$
Specifically, for odd $k\in\mathbbm N^+$, there exist $\delta_k>0$ and $M_k>0$ such that
$$
\sum_{i=1}^ka_is^\frac{2i+1}{2}-M_ks^\frac{2k+3}{2}
\le \tilde v(s)\le \sum_{i=1}^ka_is^\frac{2i+1}{2},
\quad s\in(0,\delta_k);
$$
while for even $k\in\mathbbm N^+$, there exist $\delta_k>0$ and $M_k>0$ such that
$$
\sum_{i=1}^ka_is^\frac{2i+1}{2}
\le \tilde v(s)\le
\sum_{i=1}^ka_is^\frac{2i+1}{2}+M_ks^\frac{2k+3}{2},
\quad s\in(0,\delta_k).
$$
\end{lem}
{\it\bfseries Proof.}
We prove the cases $k=2$ and $k=3$.
Other cases follow similarly.
For $k=2$, let $\beta(s):=\sqrt2s^\frac{3}{2}-\tilde v(s)$,
i.e., $\tilde v(s)=\sqrt2s^\frac{3}{2}-\beta(s)$.
The function $\tilde v(s)=\sqrt2s^\frac{3}{2}+o(s^\frac{3}{2})$
as $s\rightarrow0^+$ satisfies the differential equation \eqref{eq-v},
and then $\beta(s)$ satisfies
$$
\displaystyle
\frac{\mathrm{d}\tilde v}{\mathrm{d}s}=
\frac{3}{2}\sqrt{2s}-\frac{\mathrm{d}\beta(s)}{\mathrm{d}s}
=\frac{s}{\tilde v}-\frac{1}{\sqrt{2s}}
=\frac{s}{\sqrt2s^\frac{3}{2}-\beta(s)}-\frac{1}{\sqrt{2s}}.
$$
That is,
\begin{equation} \label{eq-beta}
\frac{\mathrm{d}\beta(s)}{\mathrm{d}s}
=\frac{3}{2}\sqrt{2s}+\frac{1}{\sqrt{2s}}-\frac{s}{\sqrt2s^\frac{3}{2}-\beta(s)}
=\frac{3\sqrt2s^\frac{5}{2}-(1+3s)\cdot\beta(s)}{2s^2-\sqrt{2s}\cdot\beta(s)}.
\end{equation}
We analyze the phase plane $(s,\beta)$
corresponding to the singular differential equation
\eqref{eq-beta} in a similar way as we solve the problem \eqref{eq-v}.

Consider two special curves in the phase plane of $(s,\beta)$
\begin{equation} \label{eq-zcurve}
\Gamma^*(s):=3\sqrt2 s^\frac{5}{2}, \quad
\Gamma_*(s):=3\sqrt2 s^\frac{5}{2}-M_2s^\frac{7}{2},
\quad s\in(0,\delta_2),
\end{equation}
where $M_2>0$ and $\delta_2>0$ are constants to be determined.

We use the same symbol $\Gamma^*$ (or $\Gamma_*$) to denote
the curve as well as the function.
At any point $(\hat s,\hat \beta)\in\Gamma^*$
(i.e., $\hat \beta=3\sqrt2 \hat s^\frac{5}{2}$),
there holds
\begin{equation} \label{eq-zGamma1}
\Gamma^*{}'(\hat s)=\frac{15}{2}\sqrt2\hat s^\frac{3}{2}
>0>\frac{3\sqrt2\hat s^\frac{5}{2}-(1+3\hat s)\cdot\hat \beta}%
{2\hat s^2-\sqrt{2\hat s}\cdot\hat \beta}=\frac{\mathrm{d}\beta(s)}{\mathrm{d}s}%
\Big|_{s=\hat s,\beta=\hat \beta}.
\end{equation}
Meanwhile, at any point $(\hat s,\hat \beta)\in\Gamma_*$
(i.e., $\hat \beta=3\sqrt2 \hat s^\frac{5}{2}-M_2\hat s^\frac{7}{2}$),
we have
\begin{equation} \label{eq-zGamma2}
\Gamma_*'(\hat s)=\frac{15}{2}\sqrt2\hat s^\frac{3}{2}-\frac{7}{2}M_2\hat s^\frac{5}{2}
<\frac{3\sqrt2\hat s^\frac{5}{2}-(1+3\hat s)\cdot\hat \beta}%
{2\hat s^2-\sqrt{2\hat s}\cdot\hat \beta}=\frac{\mathrm{d}\beta(s)}{\mathrm{d}s}%
\Big|_{s=\hat s,\beta=\hat \beta},
\end{equation}
which is equivalent to
$$
\frac{15}{2}\sqrt2\hat s^\frac{3}{2}-\frac{7}{2}M_2\hat s^\frac{5}{2}
<\frac{(M_2-9\sqrt2)\hat s^\frac{7}{2}+3M_2\hat s^\frac{9}{2}}%
{2\hat s^2-6\hat s^3+M_2\sqrt2\hat s^4}.
$$
It suffices to take $M_2=24\sqrt2$ and $\delta_2=\frac{1}{8}$.
The analysis of the trajectories according to
inequalities \eqref{eq-zGamma1} and \eqref{eq-zGamma2} show that
$$
3\sqrt2 s^\frac{5}{2}-M_2s^\frac{7}{2}\le \beta(s)
\le 3\sqrt2 s^\frac{5}{2}, \quad s\in(0,\delta_2).
$$
Therefore,
\begin{equation} \label{eq-tildev2}
\tilde v(s)=\sqrt2s^\frac{3}{2}-\beta(s)~
\begin{cases}
\ge \sqrt2s^\frac{3}{2}-3\sqrt2 s^\frac{5}{2}, \\[2mm]
\le \sqrt2s^\frac{3}{2}-3\sqrt2 s^\frac{5}{2}+M_2s^\frac{7}{2},
\end{cases}
\qquad s\in(0,\delta_2).
\end{equation}

Next, we consider the case of $k=3$.
Let $\alpha(s):=\tilde v(s)-\sqrt2s^\frac{3}{2}+3\sqrt2s^\frac{5}{2}$,
i.e., $\tilde v(s)=\sqrt2s^\frac{3}{2}-3\sqrt2s^\frac{5}{2}+\alpha(s)$.
According to the differential equation \eqref{eq-v} of $\tilde v(s)$,
we see that $\alpha(s)$ satisfies
$$
\displaystyle
\frac{\mathrm{d}\tilde v}{\mathrm{d}s}=
\frac{3}{2}\sqrt{2s}-\frac{15}{2}\sqrt2s^\frac{3}{2}
+\frac{\mathrm{d}\alpha(s)}{\mathrm{d}s}
=\frac{s}{\tilde v}-\frac{1}{\sqrt{2s}}
=\frac{s}{\sqrt2s^\frac{3}{2}-3\sqrt2s^\frac{5}{2}+\alpha(s)}-\frac{1}{\sqrt{2s}}.
$$
It is equivalent to
\begin{align} \nonumber
\frac{\mathrm{d}\alpha(s)}{\mathrm{d}s}
=&\frac{s}{\sqrt2s^\frac{3}{2}-3\sqrt2s^\frac{5}{2}+\alpha(s)}-\frac{1}{\sqrt{2s}}
-\frac{3}{2}\sqrt{2s}+\frac{15}{2}\sqrt2s^\frac{3}{2}
\\[2mm] \label{eq-alpha}
=&\frac{(-1-3s+15s^2)\cdot\alpha(s)+(24\sqrt2s^\frac{7}{2}-45\sqrt2s^\frac{9}{2})}%
{2s^2-6s^3+\sqrt{2s}\cdot\alpha(s)}.
\end{align}
The two special curves that are used to control the trajectories
in the phase plane $(s,\alpha)$ corresponding to \eqref{eq-alpha} are
\begin{equation} \label{eq-zcurve-3}
\bar\Gamma^*(s):=24\sqrt2 s^\frac{7}{2}, \quad
\bar\Gamma_*(s):=24\sqrt2 s^\frac{7}{2}-M_3s^\frac{9}{2},
\quad s\in(0,\delta_3),
\end{equation}
with $M_3=285\sqrt2$ and $\delta_3=\min\{\frac{1}{8},\frac{1}{5},\frac{77}{285}\}=\frac{1}{8}$.
Here we omit the details showing that at any point $(\hat s,\hat \alpha)\in\bar\Gamma^*$
$$
\bar\Gamma^*{}'(\hat s)>0>\frac{\mathrm{d}\alpha(s)}{\mathrm{d}s}%
\Big|_{s=\hat s,\alpha=\hat \alpha},
\quad s\in(0,\delta_3),
$$
and at any point $(\hat s,\hat \alpha)\in\bar\Gamma_*$
$$
\bar\Gamma_*'(\hat s)<\frac{\mathrm{d}\alpha(s)}{\mathrm{d}s}%
\Big|_{s=\hat s,\alpha=\hat \alpha},
\quad s\in(0,\delta_3).
$$
It follows that
$$
24\sqrt2 s^\frac{7}{2}-M_3s^\frac{9}{2}
\le\alpha(s)\le
24\sqrt2 s^\frac{7}{2}, \quad s\in(0,\delta_3),
$$
and further
\begin{equation} \label{eq-tildev3}
\tilde v(s)=\sqrt2s^\frac{3}{2}-3\sqrt2s^\frac{5}{2}+\alpha(s)~
\begin{cases}
\ge \sqrt2s^\frac{3}{2}-3\sqrt2 s^\frac{5}{2}+24\sqrt2 s^\frac{7}{2}-M_3s^\frac{9}{2}, \\[2mm]
\le \sqrt2s^\frac{3}{2}-3\sqrt2 s^\frac{5}{2}+24\sqrt2 s^\frac{7}{2},
\end{cases}
\ \ \ s\in(0,\delta_3).
\end{equation}
Generally, we can take $M_k=|a_{k+1}|$.
The proof is completed.
$\hfill\Box$

We now show the decay estimates of higher order derivatives of $\bar q(x)$.

\begin{lem} \label{le-qxx}
For $s_0<\frac{1}{8}$, let $\bar q(x)$ be the stationary solution
proved in Lemma \ref{le-sx}.
Then
\begin{equation} \label{eq-qxxdecay}
-\frac{3}{\Big(\frac{1}{\sqrt{s_0}}+\frac{x}{2\sqrt2}\Big)^4}
\le \bar q_{xx}(x)\le 0,
\quad x\in(0,+\infty),
\end{equation}
and
\begin{equation} \label{eq-qxxx}
|\bar q_{xxx}(x)|\le
\frac{6\sqrt2}{\Big(\frac{1}{\sqrt{s_0}}+\frac{x}{2\sqrt2}\Big)^5}
+\frac{C}{\Big(\frac{1}{\sqrt{s_0}}+\frac{x}{2\sqrt2}\Big)^7},
\quad x\in(0,+\infty),
\end{equation}
for some positive constant $C>0$.
For $s_0<\delta_4$ ($\delta_4$ is the positive constant in Lemma \ref{le-expansion}),
there holds
\begin{equation} \label{eq-qxxxx}
|\bar q_{xxxx}(x)|
\le \frac{222}{\Big(\frac{1}{\sqrt{s_0}}+\frac{x}{2\sqrt2}\Big)^6}
+\frac{C}{\Big(\frac{1}{\sqrt{s_0}}+\frac{x}{2\sqrt2}\Big)^8},
\end{equation}
for some positive constant $C>0$.
\end{lem}
{\it\bfseries Proof.}
According to the dynamic system \eqref{eq-sv-0},
$$
s_x(x)=-v(x), \qquad v_x(x)=-s(x)+\frac{v(x)}{\sqrt{2s(x)}}.
$$
We have
\begin{align*}
s_{xx}(x)=-v_x(x)=s(x)-\frac{v(x)}{\sqrt{2s(x)}}
=s(x)-\frac{\tilde v(s(x))}{\sqrt{2s(x)}}.
\end{align*}
Using the expansion \eqref{eq-tildev2} in Lemma \ref{le-expansion}
and the decay estimate \eqref{eq-sxdecay} in Lemma \ref{le-sx}, we deduce
$$
s_{xx}(x)\le s(x)-\frac{\sqrt2s^\frac{3}{2}(x)-3\sqrt2 s^\frac{5}{2}(x)}{\sqrt{2s(x)}}
=3s^2(x)
\le \frac{3}{\Big(\frac{1}{\sqrt{s_0}}+\frac{x}{2\sqrt2}\Big)^4},
$$
and
\begin{align*}
s_{xx}(x)\ge &s(x)-\frac{\sqrt2s^\frac{3}{2}(x)-3\sqrt2 s^\frac{5}{2}(x)
+24\sqrt2s^\frac{7}{2}(x)}{\sqrt{2s(x)}}
\\[2mm]
\ge &3s^2(x)-24s^3(x)
\ge \max\left\{\frac{3}{\Big(\frac{1}{\sqrt{s_0}}+\frac{x}{\sqrt2}\Big)^4}
-\frac{24}{\Big(\frac{1}{\sqrt{s_0}}+\frac{x}{2\sqrt2}\Big)^6},
0\right\}.
\end{align*}
Similarly, we have
\begin{align} \nonumber
s_{xxx}(x)=&-v_{xx}(x)=s_x(x)-\frac{v_x(x)}{\sqrt{2s(x)}}
+\frac{v(x)}{2\sqrt2s^\frac{3}{2}(x)}\cdot(-v(x))
\\[2mm] \label{eq-zsxxx}
=&-v(x)-\frac{-s(x)+\frac{v(x)}{\sqrt{2s(x)}}}{\sqrt{2s(x)}}
-\frac{v^2(x)}{2\sqrt2s^\frac{3}{2}(x)}
\\[2mm] \nonumber
=&-\tilde v(s(x))+\frac{\sqrt{s(x)}}{\sqrt2}-\frac{\tilde v(s(x))}{2s(x)}
-\frac{(\tilde v(s(x)))^2}{2\sqrt2s^\frac{3}{2}(x)}.
\end{align}
Utilizing the expansion \eqref{eq-tildev3} in Lemma \ref{le-expansion}, we deduce
\begin{align*}
s_{xxx}(x)\le &-(\sqrt2s^\frac{3}{2}-3\sqrt2 s^\frac{5}{2}+24\sqrt2 s^\frac{7}{2}-285\sqrt2s^\frac{9}{2})
+\frac{\sqrt{s}}{\sqrt2}
\\[2mm]
&-\frac{\sqrt2s^\frac{3}{2}-3\sqrt2 s^\frac{5}{2}+24\sqrt2 s^\frac{7}{2}-285\sqrt2s^\frac{9}{2}}{2s}
\\[2mm]
&-\frac{(\sqrt2s^\frac{3}{2}-3\sqrt2 s^\frac{5}{2}+24\sqrt2 s^\frac{7}{2}-285\sqrt2s^\frac{9}{2})^2}{2\sqrt2s^\frac{3}{2}}
\\[2mm]
\le &-6\sqrt2s^\frac{5}{2}(x)+Cs^\frac{7}{2}(x),
\end{align*}
and
\begin{align*}
s_{xxx}(x)\ge &-(\sqrt2s^\frac{3}{2}-3\sqrt2 s^\frac{5}{2}+24\sqrt2 s^\frac{7}{2})
+\frac{\sqrt{s}}{\sqrt2}
\\[2mm]
&-\frac{\sqrt2s^\frac{3}{2}-3\sqrt2 s^\frac{5}{2}+24\sqrt2 s^\frac{7}{2}}{2s}
-\frac{(\sqrt2s^\frac{3}{2}-3\sqrt2 s^\frac{5}{2}+24\sqrt2 s^\frac{7}{2})^2}{2\sqrt2s^\frac{3}{2}}
\\[2mm]
\ge &-6\sqrt2s^\frac{5}{2}(x)-Cs^\frac{7}{2}(x),
\end{align*}
where $C>0$ is a generic positive constant.
Therefore,
$$
|s_{xxx}(x)|\le 6\sqrt2s^\frac{5}{2}(x)+Cs^\frac{7}{2}(x)
\le
\frac{6\sqrt2}{\Big(\frac{1}{\sqrt{s_0}}+\frac{x}{2\sqrt2}\Big)^5}
+\frac{C}{\Big(\frac{1}{\sqrt{s_0}}+\frac{x}{2\sqrt2}\Big)^7}.
$$
We now show the estimates of $s_{xxxx}(x)$.
According to \eqref{eq-zsxxx},
\begin{align} \nonumber
s_{xxxx}(x)=&-v_x(x)+\Big(\frac{\sqrt {s(x)}}{\sqrt2}\Big)_x
-\Big(\frac{v(x)}{2s(x)}\Big)_x
-\Big(\frac{v^2(x)}{2\sqrt2s^\frac{3}{2}(x)}\Big)_x
\\[2mm] \nonumber
=&-\Big(-s(x)+\frac{v(x)}{\sqrt{2s(x)}}\Big)-\frac{v(x)}{2\sqrt{2s(x)}}
-\frac{-s(x)+\frac{v(x)}{\sqrt{2s(x)}}}{2s(x)}
-\frac{v^2(x)}{2s^2(x)}
\\[2mm] \label{eq-zsxxxx}
&-\frac{2v(x)\cdot(-s(x)+\frac{v(x)}{\sqrt{2s(x)}})}{2\sqrt2s^\frac{3}{2}(x)}
-\frac{\frac{3}{2}v^3(x)}{2\sqrt2s^\frac{5}{2}(x)}.
\end{align}
Substituting $v(x)=\tilde v(s(x))$ and the expansion
of $\tilde v(s)$ for $k=4$ in Lemma \ref{le-expansion} such that
$$
\tilde v(s)~
\begin{cases}
\ge \sqrt2s^\frac{3}{2}-3\sqrt2 s^\frac{5}{2}+24\sqrt2 s^\frac{7}{2}-285\sqrt2s^\frac{9}{2},
\\[2mm]
\le \sqrt2s^\frac{3}{2}-3\sqrt2 s^\frac{5}{2}+24\sqrt2 s^\frac{7}{2}-285\sqrt2s^\frac{9}{2}
+4284\sqrt2s^\frac{11}{2},
\end{cases}
\quad s\in(0,\delta_4),
$$
into \eqref{eq-zsxxxx} implies that
$$
|s_{xxxx}(x)|\le 222s^3(x)+Cs^4(x)
\le \frac{222}{\Big(\frac{1}{\sqrt{s_0}}+\frac{x}{2\sqrt2}\Big)^6}
+\frac{C}{\Big(\frac{1}{\sqrt{s_0}}+\frac{x}{2\sqrt2}\Big)^8},
$$
for some positive constant $C>0$.
$\hfill\Box$

\begin{lem} \label{le-u}
For $s_0=\frac{1}{2}u_-^2-\frac{1}{2}u_+^2=\frac{1}{2}u_-^2=\frac{1}{2}\delta^2
<\min\{\frac{1}{8},\delta_4\}$,
let $\bar q(x)$ be the stationary solution
proved in Lemma \ref{le-sx}
and let $\bar u(x):=-\sqrt{-2\bar q(x)}=-\sqrt{2s(x)}$.
Then we have
\begin{align*}
\frac{C_1\delta}{1+\delta x}\le-\bar u(x)\le \frac{C_2\delta}{1+\delta x},
\qquad
&\frac{C_1\delta^2}{(1+\delta x)^2}\le\bar u_x(x)\le \frac{C_2\delta^2}{(1+\delta x)^2},
\\[3mm]
\Big|\partial^k_x\bar u(x)\Big|\le \frac{C_2\delta^{k+1}}{(1+\delta x)^{k+1}},
\qquad
&\Big|\frac{\bar u_{xx}^2(x)}{\bar u_x(x)}\Big|\le \frac{C_2\delta^4}{(1+\delta x)^4},
\\[3mm]
\Big|\frac{\bar u_{xxx}^2(x)}{\bar u_x(x)}\Big|\le \frac{C_2\delta^6}{(1+\delta x)^6},
\qquad
\Big|\frac{\bar u_{xxx}(x)}{\bar u_x(x)}\Big|\le &\frac{C_2\delta^2}{(1+\delta x)^2},
\qquad
\Big|\frac{\bar u_{xxxx}(x)}{\bar u_x(x)}\Big|\le \frac{C_2\delta^3}{(1+\delta x)^3},
\end{align*}
for $k=1,2,3,4$ and some positive constants $C_1$, $C_2$.
\end{lem}
{\it\bfseries Proof.}
The lower and upper bounds of $\bar q(x)$ and $\bar q_x(x)$ in \eqref{eq-qdecay-1}
and \eqref{eq-qdecay-2} in Lemma \ref{le-sx}
show that
$$
\bar u(x)=-\sqrt{-2\bar q(x)}~
\begin{cases}
\displaystyle
\le -\frac{\sqrt2}{\frac{1}{\sqrt{s_0}}+\frac{x}{\sqrt2}},
\\[6mm]
\displaystyle
\ge -\frac{\sqrt2}{\frac{1}{\sqrt{s_0}}+\frac{x}{2\sqrt2}},
\end{cases}
\quad x\in(0,+\infty),
$$
and
$$
\bar u_x(x)=\frac{\bar q_x}{\sqrt{-2\bar q}}~
\begin{cases}
\displaystyle
\le \frac{\frac{1}{\sqrt{s_0}}+\frac{x}{\sqrt2}}%
{\Big(\frac{1}{\sqrt{s_0}}+\frac{x}{2\sqrt2}\Big)^3}
\le \frac{2}%
{\Big(\frac{1}{\sqrt{s_0}}+\frac{x}{2\sqrt2}\Big)^2},
\\[8mm]
\displaystyle
\ge \frac{\frac{1}{\sqrt{s_0}}+\frac{x}{2\sqrt2}}%
{2\Big(\frac{1}{\sqrt{s_0}}+\frac{x}{\sqrt2}\Big)^3}
\ge \frac{1}%
{4\Big(\frac{1}{\sqrt{s_0}}+\frac{x}{\sqrt2}\Big)^2},
\end{cases}
\quad x\in(0,+\infty).
$$
Utilizing the higher order estimates \eqref{eq-qxxdecay},
\eqref{eq-qxxx} and \eqref{eq-qxxxx} in Lemma \ref{le-qxx},
we have
\begin{align*}
|\bar u_{xx}(x)|=&\left|\frac{\bar q_{xx}(x)}{\sqrt{-2\bar q(x)}}
+\frac{(\bar q_x(x))^2}{(-2\bar q(x))^\frac{3}{2}}\right|
\le
\frac{3(\frac{1}{\sqrt{s_0}}+\frac{x}{\sqrt2})}%
{\sqrt2\Big(\frac{1}{\sqrt{s_0}}+\frac{x}{2\sqrt2}\Big)^4}
+\frac{\Big(\frac{1}{\sqrt{s_0}}+\frac{x}{\sqrt2}\Big)^3}%
{\sqrt2\Big(\frac{1}{\sqrt{s_0}}+\frac{x}{2\sqrt2}\Big)^6}
\\[4mm]
\le&\frac{7\sqrt2}%
{\Big(\frac{1}{\sqrt{s_0}}+\frac{x}{2\sqrt2}\Big)^3},
\end{align*}
and
$$
|\bar u_{xxx}(x)|=\left|\frac{\bar q_{xxx}}{\sqrt{-2\bar q}}
+\frac{3\bar q_x\bar q_{xx}}{(-2\bar q)^\frac{3}{2}}
+\frac{3(\bar q_x)^3}{(-2\bar q)^\frac{5}{2}}\right|
\le\frac{C}%
{\Big(\frac{1}{\sqrt{s_0}}+\frac{x}{2\sqrt2}\Big)^4},
$$
$$
|\bar u_{xxxx}(x)|=\left|
\frac{\bar q_{xxxx}}{\sqrt{-2\bar q}}+\frac{4\bar q_x\bar q_{xxx}+3(\bar q_{xx})^2}{(-2\bar q)^\frac{3}{2}}
+\frac{18(\bar q_x)^2\bar q_{xx}}{(-2\bar q)^\frac{5}{2}}
+\frac{15(\bar q_x)^4}{(-2\bar q)^\frac{7}{2}}
\right|
\le \frac{C}%
{\Big(\frac{1}{\sqrt{s_0}}+\frac{x}{2\sqrt2}\Big)^5},
$$
where $C>0$ is a generic positive constant.
Moreover,
$$
\left|\frac{\bar u_{xxx}}{\bar u_x}\right|
\le \frac{C}%
{\Big(\frac{1}{\sqrt{s_0}}+\frac{x}{2\sqrt2}\Big)^2},
\qquad
\left|\frac{\bar u_{xxxx}}{\bar u_x}\right|
\le \frac{C}%
{\Big(\frac{1}{\sqrt{s_0}}+\frac{x}{2\sqrt2}\Big)^3},
$$
and
$$
\left|\frac{\bar u_{xxx}^2}{\bar u_x}\right|
\le \frac{C}%
{\Big(\frac{1}{\sqrt{s_0}}+\frac{x}{2\sqrt2}\Big)^6},
\qquad
\left|\frac{\bar u_{xx}^2}{\bar u_x}\right|
\le \frac{C}%
{\Big(\frac{1}{\sqrt{s_0}}+\frac{x}{2\sqrt2}\Big)^4}.
$$
The proof is completed.
$\hfill\Box$

{\it\bfseries Proof of Lemma \ref{lem2.1.1}.}
The degenerate case (2): $u_-<u_+=0$ is proved through a singular phase plane analysis method
according to Lemma \ref{le-sx} and Lemma \ref{le-u}.

Next we show that this method is applicable to the case (1) $u_-<u_+<0$.
Instead of \eqref{eq-sv-0}, we have a non-degenerate dynamical system \eqref{eq-sv} for the case of $u_+<0$.
We sketch the main lines of the proof.

(1) For any $s_0:=\frac{1}{2}u_-^2-\frac{1}{2}u_+^2>0$, if $\tilde v(s)$ solves the following equation
\begin{equation} \label{eq-v-nd}
\begin{cases}
\displaystyle
\frac{\mathrm{d}\tilde v}{\mathrm{d}s}=\frac{s}{\tilde v}-\frac{1}{\sqrt{u_+^2+2s}}, \quad s\in(0,s_0),\\[6mm]
\displaystyle
\tilde v(s)>0 \mathrm{~for~}s\in(0,s_0), \quad \lim_{s\rightarrow0^+}\tilde v(s)=0,
\quad \int_0^{s_0}\frac{1}{\tilde v(s)}\mathrm{d}s=+\infty,
\end{cases}
\end{equation}
then the function $s(x)$ defined by
\begin{equation} \label{eq-sx-nd}
x=-\int_{s_0}^{s(x)}\frac{1}{\tilde v(\tau)}\mathrm{d}\tau, \qquad x\in(0,+\infty)
\end{equation}
is a solution of \eqref{eq-sv}
with $v(x):=\tilde v(s(x))>0$.

(2) Consider the approximated problem
\begin{equation} \label{eq-vk-nd}
\begin{cases}
\displaystyle
\frac{\mathrm{d}\tilde v_k}{\mathrm{d}s}=\frac{s}{\tilde v_k}-\frac{1}{\sqrt{u_+^2+2s}},
\quad s\in(0,s_0),\\[6mm]
\displaystyle
\tilde v_k(s)>0 \mathrm{~for~}s\in(0,s_0),
\quad \tilde v_k(0)=\frac{1}{k}.
\end{cases}
\end{equation}
For any $s_0>0$, the problem \eqref{eq-vk-nd} admits
a solution $\tilde v_k(s)>0$ on $s\in(0,s_0)$
such that
\\[2mm] \indent
(i) $\tilde v_k(s)$ is monotone decreasing with respect to $k$, i.e.,
$\tilde v_{k_1}(s)>\tilde v_{k_2}(s)$ for any $k_2>k_1$ on $(0,s_0)$;
\\[2mm] \indent
(ii) $\tilde v_k(s)$ has the following upper bound estimate
$$\tilde v_k(s)\le \overline\Gamma_k(s):=
\max\Big\{\frac{1}{k},s\cdot\sqrt{u_+^2+2s}\Big\},
\quad s\in (0,s_0);
$$

(iii) $\tilde v_k(s)$ has the following uniformly
lower bound estimate
$$\tilde v_k(s)\ge \lambda_0s, \quad s\in (0,s_0),$$
where $\lambda_0:=\frac{\sqrt{1+4|u_+|^2}-1}{2|u_+|}$.

(3) The limit function
$$\tilde v(s):=\lim_{k\rightarrow\infty}\tilde v_k(s), \quad s\in(0,s_0)$$
is well-defined and $\tilde v(s)$ is a solution to the problem \eqref{eq-v-nd}.
Moreover,
\begin{equation} \label{eq-tildev-nd}
\lambda_0s\le \tilde v(s)\le \lambda_0s+bs^2, \quad s\in(0,s_0),
\end{equation}
where $\lambda_0:=\frac{\sqrt{1+4|u_+|^2}-1}{2|u_+|}$ is the positive root of
$\lambda_0=\frac{1}{\lambda_0}-\frac{1}{|u_+|}$
and $b=\frac{1}{2|u_+|^3}$.

The asymptotic expansion \eqref{eq-tildev-nd} plays an essential role in the analysis of
asymptotic decay behavior of the stationary solution,
thus we present the following proof.
In the phase plane $(s,v)$, at any point $(\hat s,\hat v)$ on the curve $\hat v=\Gamma_*(\hat s):=\lambda_0\hat s$,
we have
$$
\Gamma_*'(\hat s)=\lambda_0=\frac{1}{\lambda_0}-\frac{1}{|u_+|}<
\frac{\hat s}{\hat v}-\frac{1}{\sqrt{u_+^2+2\hat s}}
=\left.\frac{\mathrm{d}\tilde v}{\mathrm{d}s}\right|_{(\hat s,\hat v)},
$$
which means the trajectory lies above $\Gamma_*$.
At any point $(\hat s,\hat v)$ on the curve $\hat v=\Gamma^*(\hat s):=\lambda_0\hat s+b\hat s^2$,
we have
$$
\Gamma^*{}'(\hat s)=\lambda_0+2b\hat s
>\frac{\hat s}{\hat v}-\frac{1}{\sqrt{u_+^2+2\hat s}}
=\left.\frac{\mathrm{d}\tilde v}{\mathrm{d}s}\right|_{(\hat s,\hat v)},
$$
since the following auxiliary function
$$
F(\hat s):=\lambda_0+2b\hat s-\frac{1}{\lambda_0+b\hat s}+\frac{1}{\sqrt{u_+^2+2\hat s}}, \quad s\in(0,s_0)
$$
is monotonically increasing as
$F'(\hat s)=2b+\frac{b}{(\lambda_0+b\hat s)^2}-\frac{1}{(|u_+|^2+2\hat s)^\frac{3}{2}}>0$ for $b=\frac{1}{2|u_+|^3}$.
This shows the trajectory lies between $\Gamma_*$ and $\Gamma^*$.

(4) Finally we show the asymptotic behavior of the stationary solution.
According to the definition of $s(x)$ and $\bar q(x)$ in \eqref{eq-sx-nd}
and the asymptotic expansion \eqref{eq-tildev-nd}, we have
\begin{align*}
x=-\int_{s_0}^{s(x)}\frac{1}{\tilde v(\tau)}\mathrm{d}\tau
\le \int_{s(x)}^{s_0}\frac{1}{\lambda_0\tau}\mathrm{d}\tau
=\frac{1}{\lambda_0}\ln\Big(\frac{s_0}{s(x)}\Big),
\quad x\in(0,+\infty),
\end{align*}
which implies
$$
s(x)\le s_0\mathrm{e}^{-\lambda_0x}, \quad \bar q(x)\ge \bar q_-\mathrm{e}^{-\lambda_0x}.
$$
On the other hand,
\begin{align*}
x=-\int_{s_0}^{s(x)}\frac{1}{\tilde v(\tau)}\mathrm{d}\tau
\ge \int_{s(x)}^{s_0}\frac{1}{\lambda_0\tau+b\tau^2}\mathrm{d}\tau
=\left.\frac{1}{\lambda_0}\ln\left(\frac{\tau}{1+\frac{b}{\lambda_0}\tau}\right)\right|_{s(x)}^{s_0},
\qquad x\in(0,+\infty).
\end{align*}
That is,
$$
s(x)\ge\frac{\frac{s_0}{1+\frac{b}{\lambda_0}s_0}\cdot \mathrm{e}^{-\lambda_0 x}}{1-\frac{b}%
{\lambda_0}\cdot\frac{s_0}{1+\frac{b}{\lambda_0}s_0}\cdot \mathrm{e}^{-\lambda_0 x}}
\ge \frac{s_0}{1+\frac{b}{\lambda_0}s_0}\cdot \mathrm{e}^{-\lambda_0 x},
$$
which shows
$$
\bar q(x)\le \frac{\bar q_-}{1-\frac{b}{\lambda_0}\bar q_-}\cdot \mathrm{e}^{-\lambda_0 x}.
$$
The decay estimates of $\bar q_x(x)$ and other higher order derivatives follow similarly,
which are all exponentially decaying, since
$$s_x(x)=-v(x)=-\tilde v(s(x)),$$
and
$$s_{xx}(x)=-v_x(x)=s(x)-\frac{v(x)}{\sqrt{u_+^2+2s(s)}}=s(x)-\frac{\tilde v(s(x))}{\sqrt{u_+^2+2s(s)}},$$
and according to the asymptotic expansion \eqref{eq-tildev-nd}.
The proof is completed.
$\hfill\Box$

\subsection{Preliminary Lemmas}

In order to show the asymptotic behavior of heat flux $q(x,t)$ which satisfies an elliptic problem,
we prove the following optimal Gagliardo-Nirenberg-Sobolev inequality.
This inequality without optimal constant is known as a special case of Gagliardo-Nirenberg-Sobolev inequality.
Here we present a primary proof based on fundamental calculus.

\begin{lem}[Optimal Gagliardo-Nirenberg-Sobolev inequality] \label{le-GN}
For any function $u(x)\in L^\infty(\mathbbm R)$ with $u_{xx}\in L^\infty(\mathbbm R)$, there holds
\begin{equation} \label{eq-GN}
\|u_x\|_{L^\infty(\mathbbm R)}\le \sqrt2 \
\|u\|_{L^\infty(\mathbbm R)}^\frac{1}{2}\cdot\|u_{xx}\|_{L^\infty(\mathbbm R)}^\frac{1}{2},
\end{equation}
and the constant $\sqrt2$ is optimal.
Moreover, for any function $u(x)\in L^\infty(\mathbbm R_+)$ with $u_{xx}\in L^\infty(\mathbbm R_+)$, there holds
\begin{equation} \label{eq-GN2}
\|u_x\|_{L^\infty(\mathbbm R_+)}\le 2 \ \|u\|_{L^\infty(\mathbbm R_+)}^\frac{1}{2}\cdot\|u_{xx}\|_{L^\infty(\mathbbm R_+)}^\frac{1}{2},
\end{equation}
and the constant $2$ is optimal.
For multi-dimensional case,
\begin{equation} \label{eq-GN-m}
\||\nabla u|\|_{L^\infty(\mathbbm R^n)}
\le \sqrt2\ \|u\|_{L^\infty(\mathbbm R^n)}^\frac{1}{2}\cdot\||D^2 u|\|_{L^\infty(\mathbbm R^n)}^\frac{1}{2},
\qquad \forall u\in W^{2,\infty}(\mathbbm R^n),
\end{equation}
and
\begin{equation} \label{eq-GN-m2}
\||\nabla u|\|_{L^\infty(\mathbbm R_+^n)}
\le 2\ \|u\|_{L^\infty(\mathbbm R_+^n)}^\frac{1}{2}\cdot\||D^2 u|\|_{L^\infty(\mathbbm R_+^n)}^\frac{1}{2},
\qquad \forall u\in W^{2,\infty}(\mathbbm R_+^n),
\end{equation}
and the constants $\sqrt2$ and $2$ are optimal,
where $|\nabla u|$ is the modulus of a vector $\nabla u$
and $|D^2u|$ is the spectral norm of a matrix $D^2u$ such that $|D^2u|:=\sup_{p\in \mathbbm R^n;|p|=1}p^\mathrm{T}(D^2u)p$.
\end{lem}
{\it\bfseries Proof.}
We prove that the inequality \eqref{eq-GN} holds for any smooth function $u\in C^2(\mathbbm R)\cap W^{2,\infty}(\mathbbm R)$
and the constant is optimal for $u\in W^{2,\infty}(\mathbbm R)$ with weak derivatives.
Then utilizing an approximation approach, we see that \eqref{eq-GN} holds for $u\in W^{2,\infty}(\mathbbm R)$
with the same optimal constant.

The inequality \eqref{eq-GN} is trivial if $\|u\|_{L^\infty(\mathbbm R)}=0$ or $\|u_{xx}\|_{L^\infty(\mathbbm R)}=0$,
since in the latter case $u(x)=c_1x+c_2$ and $c_1=0$ according to $u\in L^\infty(\mathbbm R)$.
With the observation that the inequality \eqref{eq-GN} is invariant under the scaling
$\tilde u(x):=\lambda u(\mu x)$ for any non-zero $\lambda$ and $\mu$, we only need to prove that
$\|u_x\|_{L^\infty(\mathbbm R)}\le 1$ under the condition $\|u\|_{L^\infty(\mathbbm R)}=\frac{1}{2}$ and
$\|u_{xx}\|_{L^\infty(\mathbbm R)}=1$, and further $\|u_x\|_{L^\infty(\mathbbm R)}\le 1$ is optimal.
In other words, we show that if $u_x(x_0)=1$ for some $x_0\in\mathbbm R$
and $\|u_{xx}\|_{L^\infty(\mathbbm R)}\le1$ then $\|u\|_{L^\infty(\mathbbm R)}\ge\frac{1}{2}$ and
$\frac{1}{2}$ is optimal.

According to Taylor expansion near $x_0$ for $u\in C^2(\mathbbm R)\cap W^{2,\infty}(\mathbbm R)$, we know that
\begin{align*}
u(x_0+1)&\ge u(x_0)+u_x(x_0)\cdot1+\frac{1}{2}u_{xx}(\xi_1)\cdot1^2
\ge u(x_0)+1-\frac{1}{2},
\\
u(x_0-1)&\le u(x_0)+u_x(x_0)\cdot(-1)+\frac{1}{2}u_{xx}(\xi_2)\cdot(-1)^2
\le u(x_0)-1+\frac{1}{2},
\end{align*}
with some $\xi_1\in(x_0,x_0+1)$ and $\xi_2\in(x_0-1,x_0)$.
Therefore,
$$
u(x_0+1)-u(x_0-1)\ge 1,
$$
and then $\|u\|_{L^\infty(\mathbbm R)}\ge\frac{1}{2}$.
The constant is optimal for the following $\hat u_1\in W^{2,\infty}(\mathbbm R)$
\begin{equation} \label{eq-GNoptimal}
\hat u_1(x):=
\begin{cases}
\cdots,\\[1mm]
x(1-\frac{1}{2}x), \ \ \ \ x\in[0,2),\\[1mm]
(x-2)(2-\frac{1}{2}x), \ \ \ \ x\in[2,4),\\[1mm]
\cdots
\end{cases}
\end{equation}
which is defined by extension as a $4$-periodic function.
We can verify that $\hat u_1(x)$ satisfies the following differential equation
$$
-\hat u_{1xx}=\mathrm{sign}u_1=\begin{cases}
1, \ \ \ \ u_1>0,\\[1mm]
-1, \ \ \ \ u_1<0.
\end{cases}
$$

The inequality \eqref{eq-GN2} for $u\in W^{2,\infty}(\mathbbm R_+)$ is proved by extension
$$\tilde u(x):=
\begin{cases}
u(x), \ \ \ \ x\ge0,\\[2mm]
2u(0)-u(-x), \ \ \ \ x<0,
\end{cases}
$$
such that $\|\tilde u_x\|_{L^\infty(\mathbbm R)}=\|u_x\|_{L^\infty(\mathbbm R_+)}$,
$\|\tilde u_{xx}\|_{L^\infty(\mathbbm R)}=\|u_{xx}\|_{L^\infty(\mathbbm R_+)}$,
and
$$\mathrm{osc}(u;\mathbbm R_+) \le\mathrm{osc}(\tilde u;\mathbbm R)
\le2\mathrm{osc}(u;\mathbbm R_+)\le4\|u\|_{L^\infty(\mathbbm R_+)}.$$
Here $\mathrm{osc}(f;D):=\sup_{x,y\in D}|f(x)-f(y)|$ is the oscillation of a given function $f$ and $D$ is its domain of definition.
Furthermore, applying \eqref{eq-GN}
(according to the proof, we can replace $\|\tilde u\|_{L^\infty(\mathbbm R)}$ by $\frac{1}{2}\mathrm{osc}(\tilde u;\mathbbm R)$)
$$
\|u_x\|_{L^\infty(\mathbbm R_+)}=\|\tilde u_x\|_{L^\infty(\mathbbm R)}
\le (\mathrm{osc}(\tilde u;\mathbbm R))^\frac{1}{2}\cdot\|\tilde u_{xx}\|_{L^\infty(\mathbbm R)}^\frac{1}{2}
\le 2\ \|u\|_{L^\infty(\mathbbm R)}^\frac{1}{2}\cdot\|u_{xx}\|_{L^\infty(\mathbbm R)}^\frac{1}{2}.
$$
The constant $2$ is optimal for the following $\hat u_2(x)$
\begin{equation} \label{eq-GNoptimal2}
\hat u_2(x):=
\begin{cases}
x(1-\frac{1}{2}x)-\frac{1}{4}, \ \ \ \ x\in[0,1),\\[2mm]
\frac{1}{4}, \ \ \ \ x\in[1,+\infty),
\end{cases}
\end{equation}
such that $\|\hat u_2\|_{L^\infty(\mathbbm R_+)}=\frac{1}{4}$, $\|\hat u_{2x}\|_{L^\infty(\mathbbm R_+)}=1$
and $\|\hat u_{2xx}\|_{L^\infty(\mathbbm R_+)}=1$.

For the multi-dimensional case, we note that the inequality \eqref{eq-GN-m} is invariant under the scaling
$\tilde u(x):=\lambda u(\mu x)$ for any non-zero $\lambda$ and $\mu$,
and is also invariant under the rotation of coordinates.
Therefore, for any function $u\in C^2(\mathbbm R^n)\cap W^{2,\infty}(\mathbbm R^n)$,
if $|\nabla u(x_0)|=1$ for some $x_0\in\mathbbm R^n$
and $|D^2 u(x)|\le1$ for all $x\in\mathbbm R^n$, then Taylor expansion along the direction $\pm p:=\pm \nabla u(x_0)$
shows that
\begin{align*}
u(x_0+\nabla u(x_0))&\ge u(x_0)+\nabla u(x_0)\cdot\nabla u(x_0)+\frac{1}{2}(\nabla u(x_0))^\mathrm{T}(D^2u(\xi_1))(\nabla u(x_0))
\ge u(x_0)+1-\frac{1}{2},
\\
u(x_0-\nabla u(x_0))&\le u(x_0)-\nabla u(x_0)\cdot\nabla u(x_0)+\frac{1}{2}(-\nabla u(x_0))^\mathrm{T}(D^2u(\xi_2))(-\nabla u(x_0))
\le u(x_0)-1+\frac{1}{2},
\end{align*}
for some $\xi_1=x_0+\theta_1\nabla u(x_0)$ and $\xi_2=x_0-\theta_2\nabla u(x_0)$
with $\theta_1,\theta_2\in(0,1)$.
The rest of the proof follows similarly.
$\hfill\Box$

\begin{remark}
Lemma \ref{le-GN} can be seen as a special case of Gagliardo-Nirenberg-Sobolev inequality
with the optimal constant and without the restriction of decay at infinity such that $\lim_{|x|\rightarrow+\infty}u(x)=0$.
\end{remark}

\begin{lem} \label{le-Bessel}
Assume that $f\in L^\infty(\mathbbm R_+)$, $g\in\mathbbm R$,
and $u\in W^{2,\infty}(\mathbbm R_+)$ solves the following elliptic problem
\begin{equation} \label{eq-Bessel}
\begin{cases}
-u_{xx}+u=f(x), \quad & x>0, \\[1mm]
u_x(0)=g,
\end{cases}
\end{equation}
then
\begin{equation} \label{eq-Bessel-est}
\|u\|_{L^\infty(\mathbbm R_+)}\le \|f\|_{L^\infty(\mathbbm R_+)}+|g|,
\quad
\|u_x\|_{L^\infty(\mathbbm R_+)}\le \|f\|_{L^\infty(\mathbbm R_+)}+|g|,
\quad
\|u_{xx}\|_{L^\infty(\mathbbm R_+)}\le 2\|f\|_{L^\infty(\mathbbm R_+)}+|g|,
\end{equation}
and all the above coefficients are optimal.
\end{lem}
{\it\bfseries Proof.}
Let $v(x):=u(x)+g\mathrm{e}^{-x}$. Then $v\in W^{2,\infty}(\mathbbm R_+)$ satisfies
\begin{equation} \label{eq-Bessel2}
\begin{cases}
-v_{xx}+v=f(x), \quad & x>0, \\[1mm]
v_x(0)=0.
\end{cases}
\end{equation}
Maximum principle shows that
$\|v\|_{L^\infty(\mathbbm R_+)}\le \|f\|_{L^\infty(\mathbbm R_+)}$ and then
$\|u\|_{L^\infty(\mathbbm R_+)}\le \|f\|_{L^\infty(\mathbbm R_+)}+|g|$.
Further, according to the equation \eqref{eq-Bessel} we have
$\|u_{xx}\|_{L^\infty(\mathbbm R_+)}\le 2\|f\|_{L^\infty(\mathbbm R_+)}+|g|$.
According to Lemma \ref{le-GN}, we see that
$$\|u_x\|_{L^\infty(\mathbbm R_+)}\le2\
\|u\|_{L^\infty(\mathbbm R_+)}^\frac{1}{2}\cdot\|u_{xx}\|_{L^\infty(\mathbbm R_+)}^\frac{1}{2}
\le 2\ (\|f\|_{L^\infty(\mathbbm R_+)}+|g|)^\frac{1}{2}\cdot(2\|f\|_{L^\infty(\mathbbm R_+)}+|g|)^\frac{1}{2}.$$
Here the Gagliardo-Nirenberg-Sobolev inequality in Lemma \ref{le-GN} is optimal for all $u\in W^{2,\infty}(\mathbbm R_+)$
but not for the solutions of elliptic problem \eqref{eq-Bessel}.
In order to show optimal estimates, we extend the functions $f$ and $v$ in \eqref{eq-Bessel2}
such that
$$
\tilde f(x):=\begin{cases}
f(x), \ \ \ \ x\ge0,\\
f(-x), \ \ \ \ x<0,
\end{cases}
\quad
\text{~and~}
\quad
\tilde v(x):=\begin{cases}
v(x), \ \ \ \ x\ge0,\\
v(-x), \ \ \ \ x<0.
\end{cases}
$$
Then $\tilde v$ can be solved as
$$
\tilde v(x)=\frac{1}{2}\int_{\mathbbm R}\mathrm{e}^{-|x-y|}\tilde f(y)\mathrm{d}y=
\frac{1}{2}\int_0^{+\infty}\left(\mathrm{e}^{-|x-y|}+\mathrm{e}^{-|x+y|}\right)f(y)\mathrm{d}y,
\quad x\in\mathbbm R.
$$
Therefore, for $x\ge0$,
\begin{align} \label{eq-Bessel-u}
u(x)&=\frac{1}{2}\int_0^{+\infty}\left(\mathrm{e}^{-|x-y|}+\mathrm{e}^{-|x+y|}\right)f(y)\mathrm{d}y-g\mathrm{e}^{-x},
\\[3mm] \label{eq-Bessel-ux}
u_x(x)&=\frac{1}{2}\int_0^{+\infty}\left(\mathrm{e}^{-|x-y|}(-\mathrm{sign(x-y)})+\mathrm{e}^{-|x+y|}(-1)\right)f(y)\mathrm{d}y+g\mathrm{e}^{-x},
\end{align}
and
\begin{equation} \label{eq-Bessel-uxx}
u_{xx}(x)=\frac{1}{2}\int_0^{+\infty}\left(\mathrm{e}^{-|x-y|}+\mathrm{e}^{-|x+y|}\right)f(y)\mathrm{d}y-g\mathrm{e}^{-x}-f(x).
\end{equation}
The above expressions show that the estimates \eqref{eq-Bessel-est} are valid.

Now we show that all these coefficients are optimal.
The special case of $f(x)\equiv0$, $g=-1$ and $u(x)=\mathrm{e}^{-x}$
implies that the coefficients of $|g|$ in \eqref{eq-Bessel-est} are optimal.
The case of $g=0$, $f(x)\equiv1$ and $u(x)\equiv1$ shows that the coefficient of $\|f\|_{L^\infty(\mathbbm R_+)}$
in the estimate of $\|u\|_{L^\infty(\mathbbm R_+)}$ is optimal.
For any large $x_0>0$, we set $f(x)=\mathrm{sign}(x-x_0)$ and $g=0$, then \eqref{eq-Bessel-ux} implies
\begin{align*}
u_x(x_0)=&\frac{1}{2}\int_0^{x_0}\left(\mathrm{e}^{-|x_0-y|}+\mathrm{e}^{-|x_0+y|}\right)\mathrm{d}y
+\frac{1}{2}\int_{x_0}^{+\infty}\left(\mathrm{e}^{-|x_0-y|}-\mathrm{e}^{-|x_0+y|}\right)\mathrm{d}y
\\[3mm]
=&\frac{1}{2}\int_0^{2x_0}\mathrm{e}^{-y}\mathrm{d}y
+\frac{1}{2}\int_0^{+\infty}\mathrm{e}^{-y}\mathrm{d}y
-\frac{1}{2}\int_{2x_0}^{+\infty}\mathrm{e}^{-y}\mathrm{d}y
\\[3mm]
=&1-\int_{2x_0}^{+\infty}\mathrm{e}^{-y}\mathrm{d}y
\rightarrow1, \quad \text{as~} x_0\rightarrow+\infty,
\end{align*}
which shows that the coefficient of $\|f\|_{L^\infty(\mathbbm R_+)}$
in the estimate of $\|u_x\|_{L^\infty(\mathbbm R_+)}$ is optimal.
Lastly, for $g=0$, small $\varepsilon>0$, and $f(x)=\mathrm{sign}(x-\varepsilon)$, according to \eqref{eq-Bessel-uxx}, we have
\begin{align*}
u_{xx}(0)&=\int_0^{+\infty}\mathrm{e}^{-y}f(y)\mathrm{d}y-f(0)
\\[3mm]
&=\int_0^{\varepsilon}\mathrm{e}^{-y}(-1)\mathrm{d}y+\int_{\varepsilon}^{+\infty}\mathrm{e}^{-y}(+1)\mathrm{d}y-(-1)
\\[3mm]
&=2-2\int_0^{\varepsilon}\mathrm{e}^{-y}\mathrm{d}y\rightarrow2, \quad \text{as~}~~ \varepsilon\rightarrow0,
\end{align*}
and the coefficient of $\|f\|_{L^\infty(\mathbbm R_+)}$
in the estimate of $\|u_{xx}\|_{L^\infty(\mathbbm R_+)}$ is optimal.
$\hfill\Box$

\subsection{Main Theorems}

We state our main results for the cases: $u_-<u_+\le 0$, $ 0\le u_-<u_+$ and $u_-<0<u_+$,
that the initial-boundary value problem \eqref{eq1.1}-\eqref{eq1.3} admits a unique global solution and it converges to the stationary solution, the rarefaction wave and the superposition of the nonlinear waves, respectively, as $t\to \infty$.

\begin{thm}[In the case of $0< u_-<u_+$] \label{thmxsb1}
Suppose that the boundary condition and far field states satisfy $0< u_-<u_+$, the initial data $u_0 $ satisfies $u_0-\tilde{u}_0^R\in H^2(\mathbbm{R}_+)$, where $\tilde{u}_0^R$ is defined in \eqref{xsbcsz}.
Also assume that $\delta=|u_--u_+|$ is sufficiently small. Then there is a positive constant $\epsilon_0$ such that if  $\|u_0-\tilde{u}_0^R\|_2+\delta \le \epsilon_0$, the problem \eqref{eq1.1}-\eqref{eq1.3} admits a unique solution $(u(x,t),q(x,t))$, which satisfies
\begin{equation*}
\begin{aligned}
 &u- u^R \in C^0([0,\infty);H^2)\cap C^1([0,\infty);H^1),\\[2mm]
 &q+\partial_xu^R \in C^0([0,\infty);H^3)\cap L^2(0,\infty;H^3),
\end{aligned}
\end{equation*}
and the asymptotic behavior
\begin{equation*}
\begin{split}
 &\sup_{x\in \mathbbm{R}_+} \left|\partial_x^k\left(u(x,t)- u^R(x)\right)\right|\rightarrow 0 \ \  \text{as} \ \ t\rightarrow \infty, \ \ k=0,1,\\[2mm]
 &\sup_{x\in \mathbbm{R}_+} \left|\partial_x^k\left(q(x,t)+\partial_xu^R(x)\right)\right|\rightarrow 0 \ \  \text{as} \ \ t\rightarrow \infty, \ \ k=0,1,2.
\end{split}
\end{equation*}
\end{thm}

\begin{thm}[In the case of $u_-<u_+\le0$] \label{thm2.2.1}
Suppose that the boundary condition and far field states satisfy $u_-<u_+\le0$, the initial data $u_0 $ satisfies $u_0-\bar{u}\in H^2(\mathbbm{R}_+)$, where $\bar{u}=\bar{u}_i,(i=1,2)$ is a stationary solution of Lemma \ref{lem2.1.1}.
Also assume that $\delta=|u_--u_+|$ is sufficiently small.
Then there is a positive constant $\epsilon_0$ such that if $\|u_0-\bar{u}\|_2+\delta \le \epsilon_0$, the problem \eqref{eq1.1}-\eqref{eq1.3} admits a unique solution $(u(x,t),q(x,t))$, which satisfies
\begin{equation*}
\begin{aligned}
 &u-\bar{u}\in C^0([0,\infty);H^2),\ \ u_x-\bar{u}_x \in L^2(0,\infty;H^1),\\[2mm]
 &q-\bar{q}\in C^0([0,\infty);H^3)\cap L^2(0,\infty;H^3),
\end{aligned}
\end{equation*}
and the asymptotic behavior
\begin{equation*}
\begin{split}
 &\sup_{x\in \mathbbm{R}_+} \left|\partial_x^k\left(u(x,t)-\bar{u}(x)\right)\right|\rightarrow 0 \ \  \text{as} \ \ t\rightarrow \infty, \ \ k=0,1,\\[2mm]
 &\sup_{x\in \mathbbm{R}_+} \left|\partial_x^k\left(q(x,t)-\bar{q}(x)\right)\right|\rightarrow 0 \ \  \text{as} \ \ t\rightarrow \infty, \ \ k=0,1,2.
\end{split}
\end{equation*}
\end{thm}

\begin{thm}[In the case of $u_-<0<u_+$] \label{thmfhb1}
Suppose that the boundary condition and far field states satisfy $u_-<0<u_+$, the initial data $u_0 $ satisfies $u_0-\bar{u}_2(\cdot)-\tilde{u}_4(\cdot,0)\in H^2(\mathbbm{R}_+)$,
where $\bar{u}_2$ and $\tilde{u}_4$ is a stationary solution and rarefaction wave for the cases $u_-<u_+=0$ and $u_+>u_-= 0$, respectively.
Also assume that $\delta=|u_--u_+|$ is sufficiently small.
Then there is a positive constant $\epsilon_0$ such that if  $\|u_0-\bar{u}_2(\cdot)-\tilde{u}_4(\cdot,0)\|_2+\delta \le \epsilon_0$, the problem \eqref{eq1.1}-\eqref{eq1.3} admits a unique solution $(u(x,t),q(x,t))$, which satisfies
\begin{equation*}
\begin{aligned}
 &u-\bar{u}_2-u_4^R\in C^0([0,\infty);H^2),\ \ \partial_x(u-\bar{u}_{2}-u_4^R) \in L^2(0,\infty;H^1),\\[2mm]
 &q-\bar{q}_2+\partial_xu_4^R\in C^0([0,\infty);H^3)\cap L^2(0,\infty;H^3),
\end{aligned}
\end{equation*}
and the asymptotic behavior
\begin{equation*}
\begin{split}
 &\sup_{x\in \mathbbm{R}_+} \left|\partial_x^k\left(u(x,t)-\bar{u}_2(x)-u_4^R(x) \right)\right|\rightarrow 0 \ \  \text{as} \ \ t\rightarrow \infty, \ \ k=0,1,\\[2mm]
 &\sup_{x\in \mathbbm{R}_+} \left|\partial_x^k\left(q(x,t)-\bar{q}_2(x)+\partial_xu_4^R(x)\right)\right|\rightarrow 0 \ \  \text{as} \ \ t\rightarrow \infty, \ \ k=0,1,2.
\end{split}
\end{equation*}
\end{thm}

\section{Asymptotics to Rarefaction Wave} \label{sec-3}
\subsection{Reformulation of the Problem in the Case of $u_+>u_-\ge 0$}

The special case: $0=u_-<u_+$ has been considered by Ruan and Zhu in \cite{Ruan2},
so we will focus on the case of $0<u_-<u_+$.
The case $u_->0$ means that the fluid blows in through the boundary $x=0$.
Hence, this initial boundary problem is called the \emph{in-flow} problem.
It is worth noticing that the boundary condition $u(0,t)=u_-$ is necessary for the well-posedness of the problem since the characteristic speed of the first hyperbolic equation \eqref{eq1.1}$_1$ is positive at boundary $x=0$.
Moreover, for the second elliptic equation \eqref{eq1.1}$_2$, we need boundary condition on $q(0,t)$ to ensure the well-posedness of the problem \eqref{eq1.1}.
From Lemma \ref{lemxsb} $\mathrm{(ii)}$ with $k=1$, we note that the boundary value of $q(x,t)$ can be defined as $q(0,t)=0$.
Therefore, in the case of $0<u_-<u_+$, the problem \eqref{eq1.1}-\eqref{eq1.3} is rewritten as
\begin{equation}\label{jzxsb} 
\begin{cases}
 u_t+(\frac{1}{2}u^2)_x+q_x=0, \ \ \ \ x\in \mathbbm{R}_+, \ \ \ \ t>0, \\[1mm]
 -q_{xx}+q+u_x=0, \ \ \ \ x\in \mathbbm{R}_+, \ \ \ \ t>0,\\[1mm]
 u(0,t)=u_-, \ \ \ \ q(0,t)=0, \ \ \ \ t\geq 0,\\[1mm]
 u(x,0)=u_0(x)=
   \begin{cases}
   =u_-, \ \ \ \ x=0,\\[2mm]
   \rightarrow u_+, \ \ \ \ x\rightarrow +\infty.
   \end{cases}
\end{cases}
\end{equation}
Set
\begin{equation}\label{xsbrdbl} 
    u(x,t)=\varphi(x,t)+w(x,t), \ \ \ \
    q(x,t)=\psi(x,t)+z(x,t).
\end{equation}
We note that $\psi=-\tilde{u}_x -\hat{q}= -\varphi_x- \hat{u}_x- \hat{q}$, so we can rewritten \eqref{xsbrdbl} as
\begin{equation*}
  \begin{cases}
    u(x,t)=\varphi(x,t)+w(x,t),\\[1mm]
    q(x,t)=z(x,t) -\varphi_x(x,t)- \hat{u}_x(x,t)- \hat{q}(x,t).
  \end{cases}
\end{equation*}
Then the perturbation $(w,z)$ satisfies
\begin{equation}\label{xsbrdfc} 
  \begin{cases}
    w_t+ww_x+(\varphi w)_x+z_x=R_1, \ \ \ \ x\in\mathbbm{R}_+, \ \ \ t>0,\\[1mm]
    -z_{xx}+z+w_x=R_2, \ \ \ \ x\in\mathbbm{R}_+, \ \ \ t>0,\\[1mm]
    w(0,t)=0, \ \ \ \ z(0,t)=0, \ \ \ t\geq 0, \\[1mm]
    w(x,0)=u_0(x)-\varphi_0(x)=u_0(x)-\tilde{u}_0(x)+\hat{u}(x,0), \ \ \ \ x\in \mathbbm{R}_+.
  \end{cases}
\end{equation}
We define the solution space as
\begin{equation*}
X_1(0,T) = \left\{
\begin{tabular}{c|c}
   \multirow{2}{*}{({\it w, z})}    & $w\in C^0([0,T);H^2)\cap C^1([0,T);H^1)$   \\[1mm]
             & $z\in C^0([0,T);H^3)\cap L^2(0,T;H^3)$ \ \
\end{tabular}
       \right\},
\end{equation*}
with $0<T\le +\infty$. Then the problem \eqref{xsbrdfc} can be solved globally in time as follows.

\begin{thm}\label{thmxsb} 
Suppose that the boundary condition and far field states satisfy $0<u_-<u_+$, the initial data $w_0\in H^2(\mathbbm{R}_+)$ and the wavelength $\delta=|u_--u_+|$ are sufficiently small.
Then there are the positive constants $\varepsilon_1$ and $C=C(\varepsilon_1)$ such that if $\|w_0\|_2 +\delta\le \varepsilon_1$,
the problem \eqref{xsbrdfc} admits a unique solution $(w(x,t),z(x,t))\in X_1(0,+\infty)$ satisfying
\begin{equation*}
\|w(t)\|_2^2+\|z(t)\|_3^2+\int_0^t \ (\|w_x(\tau)\|_1^2+\|z(\tau)\|_3^2 ) \,d{\tau}\le C(\|w_0\|_3^2+\delta),
\end{equation*}
and the asymptotic behavior
\begin{equation}\label{xsbdsjxw}
\begin{split}
 &\sup_{x\in \mathbbm{R}_+} |\partial_x^kw(x,t)|\rightarrow 0 \ \  \text{as} \ \ t\rightarrow \infty, \ \ k=0,1, \\
 &\sup_{x\in \mathbbm{R}_+} |\partial_x^kz(x,t)|\rightarrow 0 \ \  \text{as} \ \ t\rightarrow \infty, \ \ k=0,1,2.
\end{split}
\end{equation}
\end{thm}

The Combination of the following local existence and the \emph{a priori} estimates proves Theorem \ref{thmxsb}.

\begin{prop}[Local existence] \label{propxsbjbczx}
Suppose the boundary condition and far field states satisfy $u_-<u_+<0$, the initial data satisfy $w_0 \in H^2(\mathbbm{R}_+)$ and $\|w_0\|_{2}+ \delta\le\varepsilon_1$.
Then there are two positive constants $C=C(\varepsilon_1)$ and $T_0=T_0(\varepsilon_1)$ such that the problem \eqref{xsbrdfc} has a unique solution $(w,z)\in X_1(0,T_0)$, which satisfies
\begin{equation*}
\|w(t)\|_2^2+\|z(t)\|_3^2+\int_0^t \  ( \| w_x(\tau) \|_1^2+\| z(\tau) \|_3^2  )  \,d{\tau} \le C(\| w_0 \|_2^2+\delta),
\end{equation*}
for $t\in [0,T_0]$.
\end{prop}

\begin{prop}[\emph{A priori} estimates] \label{propxsbxygj}
Let $T$ be a positive constant. Suppose that the problem \eqref{xsbrdfc} has a unique solution $(w,z)\in X_1(0,T)$.
Then there exist two positive constants $\varepsilon_2(\le\varepsilon_1)$ and $C=C(\varepsilon_2)$ such that if
$\|w_0\|_2+\delta\le\varepsilon_2$, then we have the estimate
\begin{equation*}
\|w(t)\|_2^2+\|z(t)\|_3^2+\int_0^t \ (\| w_x(\tau)\|_1^2+\| z(\tau)\|_3^2 )  \,d{\tau}
\le C(\|w_0\|_2^2+\delta),
\end{equation*}
for $t\in[0,T]$.
\end{prop}

\subsection{\emph{A priori} Estimates}
Under the assumptions of Theorem \ref{thmxsb}, to give the proof of \emph{a priori} estimates in Proposition \ref{propxsbxygj},
we devote ourselves to the estimates on the solution $(w,z)\in X_1(0,T)$ (for some $T>0$) of \eqref{xsbrdfc} under the \emph{a priori} assumption
\begin{equation}\label{xsbxyjs} 
 |w_x(t)|_\infty\le\varepsilon_0, \ \ \ \ |w_t(t)|_\infty\le\varepsilon_0,
\end{equation}
where $0<\varepsilon_0\ll 1$. For simplicity, we divide the proof of the \emph{a priori} estimate into several lemmas.

\begin{lem}\label{lemxsbw}
There are the positive constants $\varepsilon_1(\le \varepsilon_0)$ and $C=C(\varepsilon_1)$ such that if $\|w_0\|_2+\delta \le \varepsilon_1$, then
\begin{equation}\label{xsbw}
   \|w(t)\|^2+\int_0^t \left(\|\sqrt{ \varphi_x}w(\tau)\|^2+\|z(\tau)\|_1^2\right) \,\mathrm{d}\tau \le C (\|w_0\|^2+\delta)
\end{equation}
holds for $t\in[0,T]$.
\end{lem}
{\it\bfseries Proof.}
Multiplying \eqref{xsbrdfc}$_1$ by $w$ and \eqref{xsbrdfc}$_2$ by $z$, and adding the two resulting equations up, we obtain
\begin{equation}\label{xsbjbnlgj} 
  \begin{aligned}
    \frac{1}{2}\frac{\mathrm{d}}{\mathrm{d}t}w^2+\frac{1}{2}\varphi_{x}w^2+z_x^2+z^2
    +\left\{\frac{1}{2} \varphi w^2+\frac{1}{3}w^3-z_xz+zw\right\}_x=R_1w+R_2z.
  \end{aligned}
\end{equation}
Integrating \eqref{xsbjbnlgj} over $\mathbbm{R}_+\times(0,t)$, using $w(0,t)=z(0,t)=0$, we get
\begin{equation}\label{xsbjbnlgjjf}
  \|w(t)\|^2+\int_0^t \left(\|\sqrt{\varphi_{x}}w(\tau)\|^2+\|z_x(\tau)\|^2+\|z(\tau)\|^2\right) \,\mathrm{d}\tau \le C \left( \|w_0\|^2+ \int_0^t\int_{\mathbbm{R}_+} |R_1w|+|R_2z| \,\mathrm{d}x \mathrm{d}\tau  \right).
\end{equation}
From Lemma \ref{lemjzxsb} $\mathrm{(v)}$ and $\mathrm{(vi)}$ with $k=l=0$, we have
\begin{equation}\label{xsbw1}
  \int_0^t\int_{\mathbbm{R}_+} |R_1w| \,\mathrm{d}x \mathrm{d}\tau\le \|w(t)\|\int_0^t \|R_1\| \,\mathrm{d}\tau  \le C \delta(1+\|w(t)\|^2)\int_0^t \mathrm{e}^{-c(1+t)} \,\mathrm{d}\tau \le C \delta(1+\|w(t)\|^2),
\end{equation}
and
\begin{equation}\label{xsbw2}
  \int_0^t\int_{\mathbbm{R}_+} |R_2 z| \,\mathrm{d}x \mathrm{d}\tau \le \frac{1}{4} \int_0^t \|z(\tau)\|^2 \,\mathrm{d}\tau+\int_0^t \|R_2\|^2 \,\mathrm{d}\tau  \le  \frac{1}{4} \int_0^t \|z(\tau)\|^2 \,\mathrm{d}\tau+C \delta (1+t)^{-\frac{3}{2} }.
\end{equation}
Substituting \eqref{xsbw1} and \eqref{xsbw2} into \eqref{xsbjbnlgjjf}, we conclude \eqref{xsbw}.
$\hfill\Box$

\begin{lem}\label{lemxsbwx}
There are two positive constants $\varepsilon_2(\le \varepsilon_1)$ and $C=C(\varepsilon_2)$
such that if $\|w_0\|_2+\delta \le \varepsilon_1$, then
\begin{equation}\label{xsbwx}
     \|w_x(t)\|^2+\int_0^t \left(\|\sqrt{\varphi_x}w_x(\tau)\|^2+\|z_{xx}(\tau)\|^2+\|z_x(\tau)\|^2\right) \,\mathrm{d}\tau \le C(\|w_0\|_1^2+\delta)
\end{equation}
holds for $t\in[0,T]$.
\end{lem}
{\it\bfseries Proof.}
We differentiate \eqref{xsbrdfc}$_1$ with respect to $x$ and multiply it by $w_x$,
and multiply \eqref{xsbrdfc}$_2$ by $-z_{xx}$. Then, adding these two equations up, we have
\begin{equation}\label{xsbwx1}
     \begin{aligned}
    \frac{1}{2}\frac{\mathrm{d}}{\mathrm{d}t}w_x^2+\frac{3}{2}\varphi_{x}w_x^2+z_{xx}^2+z_x^2
    +\left\{\frac{1}{2} \varphi w_x^2+\frac{1}{2}ww_x^2-z_xz\right\}_x=- \varphi_{xx}ww_x- \frac{1}{2} w_x^3+ R_{1x}w_x-R_2z_{xx}.
      \end{aligned}
\end{equation}
Integrating \eqref{xsbwx1} over $\mathbbm{R}_+\times(0,t)$, combining it with $\varphi(0,t)=u_->0$ and $w(0,t)=z(0,t)=0$, we get
\begin{equation}\label{xsbwxjf}
  \begin{aligned}[b]
    &\|w_x(t)\|^2+\int_0^t \left(\|\sqrt{\varphi_x}w_x(\tau)\|^2+\|z_{xx}(\tau)\|^2+\|z_x(\tau)\|^2\right) \,\mathrm{d}\tau \\[3mm]
    &\le C\left(\|w_0\|_1^2+u_-\int_0^t w_x^2(0,\tau) \,\mathrm{d}\tau+ \int_0^t \int_{\mathbbm{R}_+} |\varphi_{xx}ww_x|+ | w_x|^3+ |R_{1x}w_x|+|R_2z_{xx}| \,\mathrm{d}x \mathrm{d}\tau  \right).
  \end{aligned}
\end{equation}
Since the equation \eqref{xsbrdfc}$_1$ implies
\begin{equation*}
  u_-w_x(0,t)=-z_x(0,t)+R_1(0,t),
\end{equation*}
we can estimate the integral on the boundary as follows:
\begin{equation}\label{xsbwxbj}
  \begin{aligned}[b]
    u_-\int_0^t w_x^2(0,\tau) \,\mathrm{d}\tau
    &\le C\int_0^t \left(z_x^2(0,\tau)+R_1^2(0,\tau)\right) \,\mathrm{d}\tau\\[3mm]
    &\le C\int_0^t \left(|z_x|_\infty^2+|R_1|_\infty^2\right) \,\mathrm{d}\tau \\[3mm]
    &\le \frac{1}{8} \int_0^t \|z_{xx}(\tau)\|^2 \,\mathrm{d}\tau+ C \int_0^t \|z_x(\tau)\|^2 \,\mathrm{d}\tau+C \delta.
  \end{aligned}
\end{equation}
Next, we estimate the last four terms on the right-hand side of \eqref{xsbwxjf}. Using Lemma \ref{lemjzxsb} $\mathrm{(iv)}$ with $k=2$ and $l=0$, we have
\begin{equation*}
  \begin{aligned}
    \int_0^t \int_{\mathbbm{R}_+} |\varphi_{xx}ww_x| \,\mathrm{d}x \mathrm{d}\tau
    &\le \frac{1}{24} \int_0^t \|w_x(\tau)\|^2 \,\mathrm{d}\tau+ C|w(t)|_\infty^2 \int_0^t \|\varphi_{xx}(\tau)\|^2 \,\mathrm{d}\tau \\[3mm]
    &\le  \frac{1}{24} \int_0^t \|w_x(\tau)\|^2 \,\mathrm{d}\tau+ C \delta (1+t)^{-\frac{1}{2} },
  \end{aligned}
\end{equation*}
and
\begin{equation*}
  \begin{aligned}
    \int_0^t \int_{\mathbbm{R}_+}  | w_x(\tau)|^3 \,\mathrm{d}x \mathrm{d}\tau \le \varepsilon_0 \int_0^t \|w_x(\tau)\|^2 \,\mathrm{d}\tau.
  \end{aligned}
\end{equation*}
From Lemma \ref{lemjzxsb} $\mathrm{(v)}$ with $k=1$ and $\mathrm{(vi)}$ with $k=0, l=0$, we get
\begin{equation*}
  \begin{aligned}
    \int_0^t \int_{\mathbbm{R}_+}  |R_{1x}w_x| \,\mathrm{d}x \mathrm{d}\tau \le \frac{1}{24} \int_0^t \|w_x(\tau)\|^2 \,\mathrm{d}\tau+ C \delta ,
  \end{aligned}
\end{equation*}
and
\begin{equation}\label{xsbwx2}
  \begin{aligned}
    \int_0^t \int_{\mathbbm{R}_+} |R_2z_{xx}| \,\mathrm{d}x \mathrm{d}\tau \le \frac{1}{8}\int_0^t \|z_{xx}(\tau)\|^2 \,\mathrm{d}\tau +C \delta (1+t)^{-\frac{3}{2} }.
  \end{aligned}
\end{equation}

On the other hand, from \eqref{xsbrdfc}$_2$, we have
\begin{equation}\label{xsbwz}
  \|w_x(t)\|^2\le 3(\|z_{xx}(t)\|^2+\|z(t)\|^2+\|R_2(t)\|^2).
\end{equation}
In deriving the equation $\eqref{xsbwz}$ we have used the fact that for any $a,b,c\in \mathbbm{R}$,
$$(a+b+c)^2=a^2+b^2+c^2+2ab+2bc+2ac\le 3(a^2+b^2+c^2).$$
Substituting \eqref{xsbwxbj}-\eqref{xsbwx2} into \eqref{xsbwxjf} and using \eqref{xsbwz}, for some small $\delta$ and $\varepsilon_0(<\frac{1}{6} )$, we have
\begin{equation*}
  \begin{aligned}
    \|w_x(t)\|^2+\int_0^t \left(\|\sqrt{\varphi_x}w_x(\tau)\|^2+\|z_{xx}(\tau)\|^2+\|z_x(\tau)\|^2\right) \,\mathrm{d}\tau \le C(\|w_0\|_1^2+\delta).
  \end{aligned}
\end{equation*}
This completes the proof of Lemma \ref{lemxsbwx}.
$\hfill\Box$

For \eqref{xsbwxbj} and \eqref{xsbwz}, combining the results of Lemma \ref{lemxsbw} and \ref{lemxsbwx}, we can easily show the following Corollary \ref{corxsbwx} and Corollary \ref{corxsbwxbj}.

\begin{cor}\label{corxsbwx}
Under the assumptions of Lemma \ref{lemxsbwx}, there exists a positive constant $C$ such that
\begin{equation*}
    \int_0^t \|w_x(\tau)\|^2 \,\mathrm{d}\tau \le C(\|w_0\|_1^2+\delta), \quad \forall t\in[0,T].
\end{equation*}
\end{cor}

\begin{cor}\label{corxsbwxbj}
Under the assumptions of Lemma \ref{lemxsbwx}, there exists a positive constant $C$ such that
\begin{equation*}
    \int_0^t w_x^2(0,\tau) \,\mathrm{d}\tau \le C(\|w_0\|_1^2+\delta),  \quad \forall t\in[0,T].
\end{equation*}
\end{cor}

Next, we try to give the estimate for $w_{xx}$. When estimating $w_{xx}$, we need to deal with the boundary term $w_xx(0,t)$ (see $\eqref{xsbwxxjf}$). It is quite difficult to estimate the boundary term $w_xx(0,t)$ directly. However, we can get the estimate of the boundary term $w_{xt}(0,t)$ owing to $w(0,t)=w_t(0,t)=w_{tt}(0,t)=0$, and then the estimate of $wxx(0,t)$ is obtained through the equation $\eqref{xsbwxxzxxbj}$. Thus, to give the estimate for $w_{xx}$, we firstly proceed to the \emph{a priori} estimate for the derivatives $w_t$ and $w_{xt}$.

\begin{lem}\label{xsbwt1}
Under the assumptions of Lemma \ref{lemxsbwx}, there is a positive constant $C$ such that
\begin{equation*}
  \int_0^t \|w_t(\tau)\|^2 \,\mathrm{d}\tau \le C(\|w_0\|_1^2+\delta).
\end{equation*}
\end{lem}
{\it\bfseries Proof.}
With the help of Lemma \ref{lemxsbw}, Corollary \ref{corxsbwxbj} and Lemma \ref{lemjzxsb} $\mathrm{(v)}$ with $k=l=0$,
we see from \eqref{xsbrdfc}$_1$ that
\begin{equation*}
\begin{aligned}
    \int_0^t \|w_t(\tau)\|^2 \,\mathrm{d}\tau
  &\le C \int_0^t\int_{\mathbbm{R}_+} (w^2w_x^2+\varphi_x^2w^2+\varphi^2w_x^2+z_x^2+R_1^2) \,\mathrm{d}x \mathrm{d}\tau \\[3mm]
  &\le C \left( (|w(t)|_\infty^2  \!+\! |\varphi(t)|_\infty^2 )\int_0^t \|w_x(\tau)\|^2 \,\mathrm{d}\tau  \!+\! |\varphi_x(t)|_\infty \int_0^t \|\sqrt{\varphi_x}w(\tau)\|^2 \,\mathrm{d}\tau \!+\! \int_0^t \|z_x(\tau)\|^2 \,\mathrm{d}\tau \!+\!  \delta \right)\\[3mm]
  &\le C(\|w_0\|_1^2+\delta).
\end{aligned}
\end{equation*}
Thus, the proof of Lemma \ref{xsbwt1} is completed.
$\hfill\Box$

\begin{lem}\label{lemxsbwt2}
Under the assumptions of Lemma \ref{lemxsbwx}, there is a positive constant $C$ such that
\begin{equation}\label{xsbwt5}
    \|w_t(t)\|^2+\int_0^t \left(\|\varphi_xw_t(\tau)\|^2+\|z_{xt}(\tau)\|^2+\|z_t(\tau)\|^2\right) \,\mathrm{d}\tau \le C(\|w_0\|_1^2+\delta)
\end{equation}
holds for $t\in[0,T]$.
\end{lem}
{\it\bfseries Proof.}
We differentiate \eqref{xsbrdfc} with respect to $t$ and multiply the first and the second resulting equations by $w_t$ and $z_t$ respectively. Then, adding these two equations up, we have
\begin{equation}\label{xsbwt2}
     \begin{aligned}
    \frac{1}{2}\frac{\mathrm{d}}{\mathrm{d}t}&w_t^2+\frac{1}{2}\varphi_{x}w_t^2+z_{xt}^2+z_t^2
    +\left\{\frac{1}{2} \varphi w_t^2+\frac{1}{2}ww_t^2-z_{xt}z_t+z_tw_t\right\}_x\\[2mm]
    &=- \varphi_{xt}ww_t- \varphi_tw_xw_t- \frac{1}{2} w_xw_t^2+ R_{1t}w_t+R_{2t}z_{t}.
      \end{aligned}
\end{equation}
We note that $w_t(0,t)=z_t(0,t)=0$ due to $w(0,t)=z(0,t)=0$. Integrating \eqref{xsbwt2} over $\mathbbm{R}_+\times(0,t)$, we have
\begin{equation}\label{xsbwt2jf}
   \begin{aligned}
    \|w_t(t)\|^2&+\int_0^t \left(\|\varphi_xw_t(\tau)\|^2+\|z_{xt}(\tau)\|^2+\|z_t(\tau)\|^2\right) \,\mathrm{d}\tau \\[2mm]
     &\le \int_0^t\int_{\mathbbm{R}_+} (|\varphi_{xt}ww_t|+ |\varphi_tw_xw_t|+ |w_xw_t^2|+ |R_{1t}w_t|+|R_{2t}z_{t}|) \,\mathrm{d}x \mathrm{d}\tau.
   \end{aligned}
\end{equation}
Combining Lemma \ref{xsbwt1} and Lemma \ref{lemjzxsb}, using Cauchy-Schwarz inequality and \eqref{xsbxyjs}, we can estimate the terms on the right-hand side of \eqref{xsbwt2jf} as follows:
\begin{equation}\label{xsbwt3}
  \begin{aligned}
    \int_0^t\int_{\mathbbm{R}_+} |\varphi_{xt}ww_t| \,\mathrm{d}x \mathrm{d}\tau
    \le \int_0^t \|w_t(\tau)\|^2 \,\mathrm{d}\tau +C \delta|w(t)|_\infty^2 (1+t)^{-\frac{1}{2} },
  \end{aligned}
\end{equation}
\begin{equation*}
  \begin{aligned}
    \int_0^t\int_{\mathbbm{R}_+} |\varphi_tw_xw_t| \,\mathrm{d}x \mathrm{d}\tau
    \le \int_0^t \|w_t(\tau)\|^2 \,\mathrm{d}\tau+|\varphi_t(t)|_\infty^2\int_0^t \|w_x(\tau)\|^2 \,\mathrm{d}\tau,
  \end{aligned}
\end{equation*}
\begin{equation*}
  \begin{aligned}
    \int_0^t\int_{\mathbbm{R}_+} |w_xw_t^2| \,\mathrm{d}x \mathrm{d}\tau
    \le \varepsilon_0\int_0^t \|w_t(\tau)\|^2 \,\mathrm{d}\tau,
  \end{aligned}
\end{equation*}
\begin{equation*}
  \begin{aligned}
    \int_0^t\int_{\mathbbm{R}_+}  |R_{1t}w_t| \,\mathrm{d}x \mathrm{d}\tau
    \le \int_0^t \|w_t(\tau)\|^2 \,\mathrm{d}\tau+ C \delta,
  \end{aligned}
\end{equation*}
and
\begin{equation}\label{xsbwt4}
  \begin{aligned}
    \int_0^t\int_{\mathbbm{R}_+} |R_{2t}z_{t}| \,\mathrm{d}x \mathrm{d}\tau
    \le \frac{1}{4} \int_0^t \|z_t(\tau)\|^2 \,\mathrm{d}\tau+ C \delta(1+t)^{-\frac{5}{2} }.
  \end{aligned}
\end{equation}
Substituting \eqref{xsbwt3}-\eqref{xsbwt4} into \eqref{xsbwt2jf} yields \eqref{xsbwt5}.
This proves Lemma \ref{lemxsbwt2}.
$\hfill\Box$

Next, we show the estimate for $w_{xt}$ in the following Lemma \ref{lemxsbwxt}.

\begin{lem}\label{lemxsbwxt}
There are two positive constants $\varepsilon_3(\le \varepsilon_2)$ and $C=C(\varepsilon_3)$
such that if $\|w_0\|_2+\delta \le \varepsilon_3$, then
\begin{equation}\label{xsbwxt}
     \|w_{xt}(t)\|^2+\int_0^t \left(\|\sqrt{\varphi_x}w_{xt}(\tau)\|^2+\|z_{xxt}(\tau)\|^2+\|z_{xt}(\tau)\|^2\right) \,\mathrm{d}\tau \le C(\|w_0\|_1^2+\delta)+ C(\delta+\varepsilon_0) \int_0^t \|w_{xx}(\tau)\|^2 \,\mathrm{d}\tau
\end{equation}
holds for $t\in[0,T]$.
\end{lem}
{\it\bfseries Proof.}
Differentiate \eqref{xsbrdfc}$_1$ with respect to $x$ and $t$, then multiply it by $w_{xt}$.
Differentiate \eqref{xsbrdfc}$_2$ with respect to $t$ and multiply it by $-z_{xxt}$.
Finally, adding these two equations up, we have
\begin{equation}\label{xsbwxt1}
\begin{aligned}
       &\frac{1}{2}\frac{\mathrm{d}}{\mathrm{d}t}w_{xt}^2+\frac{3}{2}\varphi_{x}w_{xt}^2+z_{xxt}^2+z_{xt}^2
    +\left\{\frac{1}{2} \varphi w_{xt}^2+\frac{1}{2}ww_{xt}^2-z_{xt}z_t\right\}_x\\
    =&- \varphi_{xxt}ww_{xt}- \varphi_{xx}w_tw_{xt}- 2\varphi_{xt}w_xw_{xt}- \varphi_{t}w_{xx}w_{xt} -w_tw_{xx}w_{xt} -\frac{3}{2} w_xw_{xt}^2+ R_{1xt}w_{xt}-R_{2t}z_{xxt}.
    \end{aligned}
\end{equation}
Integrating \eqref{xsbwxt1} over $\mathbbm{R}_+\times(0,t)$, using $\varphi(0,t)=u_->0$ and $w(0,t)=z_t(0,t)=0$, we have
\begin{equation}\label{xsbwxtjif}
    \begin{aligned}[b]
      \|w_{xt}\|^2&+\int_0^t \left(\|\sqrt{\varphi_x}w_{xt}(\tau)\|^2+\|z_{xxt}(\tau)\|^2+\|z_{xt}(\tau)\|^2\right) \,\mathrm{d}\tau \\[2mm]
      &\le C\left(\|w_0\|_1^2 +\int_0^t w_{xt}^2(0,\tau) \,\mathrm{d}\tau +\int_0^t\int_{\mathbbm{R}_+} (|\varphi_{xxt}ww_{xt}|+ |\varphi_{xx}w_tw_{xt}|+ |\varphi_{xt}w_xw_{xt}|) \,\mathrm{d}x \mathrm{d}\tau  \right.\\[2mm]
      &\left. \ \ \ \ \ \ \ \ +\int_0^t\int_{\mathbbm{R}_+} (|\varphi_{t}w_{xx}w_{xt}|+ |w_tw_{xx}w_{xt}| + |w_xw_{xt}^2|+ |R_{1xt}w_{xt}|+|R_{2t}z_{xxt}|) \,\mathrm{d}x \mathrm{d}\tau \right).
    \end{aligned}
\end{equation}

Firstly, from \eqref{xsbrdfc}$_1$, we have the following equation at the boundary $x=0$,
\begin{equation}\label{xsbwxtzxtbj}
    u_-w_{xt}(0,t)=R_{1t}(0,t)-z_{xt}(0,t),
\end{equation}
which plays an important role in estimating boundary terms.
In fact, we differential \eqref{xsbrdfc}$_1$ with respect to $t$, then we get
\begin{equation}\label{xsbwxtzxt1}
      w_{tt}+w_tw_x+ww_{xt}+\varphi_{xt}w+\varphi_tw_x+\varphi w_{xt}+z_{xt}=R_{1t}.
\end{equation}
For $\varphi(0,t)=u_-$ and $w(0,t)=0$, the boundary values at $x=0$ of \eqref{xsbwxtzxt1} yields \eqref{xsbwxtzxtbj}.
According to \eqref{xsbwxtzxtbj}, using Lemma \ref{lemxsbwt2} and Lemma \ref{lemjzxsb} $\mathrm{(v)}$ with $k=0,l=1$, we can get
\begin{equation*}
      \begin{aligned}[b]
        \int_0^t w_{xt}^2(0,\tau) \,\mathrm{d}\tau
        &\le C \int_0^t (R_{1t}^2(0,\tau)+z_{xt}^2(0,\tau)) \,\mathrm{d}\tau \\[2mm]
        &\le C \int_0^t (|R_{1t}|_\infty^2+|z_{xt}|_\infty^2) \,\mathrm{d}\tau \\[2mm]
        &\le \frac{1}{8}\int_0^t \|z_{xxt}(\tau)\|^2 \,\mathrm{d}\tau +C \int_0^t \|z_{xt}(\tau)\|^2 \,\mathrm{d}\tau +C \delta \\[2mm]
        &\le \frac{1}{8}\int_0^t \|z_{xxt}(\tau)\|^2 \,\mathrm{d}\tau+ C(\|w_0\|_1^2+\delta).
      \end{aligned}
\end{equation*}

Next, the terms on the right-hand side of \eqref{xsbwxtjif} can be estimated as follows:
\begin{equation*}
      \begin{aligned}[b]
        &\int_0^t\int_{\mathbbm{R}_+} (|\varphi_{xxt}ww_{xt}|+|\varphi_{xx}w_tw_{xt}|+|\varphi_{xt}w_xw_{xt}|) \,\mathrm{d}x \mathrm{d}\tau\\[2mm]
        \le& \frac{1}{24} \int_0^t \|w_{xt}(\tau)\|^2 \,\mathrm{d}\tau +C \delta|w|_\infty^2 (1+t)^{-\frac{3}{2} }+C \delta(|w_t|_\infty^2 +|w_x|_\infty^2 ) (1+t)^{-\frac{1}{2} }.
      \end{aligned}
\end{equation*}
Using the \emph{a priori} assumption \eqref{xsbxyjs} and Cauchy-Schwarz inequality, we get
\begin{equation*}
      \begin{aligned}[b]
        &\int_0^t\int_{\mathbbm{R}_+} (|\varphi_{t}w_{xx}w_{xt}|+ |w_tw_{xx}w_{xt}|) \,\mathrm{d}x \mathrm{d}\tau
        \\
\le& |\varphi_t(t)|_\infty\int_0^t \|w_{xx}(\tau)\|\|w_{xt}(\tau)\| \,\mathrm{d}\tau+|w_t(t)|_\infty\int_0^t \|w_{xx}(\tau)\|\|w_{xt}(\tau)\| \,\mathrm{d}\tau\\[2mm]
\le& C(\delta+\varepsilon_0)\int_0^t \|w_{xx}(\tau)\|^2 \,\mathrm{d}\tau + ( \delta+ \varepsilon_0)\int_0^t \|w_{xt}(\tau)\|^2 \,\mathrm{d}\tau,
      \end{aligned}
\end{equation*}
and
\begin{equation*}
      \begin{aligned}
        \int_0^t\int_{\mathbbm{R}_+} |w_xw_{xt}^2| \,\mathrm{d}x \mathrm{d}\tau
        \le |w_x(t)|_\infty \int_0^t \|w_{xt}(\tau)\|^2 \,\mathrm{d}\tau
        \le \varepsilon_0\int_0^t \|w_{xt}(\tau)\|^2 \,\mathrm{d}\tau.
      \end{aligned}
\end{equation*}
From Lemma \ref{lemjzxsb}, the last two estimates can be given as follows:
\begin{equation*}
      \begin{aligned}
        \int_0^t\int_{\mathbbm{R}_+} |R_{1xt}w_{xt}| \,\mathrm{d}x \mathrm{d}\tau
        \le \frac{1}{24} \int_0^t \|w_{xt}(\tau)\|^2 \,\mathrm{d}\tau+ C \delta,
      \end{aligned}
\end{equation*}
and
\begin{equation}\label{xsbwxt2}
      \begin{aligned}
        \int_0^t\int_{\mathbbm{R}_+} |R_{2t}z_{xxt}| \,\mathrm{d}x \mathrm{d}\tau
        \le \frac{1}{8} \int_0^t \|z_{xxt}(\tau)\|^2 \,\mathrm{d}\tau +C \delta(1+t)^{-\frac{5}{2} }.
      \end{aligned}
\end{equation}
On the other hand, from \eqref{xsbrdfc}$_2$, we have
\begin{equation}\label{xsbwxtzxt}
  \|w_{xt}(t)\|^2\le 3(\|z_{xxt}(t)\|^2+\|z_t(t)\|^2+\|R_{2t}\|^2).
\end{equation}
Substituting \eqref{xsbwxtzxtbj}-\eqref{xsbwxt2} into \eqref{xsbwxtjif} and using \eqref{xsbwxtzxt}, for some small $\delta$ and $\varepsilon_0$, we have
\begin{equation*}
  \begin{aligned}
    \|w_x(t)\|^2+\int_0^t \left(\|\sqrt{\varphi_x}w_x(\tau)\|^2+\|z_{xx}(\tau)\|^2+\|z_x(\tau)\|^2\right) \,\mathrm{d}\tau \le C(\|w_0\|_1^2+\delta)+ C(\delta+\varepsilon_0) \int_0^t \|w_{xx}(\tau)\|^2 \,\mathrm{d}\tau.
  \end{aligned}
\end{equation*}
This completes the proof of Lemma \ref{lemxsbwxt}.
$\hfill\Box$

\begin{cor}\label{corwxtbj}
Under the same assumptions of Lemma \ref{lemxsbwxt}, there is a positive constant $C$ such that
\begin{equation*}
    \int_0^t w_{xt}^2(0,\tau) \,\mathrm{d}\tau \le C(\|w_0\|_1^2+\delta)+ C(\delta+\varepsilon_0) \int_0^t \|w_{xx}(\tau)\|^2 \,\mathrm{d}\tau,
    \quad \forall t\in[0,T].
\end{equation*}
\end{cor}

Finally, combining it with Lemma \ref{lemxsbwxt}, we show the estimate for $w_{xx}$.

\begin{lem}\label{lemxsbwxx}
There are two positive constants $\varepsilon_4(\le \varepsilon_3)$ and $C=C(\varepsilon_4)$
such that if $\|w_0\|_2+\delta \le \varepsilon_4$, then
\begin{equation}\label{xsbwxx}
     \|w_{xx}(t)\|^2+\int_0^t \left(\|\sqrt{\varphi_x}w_{xx}(\tau)\|^2+\|z_{xxx}(\tau)\|^2+\|z_{xx}(\tau)\|^2\right) \,\mathrm{d}\tau \le C(\|w_0\|_2^2+\delta)
\end{equation}
holds for $t\in[0,T]$.
\end{lem}
{\it\bfseries Proof.}
Differentiate \eqref{xsbrdfc}$_1$ twice with respect to $x$, then multiply it by $w_{xx}$.
Differentiate \eqref{xsbrdfc}$_2$ with respect to $x$ and multiply it by $-z_{xxx}$.
In the end, adding these two equations up, we have
\begin{equation}\label{xsbwxx1}
    \begin{aligned}[b]
       &\frac{1}{2}\frac{\mathrm{d}}{\mathrm{d}t}w_{xx}^2+\frac{5}{2}\varphi_{x}w_{xx}^2+z_{xxx}^2+z_{xx}^2
    +\left\{\frac{1}{2} \varphi w_{xx}^2+\frac{1}{2}ww_{xx}^2-z_{xx}z_x\right\}_x\\[1mm]
    =&- \varphi_{xxx}ww_{xx}- 3\varphi_{xx}w_xw_{xx} - \frac{5}{2}  w_xw_{xx}^2 + R_{1xx}w_{xx}-R_{2x}z_{xxx}.
    \end{aligned}
\end{equation}
Integrating \eqref{xsbwxx1} over $\mathbbm{R}_+\times(0,t)$, using $\varphi(0,t)=u_->0$ and $w(0,t)=0$, we obtain
\begin{equation}\label{xsbwxxjf}
  \begin{aligned}
    \|w_{xx}(t)\|^2&+\int_0^t \left(\|\sqrt{\varphi_x}w_{xx}(\tau)\|^2+\|z_{xxx}(\tau)\|^2+\|z_{xx}(\tau)\|^2\right) \,\mathrm{d}\tau \\[2mm]
    &\le C\left(\|w_0\|_2^2+\int_0^t w_{xx}^2(0,\tau) \,\mathrm{d}\tau+\int_0^t |z_{xx}(0,t)||z_x(0,t)| \,\mathrm{d}\tau  \right.\\[2mm]
    &\left. \ \ \ \ \ \ \ \  +\int_0^t\int_{\mathbbm{R}_+} (|\varphi_{xxx}ww_{xx}|+ |\varphi_{xx}w_xw_{xx}| +  |w_xw_{xx}^2| + |R_{1xx}w_{xx}|+|R_{2x}z_{xxx}|) \,\mathrm{d}x \mathrm{d}\tau \right).
  \end{aligned}
\end{equation}
Firstly, from \eqref{xsbrdfc}$_1$, we have an equation at the boundary $x=0$,
\begin{equation}\label{xsbwxxzxxbj}
    u_-w_{xx}(0,t)=R_{1x}(0,t)-z_{xx}(0,t)-w_{xt}(0,t)-w_x^2(0,t)-2 \varphi_x(0,t)w_x(0,t),
\end{equation}
which is important to estimate the boundary terms. In fact, we differential \eqref{xsbrdfc}$_1$ with respect to $x$, then we get
\begin{equation}\label{xsbwxxzxx1}
      w_{xt}+w_x^2+ww_{xx}+\varphi_{xx}w+2\varphi_xw_x+\varphi w_{xx}+z_{xx}=R_{1x}.
\end{equation}
For $\varphi(0,t)=u_-$ and $w(0,t)=0$, the boundary value at $x=0$ of \eqref{xsbwxxzxx1} yields \eqref{xsbwxxzxxbj}.
Using \eqref{xsbwxxzxxbj} and Cauchy-Schwarz inequality, combining Corollary \ref{corxsbwxbj} and Corollary \ref{corwxtbj}, we have
\begin{equation}\label{xsbwxxbjjf}
  \begin{aligned}[b]
    C \int_0^t w_{xx}^2(0,\tau) \,\mathrm{d}\tau
    &\le C \int_0^t \left(R_{1x}^2(0,\tau)+z_{xx}^2(0,\tau)+w_{xt}^2(0,\tau)+w_x^4(0,\tau)+ \varphi_x^2(0,\tau)w_x^2(0,\tau) \right) \,\mathrm{d}\tau \\[2mm]
    &\le C(\|w_0\|_1^2+\delta)+ C(\delta+\varepsilon_0) \int_0^t \|w_{xx}(\tau)\|^2 \,\mathrm{d}\tau++\frac{1}{8} \int_0^t \|z_{xxx}(\tau)\|^2 \,\mathrm{d}\tau.
  \end{aligned}
\end{equation}
Combining the results of Lemmas \ref{lemxsbw} and \ref{lemxsbwx}, we can get
\begin{equation*}
	\begin{aligned}
		\int_0^t |z_{xx}(0,\tau)||z_x(0,\tau)| \,\mathrm{d}\tau
		\le \frac{1}{8}\int_0^t \|z_{xxx}(\tau)\|^2 \,\mathrm{d}\tau  +C \int_0^t \left(\|z_{xx}(\tau)\|^2+\|z_x(\tau)\|^2\right) \,\mathrm{d}\tau.
	\end{aligned}
\end{equation*}

Next, the rest terms on the right-hand side of \eqref{xsbwxxjf} can be estimated as follows. According to the \emph{a priori} assumption \eqref{xsbxyjs},
\begin{equation*}
  \begin{aligned}
    \int_0^t\int_{\mathbbm{R}_+}  |w_xw_{xx}^2|  \,\mathrm{d}x \mathrm{d}\tau
    \le \varepsilon_0 \int_0^t \|w_{xx}(\tau)\|^2 \,\mathrm{d}\tau.
  \end{aligned}
\end{equation*}
From Lemma \ref{lemjzxsb}, using Cauchy-Schwarz inequality, we have
\begin{equation*}
  \begin{aligned}
    \int_0^t\int_{\mathbbm{R}_+} |\varphi_{xxx}ww_{xx}| \,\mathrm{d}x \mathrm{d}\tau
    \le \frac{1}{24} \int_0^t \|w_{xx}(\tau)\|^2 \,\mathrm{d}\tau + C \delta|w(t)|_\infty^2 (1+t)^{-\frac{3}{2} },
  \end{aligned}
\end{equation*}
\begin{equation*}
  \begin{aligned}
    \int_0^t\int_{\mathbbm{R}_+} |\varphi_{xx}w_xw_{xx}| \,\mathrm{d}x \mathrm{d}\tau
    \le \frac{1}{24} \int_0^t \|w_{xx}(\tau)\|^2 \,\mathrm{d}\tau + + C \delta|w_x(t)|_\infty^2 (1+t)^{-\frac{1}{2} },
  \end{aligned}
\end{equation*}
\begin{equation*}
  \begin{aligned}
    \int_0^t\int_{\mathbbm{R}_+}  |R_{1xx}w_{xx}| \,\mathrm{d}x \mathrm{d}\tau
    \le \frac{1}{24} \int_0^t \|w_{xx}(\tau)\|^2 \,\mathrm{d}\tau+C \delta,
  \end{aligned}
\end{equation*}
and
\begin{equation}\label{xsbwxx2}
  \begin{aligned}
    \int_0^t\int_{\mathbbm{R}_+} |R_{2x}z_{xxx}| \,\mathrm{d}x \mathrm{d}\tau
    \le \frac{1}{8} \int_0^t \|z_{xxx}(\tau)\|^2 \,\mathrm{d}\tau +C \delta(1+t)^{-\frac{5}{2} }.
  \end{aligned}
\end{equation}
In the end, from \eqref{xsbrdfc}$_2$, we have
\begin{equation}\label{xsbwxxzxx2}
  \|w_{xx}(t)\|^2\le 3\left(\|z_{xxx}(t)\|^2+\|z_x(t)\|^2+\|R_{2x}(t)\|^2\right).
\end{equation}
Substituting \eqref{xsbwxxbjjf}-\eqref{xsbwxx2} into \eqref{xsbwxxjf} and using \eqref{xsbwxxzxx2}, for some small $\delta$ and $\varepsilon_0$, we get
\begin{equation*}
  \begin{aligned}[b]
   & \|w_{xx}(t)\|^2+\int_0^t \left( \|\sqrt{\varphi_x}w_{xx}(\tau)\|^2+\|z_{xxx}(\tau)\|^2+\|z_{xx}(\tau)\|^2\right) \,\mathrm{d}\tau \\[2mm]
     \le & C(\|w_0\|_2^2+\delta)+\left(\frac{3}{4}+C \varepsilon_0+C \delta \right)\int_0^t \|z_{xxx}(\tau)\|^2 \,\mathrm{d}\tau,
  \end{aligned}
\end{equation*}
which yields \eqref{xsbwxx}.
$\hfill\Box$

Substituting \eqref{xsbwxx} into \eqref{xsbwxxbjjf}, using Lemma \ref{lemxsbwxt} and Lemma \ref{corwxtbj},
we can get the following Corollaries \ref{corxsbwxxbj}-\ref{corcorwxtbj}.

\begin{cor}\label{corxsbwxxbj}
Under the same assumptions of Lemma \ref{lemxsbwxx}, there exists a positive constant $C$ such that
\begin{equation*}
\begin{aligned}
    \int_0^t \|w_{xx}(\tau)\|^2 \,\mathrm{d}\tau\le C(\|w_0\|_2^2+\delta), \quad \forall t\in[0,T],\\
    \int_0^t w_{xx}^2(0,\tau) \,\mathrm{d}\tau\le C(\|w_0\|_2^2+\delta), \quad \forall t\in[0,T].
\end{aligned}
\end{equation*}	
\end{cor}

\begin{cor}\label{corlemxsbwxt}
Under the same assumptions of Lemma \ref{lemxsbwxx}, for some small $\varepsilon_0$, there exists a positive constant $C$ such that
\begin{equation*}
    \|w_{xt}(t)\|^2+\int_0^t \left(\|\sqrt{\varphi_x}w_{xt}(\tau)\|^2+\|z_{xxt}(\tau)\|^2+\|z_{xt}(\tau)\|^2\right) \,\mathrm{d}\tau \le C(\|w_0\|_2^2+\delta), \quad \forall t\in[0,T].
\end{equation*}
\end{cor}

\begin{cor}\label{corcorwxtbj}
Under the same assumptions of Lemma \ref{lemxsbwxx}, for some small $\varepsilon_0$, there exists a positive constant $C$ such that
\begin{equation*}
\begin{aligned}
    \int_0^t \|w_{xt}(\tau)\|^2 \,\mathrm{d}\tau \le C(\|w_0\|_2^2+\delta), \quad \forall t\in[0,T],\\
    \int_0^t w_{xt}^2(0,\tau) \,\mathrm{d}\tau \le C(\|w_0\|_2^2+\delta), \quad \forall t\in[0,T].	
\end{aligned}
\end{equation*}
\end{cor}

In the end, using the above estimates and the relation between $w$ and $z$, we can easily get the estimate for $\|z(t)\|_3$.

\begin{lem}\label{lemxsbz}
Under the assumptions of Lemma \ref{lemxsbwxx}, there holds
\begin{equation}\label{xsbzH3}
    \|z(t)\|_3^2\le C(\|w_0\|_2^2+\delta), \quad \forall t\in[0,T].
\end{equation}
\end{lem}
{\it\bfseries Proof.}
Firstly, from \eqref{xsbrdfc}$_1$, we have the following equation at the boundary $x=0$,
\begin{equation*}
    z_x(0,t)=-u_-(w_x(0,t)+R_1(0,t)),
\end{equation*}
which can help us to estimate the boundary terms. From \eqref{xsbrdfc}$_2$, we obtain
\begin{equation}\label{xsbzpf}
   z_{xx}^2+z^2+2z_x^2=w_x^2+R_2^2+2(z_xz)_x+2R_2w_x.
\end{equation}
Integrating \eqref{xsbzpf} over $\mathbbm{R}_+$, using Cauchy-Schwarz inequality and $z(0,t)=0$, we get
\begin{equation}\label{xsbzxx}
  \begin{aligned}[b]
     \|z_{xx}(t)\|^2+2\|z_x(t)\|^2+\|z(t)\|^2
     \le 2\|w_x(t)\|^2+2\|R_2\|^2
     \le 2\|w_x(t)\|^2+C\delta
    \le C(\|w_0\|_2^2+\delta).
  \end{aligned}
\end{equation}

Secondly, differentiating \eqref{xsbrdfc}$_2$ with respect to $x$ and integrating these equations over $\mathbbm{R}_+$, combining \eqref{xsbzxx}, Lemma \ref{lemxsb} and Lemma \ref{lemxsbwxx}, we get
\begin{equation}\label{xsbzxxx}
  \begin{aligned}
    \|z_{xxx}(t)\|^2
    \le C\left(\|z_x(t)\|^2+\|w_{xx}(t)\|^2+\|R_{2x}(t)\|^2\right)
    \le C\left(\|w_0\|_2^2+\delta\right).
  \end{aligned}
\end{equation}
Combining \eqref{xsbzxx} with \eqref{xsbzxxx}, we finish the proof of Lemma \ref{lemxsbz}.
$\hfill\Box$

\subsection{Asymptotic Behavior toward the Rarefaction Wave}
By combining the local existence, Proposition \ref{propxsbjbczx} and the \emph{a priori} estimates, we can get the global in time solution
\begin{equation*}
\begin{cases}
	w\in C^0([0,\infty);H^2)\cap C^1([0,\infty);H^1),\ \ w_t\in C^0([0,\infty);H^1)\cap L^2(0,\infty;H^1)\\[1mm]
	z\in C^0([0,\infty);H^3)\cap L^2(0,\infty;H^3),\ \ z_t\in L^2(0,\infty;H^2),
\end{cases}
\end{equation*}
such that
\begin{equation}
\begin{cases}\label{xsbxyfjjg}
	\sup\limits_{t\ge0} (\|w(t)\|_2^2+\|w_t(t)\|_1^2+\|z(t)\|_3^2)< \infty,\\[4mm]
	\int_0^\infty (\|w_x(t)\|_1^2+\|w_t(t)\|_1^2+\|z(t)\|_3^2+\|z_t(t)\|_2^2) \,\mathrm{d}t < \infty .
\end{cases}
\end{equation}
In order to show the large-time behavior \eqref{xsbdsjxw} in Theorem \ref{thmxsb}, using the Sobolev inequality
\begin{equation*}
	\begin{cases}
	\sup\limits_{x\in \mathbbm{R}_+}|f(x,t)|\le \sqrt{2}\|f(t)\|^{\frac{1}{2} }\|f_x(t)\|^{\frac{1}{2} }, \\[4mm]
	\sup\limits_{x\in \mathbbm{R}_+}|f_x(x,t)|\le \sqrt{2}\|f_x(t)\|^{\frac{1}{2} }\|f_{xx}(t)\|^{\frac{1}{2} },\\[4mm]
	\sup\limits_{x\in \mathbbm{R}_+}|f_{xx}(x,t)|\le \sqrt{2}\|f_{xx}(t)\|^{\frac{1}{2} }\|f_{xxx}(t)\|^{\frac{1}{2} },
	\end{cases}
\end{equation*}
we just need to prove
\begin{equation}\label{xsbdsjxw1}
 	\|w_x(t)\|\to 0,\ \ \|z(t)\|\to 0, \ \ \|z_{xx}(t)\|\to 0, \ \ \text{as} \ \ t\to \infty.
 \end{equation}
According to \eqref{xsbxyfjjg}, we only need to show
\begin{equation}\label{xsbdsjxw2}
	\int_0^\infty \left|\frac{\mathrm{d}}{\mathrm{d}t}\|w_x(t)\|^2\right| \,\mathrm{d}t <\infty, \ \ \ \ \int_0^\infty \left|\frac{\mathrm{d}}{\mathrm{d}t}\|z(t)\|^2\right|  \,\mathrm{d}t <\infty, \ \ \int_0^\infty \left|\frac{\mathrm{d}}{\mathrm{d}t}\|z_{xx}(t)\|^2\right| \,\mathrm{d}t <\infty.
\end{equation}
Here, we give the proof of \eqref{xsbdsjxw2} as follows.

Combining with the results of \eqref{xsbxyfjjg}, we can easily get
\begin{equation*}
\begin{aligned}
	\int_0^\infty \left|\frac{\mathrm{d}}{\mathrm{d}t}\|w_x(t)\|^2\right| \, \mathrm{d}t
	&\le 2\int_0^\infty \int_{\mathbbm{R}_+} |w_x||w_{xt}| \,\mathrm{d}x\mathrm{d}t \\[3mm]
	&\le \int_0^\infty \|w_x(t)\|^2 \,\mathrm{d}t  +\int_0^\infty \|w_{xt}(t)\|^2 \,\mathrm{d}t <\infty.	
\end{aligned}
\end{equation*}
Similarly, using \eqref{xsbxyfjjg}, we can obtain
\begin{equation*}
	\begin{aligned}
		\int_0^\infty \left|\frac{\mathrm{d}}{\mathrm{d}t}\|z(t)\|^2\right|  \,\mathrm{d}t
		&\le 2 \int_0^\infty \int_{\mathbbm{R}_+} |z(t)||z_t(t)| \,\mathrm{d}x \mathrm{d}t\\[3mm]
		&\le \int_0^\infty \|z(t)\|^2 \,\mathrm{d}t  +\int_0^\infty \|z_t(t)\|^2 \,\mathrm{d}t <\infty ,
	\end{aligned}
\end{equation*}
and
\begin{equation*}
	\begin{aligned}
		\int_0^\infty \left|\frac{\mathrm{d}}{\mathrm{d}t}\|z_{xx}(t)\|^2\right|  \,\mathrm{d}t
		&\le 2 \int_0^\infty \int_{\mathbbm{R}_+} |z_{xx}(t)||z_{xxt}(t)| \,\mathrm{d}x \mathrm{d}t\\[3mm]
		&\le \int_0^\infty \|z_{xx}(t)\|^2 \,\mathrm{d}t  +\int_0^\infty \|z_{xxt}(t)\|^2 \,\mathrm{d}t  <\infty.
	\end{aligned}
\end{equation*}
Therefore, we finish the proof of the large-time behavior \eqref{xsbdsjxw}. That is, the proof of Theorem \ref{thmxsb} is completed.

\section{Asymptotics to Stationary Solution} \label{sec-4}
\subsection{Reformulation of the Problem in the Case of $u_-<u_+\le 0$}
In the cases (1): $u_-<u_+<0$ and (2): $u_-<u_+=0$, the IBVP admits a stationary solution $(\bar{u},\bar{q})=(\bar{u}_i(x),\bar{q}_i(x)), i=1,2$, respectively. The stationary solution $(\bar{u},\bar{q})$ satisfies the following ordinary differential equations
\begin{equation*}
	\begin{cases}
		\bar{u}\bar{u}_x+\bar{q}_x=0, \ \ \ \ x\in \mathbbm{R}_+,\\[1mm]
		-\bar{q}_{xx}+\bar{q}+\bar{u}_x=0, \ \ \ \ x\in \mathbbm{R}_+,\\[1mm]
		\bar{u}(0)=u_-, \ \ \ \ \bar{u}(+\infty)=u_+, \ \ \ \ \bar{q}(+\infty)=0.
	\end{cases}
\end{equation*}
Put
\begin{equation*}
u(x,t)=\bar{u}(x)+w(x,t), \qquad q(x,t)=\bar{q}(x)+z(x,t).
\end{equation*}
The equation \eqref{eq1.1} can be reformulated as
\begin{equation}\label{eq2.8}
\begin{cases}
 w_t+ww_x+(\bar{u}w)_x+z_x=0, \ \ \ \ x\in \mathbbm{R}_+, \ \ \ t>0,\\[1mm]
 -z_{xx}+z+w_x=0, \ \ \ \ x\in \mathbbm{R}_+, \ \ \ t>0,\\[1mm]
 w(0,t)=0, \ \ \ \ w(+\infty,t)=0, \ \ \ t>0,\\[1mm]
 w(x,0)=w_0(x)=u_0(x)-\bar{u}(x), \ \ \ \ x\in \mathbbm{R}_+.
 \end{cases}
\end{equation}
Define the solution space of \eqref{eq2.8} by
$$X_2(0,T)=\left\{ w\in C^0([0,T);H^2), w_x \in L^2(0,T;H^1); z\in C^0([0,T);H^3)\cap L^2(0,T;H^3)\right\}$$
with $0<T\le +\infty$. Then the problem \eqref{eq2.8} can be solved globally in time as follows.

\begin{thm}\label{thm2.2}
Suppose that the boundary condition and far field states satisfy $u_-<u_+\le 0$, the initial data $w_0\in H^2(\mathbbm{R}_+)$ and the wavelength $\delta=|u_--u_+|$ in Lemma \ref{lem2.1.1} both are sufficiently small.
Then there are positive constants $\varepsilon_1$ and $C=C(\varepsilon_1)$ such that if $\|w_0\|_2 +\delta\le \varepsilon_1$, the problem \eqref{eq2.8} admits a unique solution $(w(x,t),z(x,t))\in X_2(0,+\infty)$ satisfying
\begin{equation*}
\|w(t)\|_2^2+\|z(t)\|_3^2+\int_0^t \ (\|w_x(\tau)\|_1^2+\|z(\tau)\|_3^2 ) \,d{\tau}\le C\|w_0\|_2^2,
\end{equation*}
and the asymptotic behavior
\begin{equation}
\begin{split}\label{wtjdsjxw}
 &\sup_{x\in \mathbbm{R}_+} |\partial_x^kw(x,t)|\rightarrow 0 \ \  \text{as} \ \ t\rightarrow \infty, \ \ \ k=0,1, \\[2mm]
 &\sup_{x\in \mathbbm{R}_+} |\partial_x^kz(x,t)|\rightarrow 0 \ \  \text{as} \ \ t\rightarrow \infty, \ \ \ k=0,1,2.
\end{split}
\end{equation}
\end{thm}

Theorem \ref{thm2.2} is proved by combining the local existence of the solution together with the \emph{a priori} estimates.

\begin{prop}[Local existence]
Suppose the boundary condition and far field states satisfy $u_-<u_+\le 0$, the initial data satisfies $w_0 \in H^2(\mathbbm{R}_+)$ and $\|w_0\|_{2}+ \delta\le\varepsilon_2$. Then there are two positive constants $C=C(\varepsilon_2)$ and $T_0=T_0(\varepsilon_2)$ such that the problem \eqref{eq2.8} has a unique solution $(w,z)\in X_2(0,T_0)$, which satisfies
\begin{equation*}
\|w(t)\|_2^2+\|z(t)\|_3^2+\int_0^t \  ( \| w_x(\tau) \|_1^2+\| z(\tau) \|_3^2  )  \,d{\tau} \le C\| w_0 \|_2^2, \quad \forall t\in [0,T_0].
\end{equation*}
\end{prop}

\begin{prop}[A priori estimates] \label{prop2.2}
Let $T$ be a positive constant. Suppose that the problem \eqref{eq2.8} has a unique solution $(w,z)\in X_2(0,T)$. Then there exist positive constants $\varepsilon_1(\le\varepsilon_2)$ and $C=C(\varepsilon_1)$ such that if $\|w_0\|_2+\delta\le\varepsilon_1$, then we have the estimate
\begin{equation*}
\|w(t)\|_2^2+\|z(t)\|_3^2+\int_0^t \ ( \| w_x(\tau)\|_1^2+\| z(\tau)\|_3^2 )  \,d{\tau}
\le C\|w_0\|_2^2, \quad \forall t\in [0,T_0].
\end{equation*}
\end{prop}

\subsection{\emph{A priori} Estimates}
Under the assumptions of Theorem \ref{thm2.2}, we want to give the proof of the \emph{a priori} estimate in Proposition \ref{prop2.2}.
To do this, we devote ourselves to the estimates on the solution $(w,z)\in X_2(0,T)$ (for some $T>0$) of \eqref{eq2.8} under the
\emph{a priori} assumption
\begin{equation}\label{eq2.37}
|w_x(t)|_\infty\le\varepsilon_0,
\end{equation}
where $0<\varepsilon_0\ll 1$. For simplicity, we divide the proof of the \emph{a priori} estimate into the following lemmas.

\begin{lem}\label{lem2.6}
There are positive constants $\varepsilon_1(\le \varepsilon_0)$ and $C=C(\varepsilon_1)$ such that if $\|w_0\|_2+\delta \le \varepsilon_1$, then
\begin{equation}\label{eq2.24}
\|w(t)\|_1^2+\int_0^t \ \left(\|\sqrt{\bar{u}_x}w(\tau) \|^2 + \|\sqrt{\bar{u}_x}w_x(\tau) \|^2+\|w_x(\tau)\|^2\right) \,\mathrm{d}{\tau}+\int_0^t w_x^2(0,\tau) \, \mathrm{d}{\tau}\le C\|w_0\|_1^2,
\end{equation}
for $t\in [0,T]$.
\end{lem}
{\it\bfseries Proof.}
Multiplying \eqref{eq2.8}$_1$ by $w$ and \eqref{eq2.8}$_2$ by $z$, and adding the two resulting equations up, we obtain
\begin{equation}\label{eq2.81}
\frac{1}{2}\frac{\mathrm{d}}{\mathrm{d}t}w^2+\frac{1}{2}\bar{u}_xw^2+z_x^2+z^2+\left\{\frac{1}{2}\bar{u}w^2+\frac{1}{3}w^3-z_xz+zw\right\}_x=0.
\end{equation}
We differentiate \eqref{eq2.8}$_1$ with respect to $x$ and multiply it by $w_x$, and multiply \eqref{eq2.8}$_2$ by $-z_{xx}$.
Adding these two equations up, we obtain
\begin{equation}\label{eq2.82}
\frac{1}{2}\frac{\mathrm{d}}{\mathrm{d}t}w_x^2+\frac{3}{2}\bar{u}_xw_x^2+z_{xx}^2+z_x^2+\left\{\frac{1}{2}\bar{u}w_{x}^2+\frac{1}{2}ww_x^2-z_{x}z\right\}_x=-\frac{1}{2}w_x^3-\bar{u}_{xx}ww_x.
\end{equation}
On the other hand, rewriting \eqref{eq2.8}$_2$ in the form $w_{x}=z_{xx}-z$ and squaring this equation, we get
\begin{equation}\label{eq2.83}
w_x^2=z_{xx}^2+z^2+2z_x^2-2(z_xz)_x.
\end{equation}
Adding \eqref{eq2.81}, \eqref{eq2.82} and \eqref{eq2.83} up, we get
\begin{equation}\label{eq2.84}
\begin{aligned}[b]
	\frac{1}{2}\frac{\mathrm{d}}{\mathrm{d}t}(w^2+w_x^2)&+\frac{1}{2}\bar{u}_xw^2+\frac{3}{2}\bar{u}_xw_x^2+w_x^2\\[3mm]
	&+\left\{\frac{1}{2}\bar{u}w^2+\frac{1}{3}w^3+zw+\frac{1}{2}\bar{u}w_{x}^2+\frac{1}{2}ww_x^2\right\}_x
	=-\frac{1}{2}w_x^3-\bar{u}_{xx}ww_x.
\end{aligned}
\end{equation}
Integrating \eqref{eq2.84} over $\mathbbm{R}_+\times(0,t)$, combining it with $u_-<0$ and $w(0,t)=0$, we have
\begin{equation}\label{wtjwxjf}
\begin{aligned}[b]
	\frac{1}{2} \|w(t)\|_1^2&+ \int_0^t \left(\frac{1}{2}\|\sqrt{\bar{u}_x}w(\tau)\|^2+\frac{3}{2} \|\sqrt{\bar{u}_x}w_x(\tau)\|^2+\|w_{x}(\tau)\|^2\right) \,\mathrm{d}\tau-u_-\int_0^t w_x^2(0,\tau) \,\mathrm{d}\tau\\[3mm]
  &\le \frac{1}{2} \|w_0\|_1^2 + \int_0^t\int_{\mathbbm{R}_+} (|w_x|^3+|\bar{u}_{xx}ww_x|) \,\mathrm{d}x \mathrm{d}\tau.
\end{aligned}
\end{equation}
From \eqref{eq2.37} and Lemma \ref{lem2.1.1}, we get
\begin{equation}\label{eq2.84.1}
	\int_0^t\int_{\mathbbm{R}_+} |w_x|^3 \,\mathrm{d}x \mathrm{d}\tau \le |w_x(t)|_\infty \int_0^t \|w_x(\tau)\|^2 \,\mathrm{d}\tau\le \varepsilon_0 \int_0^t \|w_x(\tau)\|^2 \,\mathrm{d}\tau,
\end{equation}
and
\begin{equation}\label{eq2.84.2}
\begin{aligned}[b]
  \int_0^t\int_{\mathbbm{R}_+} |\bar{u}_{xx}ww_x| \,\mathrm{d}x \mathrm{d}\tau
  &\le \left|\frac{\bar{u}^2_{xx}}{\bar{u}_x}\right|_\infty \int_0^t \|\sqrt{\bar{u}_x}w(\tau)\|^2 \,\mathrm{d}\tau +\frac{1}{8}\int_0^t \|w_x(\tau)\|^2 \,\mathrm{d}\tau \\[3mm]
  &\le C\delta \int_0^t \|\sqrt{\bar{u}_x}w(\tau)\|^2 \,\mathrm{d}\tau+\frac{1}{8}\int_0^t \|w_x(\tau)\|^2 \,\mathrm{d}\tau.
\end{aligned}
\end{equation}
Substituting \eqref{eq2.84.1}-\eqref{eq2.84.2} into \eqref{wtjwxjf}, for some small $\delta$ and $\varepsilon_0(<\frac{1}{8} )$, we get
\begin{equation}\label{wtjwx}
\begin{aligned}
  \|w(t)\|_1^2+\int_0^t \left(\|\sqrt{\bar{u}_x}w(\tau)\|^2+\|\sqrt{\bar{u}_x}w_x(\tau)\|^2+\|w_x(\tau)\|_2^2\right) \,\mathrm{d}\tau+\int_0^t w_x^2(0,\tau) \,\mathrm{d}\tau\le C\|w_0\|_1^2.
\end{aligned}
\end{equation}
Hence, Lemma \ref{lem2.6} is proved.
$\hfill\Box$

The combination of Lemma \ref{lem2.6} and equation \eqref{eq2.83} yields the following Lemma \ref{lemwtjz1}.

\begin{lem}\label{lemwtjz1}
Under the same assumptions of Lemma \ref{lem2.6}, there exists a positive constant C such that
\begin{equation*}
    \int_0^t \|z(\tau)\|_2^2 \,\mathrm{d}\tau \le C \|w_0\|_1^2
\end{equation*}
holds for $t\in [0,T]$.
\end{lem}
{\it\bfseries Proof.}
Integrating \eqref{eq2.83} over $\mathbbm{R}_+\times(0,t)$, we can easily get the following estimate for $z$,
\begin{equation}\label{wtjzH2}
  	\int_0^t (\|z_{xx}(\tau)\|^2+2\|z_x(\tau)\|^2+\|z(\tau)\|^2) \,\mathrm{d}\tau \le \int_0^t (\|w_x(\tau)\|^2+2|z_x(0,\tau)||z(0,\tau)|) \,\mathrm{d}\tau  .
\end{equation}
From \eqref{eq2.8}$_1$, we have the equation at the boundary $x=0$,
\begin{equation}\label{wtjwzbj}
    -u_-w_x(0,t)=z_x(0,t).
\end{equation}
Using \eqref{wtjwzbj} and Cauchy-Schwarz inequality, we have
\begin{equation}\label{wtjzxbj}
  \begin{aligned}[b]
    2 \int_0^t  |z_x(0,\tau)||z(0,\tau)|  \,\mathrm{d}\tau
    \le C \int_0^t |w_x(0,\tau)|^2 \,\mathrm{d}\tau +\frac{1}{2} \int_0^t |z|_\infty^2 \,\mathrm{d}\tau
    \le C\|w_0\|_1^2+\frac{1}{2} \int_0^t (\|z_x(\tau)\|^2+\|z(\tau)\|^2) \,\mathrm{d}\tau.
  \end{aligned}
\end{equation}
Substituting \eqref{wtjzxbj} into \eqref{wtjzH2}, we finish the proof of Lemma \ref{lemwtjz1}.
$\hfill\Box$

\begin{lem}\label{lem2.8}
There are positive constants $\varepsilon_2(\le \varepsilon_1)$ and $C=C(\varepsilon_2)$
such that if $\|w_0\|_2+\delta \le \varepsilon_2$, then
\begin{equation}\label{eq2.86}
\|w_{xx}(t)\|^2+\int_0^t \left(\|\sqrt{\bar{u}_x}w_{xx}(\tau)\|^2+\|z_{xxx}(\tau)\|^2+\|z_{xx}(\tau)\|^2+w_{xx}^2(0,\tau)\right) \,\mathrm{d}\tau \le C\|w_0\|_2^2,
\end{equation}
for $t\in [0,T]$.
\end{lem}
{\it\bfseries Proof.}
We differentiate \eqref{eq2.8} twice with respect to $x$ and multiply the first and the second resulting equations by $w_{xx}$ and $z_{xx}$ respectively. Then, adding these two equations up, we have
\begin{equation}\label{eq2.85.2}
    \begin{aligned}[b]
      \frac{1}{2}\frac{\mathrm{d}}{\mathrm{d}t}w_{xx}^2&+\frac{5}{2}\bar{u}_xw_{xx}^2+z_{xxx}^2+z_{xx}^2
      +\left\{\frac{1}{2}\bar{u}w_{xx}^2+\frac{1}{2}ww_{xx}^2+z_{xx}w_{xx}-z_{xx}z_{xxx}\right\}_x \\[3mm]
      &=-\bar{u}_{xxx}ww_{xx}-3\bar{u}_{xx}w_{x}w_{xx}-\frac{5}{2}w_xw_{xx}^2.
    \end{aligned}
\end{equation}
Integrating \eqref{eq2.85.2} over $\mathbbm{R}_+\times(0,t)$, using $u_-<0$ and $z_{xx}w_{xx}-z_{xx}z_{xxx}=-z_{xx}z_x$ due to \eqref{eq2.8}$_2$, we have after some calculations that
\begin{equation}\label{eq2.85.5}
    \begin{aligned}[b]
      \|w_{xx}(t)\|^2&+\int_0^t (\|\sqrt{\bar{u}_x}w_{xx}(\tau)\|+\|z_{xxx}(\tau)\|^2+\|z_{xx}(\tau)\|^2) \,\mathrm{d}\tau-u_-\int_0^t w_{xx}^2(0,\tau) \,\mathrm{d}\tau \\[3mm]
      &\le C\left(\|w_0\|_2^2+\int_0^t |z_{xx}(0,\tau)||z_{x}(0,\tau)| \,\mathrm{d}\tau \right. \\[3mm]
      &\left. \ \ \ \ \ \ \ \ \ \ \ \ \ \    +\int_0^t\int_{\mathbbm{R}_+} (|\bar{u}_{xxx}ww_{xx}|+|\bar{u}_{xx}w_xw_{xx}|+|w_x|w_{xx}^2) \,\mathrm{d}x \mathrm{d}\tau \right).
    \end{aligned}
\end{equation}
Firstly, from \eqref{eq2.24}, we have
\begin{equation}\label{eq2.85.6}
  \begin{aligned}[b]
    C \int_0^t |z_{xx}(0,\tau)||z_x(0,\tau)| \,\mathrm{d}\tau
    &\le C \int_0^t |z_{xx}|_\infty|z_x|_\infty \,\mathrm{d}\tau\\[3mm]
    &\le \frac{1}{4}  \int_0^t \|z_{xxx}(\tau)\|^2 \,\mathrm{d}\tau +C \int_0^t (\|z_{xx}(\tau)\|^2+\|z_x(\tau)\|^2) \,\mathrm{d}\tau \\[3mm]
    &\le \frac{1}{4}  \int_0^t \|z_{xxx}(\tau)\|^2 \,\mathrm{d}\tau +C \|w_0\|_1^2.
  \end{aligned}
\end{equation}
From Lemma \ref{lem2.1.1} and \eqref{eq2.37}, and using Cauchy-Schwarz inequality, we can get the estimates on the right-hand side of the equation \eqref{eq2.85.5}
\begin{equation}\label{eq2.85.7}
    \begin{aligned}[b]
      \int_0^t\int_{\mathbbm{R}_+} |\bar{u}_{xxx}ww_{xx}| \,\mathrm{d}x \mathrm{d}\tau
      &\le \frac{1}{8}\int_0^t \|w_{xx}(\tau)\|^2 \,\mathrm{d}\tau +\left|\frac{\bar{u}_{xxx}^2}{\bar{u}_x}\right|_\infty \int_0^t \|\sqrt{\bar{u}_x}w(\tau)\|^2 \,\mathrm{d}\tau  \\[3mm]
      &\le \frac{1}{8}\int_0^t \|w_{xx}(\tau)\|^2 \,\mathrm{d}\tau +C\delta \int_0^t \|\sqrt{\bar{u}_x}w(\tau)\|^2 \,\mathrm{d}\tau,
    \end{aligned}
\end{equation}
and
\begin{equation}\label{eq2.85.8}
    \begin{aligned}
      \int_0^t\int_{\mathbbm{R}_+} |\bar{u}_{xx}w_xw_{xx}| \,\mathrm{d}x \mathrm{d}\tau
      &\le  \frac{1}{8}\int_0^t \|w_{xx}(\tau)\|^2 \,\mathrm{d}\tau +|\bar{u}_{xx}|_\infty \int_0^t\int_{\mathbbm{R}_+} w_x^2 \,\mathrm{d}x \mathrm{d}\tau\\[3mm]
      &\le  \frac{1}{8}\int_0^t \|w_{xx}(\tau)\|^2 \,\mathrm{d}\tau +C \int_0^t \|w_x(\tau)\|^2 \,\mathrm{d}\tau,
    \end{aligned}
\end{equation}
\begin{equation}\label{eq2.85.9}
    \begin{aligned}
      \int_0^t\int_{\mathbbm{R}_+} |w_x|w_{xx}^2 \,\mathrm{d}x \mathrm{d}\tau
      \le |w_x|_\infty \int_0^t\int_{\mathbbm{R}_+} w_{xx}^2 \,\mathrm{d}x \mathrm{d}\tau
      \le \varepsilon_0 \int_0^t \|w_{xx}(\tau)\|^2 \, \mathrm{d}\tau.
    \end{aligned}
\end{equation}
In the end, from \eqref{eq2.8}$_2$, the relation between $w$ and $z$ satisfies
\begin{equation}\label{wtjwz2}
    \|w_{xx}(t)\|^2\le 2(\|z_{xxx}(t)\|^2+\|z_x(t)\|^2).
\end{equation}
Substituting \eqref{eq2.85.6}-\eqref{eq2.85.9} and \eqref{wtjwz2} into \eqref{eq2.85.5}, combining Lemma \ref{lem2.6} and Lemma \ref{lemwtjz1}, this yields \eqref{eq2.86} for some small $\delta$ and $\varepsilon_0(<\frac{1}{8})$.
$\hfill\Box$

Corollary \ref{corwtjwxx} is given by the combination of Lemma \ref{lem2.8} and \eqref{wtjwz2}.

\begin{cor}\label{corwtjwxx}
Under the same assumptions of Lemma \ref{lem2.8}, the estimate
\begin{equation*}
    \int_0^t \|w_{xx}(\tau)\|^2 \,\mathrm{d}\tau \le C\|w_0\|_2^2
\end{equation*}
holds for $t\in [0,T]$.
\end{cor}

The final lemma we need for Proposition \ref{prop2.2} is the following one.

\begin{lem}\label{lem2.9}
Under the same assumptions of Lemma \ref{lem2.8}.
there is a positive constant $C$ independent on $\varepsilon_2$ such that if $\|w_0\|_2+\delta \le \varepsilon_2$, then
\begin{equation*}
		\|z(t)\|_3^2\le C\|w_0\|_2^2, \quad \forall t\in[0,T].
\end{equation*}
\end{lem}
{\it\bfseries Proof.}
Rewriting the equation \eqref{eq2.8}$_2$ as $z_{xx}-z=w_x$, and squaring this equation, we have
\begin{equation}\label{wtjzxx}
  z_{xx}^2+2z_x^2+z^2=w_x^2+(z_xz)_x.
\end{equation}
Integrating \eqref{wtjzxx} over $\mathbbm{R}_+$, combining it with \eqref{wtjwzbj},
and using Cauchy-Schwarz inequality, we get
\begin{equation*}
  \begin{aligned}[b]
  \|z(t)\|_2^2
   &\le \|w_{x}(t)\|^2+C|w_x(0,t)|^2 +\frac{1}{2}\|z(t)\|_\infty^2 \\[2mm]
   &\le \|w_{x}(t)\|^2+C|w_x(t)|_\infty^2 +\frac{1}{2} (\|z(t)\|^2+ \|z_x(t)\|^2) \\[2mm]
   &\le C\|w_{x}(t)\|^2+C \|w_{xx}(t)\|^2+\frac{1}{2} (\|z(t)\|^2+ \|z_x(t)\|^2),
  \end{aligned}
\end{equation*}
which yields $\|z\|_2^2\le C\|w_0\|_2^2$. To get the $L^2$-estimate on $z_{xxx}$, we differentiate \eqref{eq2.8}$_2$ with respect to $x$, then
\begin{equation*}
   \|z_{xxx}(t)\|^2\le 2( \|w_{xx}(t)\|^2+\|z_x(t)\|^2)\le C\|w_0\|_2^2.
\end{equation*}
This completes the proof of Lemma \ref{lem2.9}.
$\hfill\Box$

\subsection{Asymptotic Behavior toward the Stationary Solution}
The global existence of the unique solution for problem \eqref{eq2.8} and its large time behavior is an immediate consequence of Proposition \ref{prop2.2}. Indeed, combining the standard theory of the existence and uniqueness of the local solution with the \emph{a priori} estimates, one can extend the local solution for problem \eqref{eq2.8} globally, that is
\begin{equation*}
\begin{cases}
	w\in C^0([0,\infty);H^2), \quad w_x \in L^2(0,\infty;H^1),\\[1mm]
	z\in C^0([0,\infty);H^3)\cap L^2(0,\infty;H^3).
\end{cases}
\end{equation*}
Then, the \emph{a priori} estimates again assert that
\begin{equation}
\begin{cases}\label{wtjxyfjjg}
	\sup\limits_{t\ge0} (\|w(t)\|_2^2+\|z(t)\|_3^2)< \infty,\\[4mm]
	\int_0^t (\|w_x(\tau)\|_1^2+\|z(\tau)\|_3^2) \,\mathrm{d}\tau < \infty .
\end{cases}
\end{equation}

To complete the proof of Theorem \ref{thm2.2}, by using \eqref{wtjxyfjjg}, we can easily get
\begin{equation}\label{wtjdsjxw1}
	\int_0^\infty \left|\frac{\mathrm{d}}{\mathrm{d}t}\|w_x\|^2 \right| \,\mathrm{d}t < \infty.
\end{equation}
Therefore, it follows from \eqref{wtjxyfjjg} and \eqref{wtjdsjxw1} that
\begin{equation}\label{wtjwxqy0}
	\|w_x\|\rightarrow 0,\ \ {\rm as} \ t \rightarrow 0.
\end{equation}
From the Sobolev inequality, the desired asymptotic behavior in Theorem \ref{thm2.2} can be obtained as
\begin{equation}\label{wtjLwq}
\begin{cases}
\sup\limits_{x\in \mathbbm{R}_+}|w(x,t)|\le \sqrt{2}\|w\|^{\frac{1}{2} }\|w_x\|^{\frac{1}{2} } \rightarrow 0, \ \ {\rm as} \ t \rightarrow 0,\\[4mm]
\sup\limits_{x\in \mathbbm{R}_+}|w_x(x,t)|\le \sqrt{2}\|w_x\|^{\frac{1}{2} }\|w_{xx}\|^{\frac{1}{2} } \rightarrow 0, \ \ {\rm as} \ t \rightarrow 0.
\end{cases}
\end{equation}
The combination of \eqref{wtjwxqy0} and \eqref{wtjLwq} completes the proof of Theorem \ref{thm2.2}.

Here, we give the proof of \eqref{wtjdsjxw1}. In fact, from $\eqref{eq2.82}$, combining $\eqref{wtjxyfjjg}$, we can conclude that
\begin{equation*}
 \begin{aligned}[b]
 	\int_0^\infty \left|\frac{\mathrm{d}}{\mathrm{d}t}\|w_x\|^2 \right| \,\mathrm{d}t
 	&\le C \int_0^\infty\int_{\mathbbm{R}_+} (|w_x^3|+|\bar{u}_{xx}ww_x|) \,\mathrm{d}x \mathrm{d}t	\\
 	&\le C \int_0^\infty (\| w_x(t) \|^2+\| \sqrt{\bar{u}_x}w(t) \|^2 ) \,\mathrm{d}t \\
 	&< \infty.
 \end{aligned}
\end{equation*}
Thus, we finish the proof of \eqref{wtjdsjxw}$_1$ in Theorem \ref{thm2.2}. Finally, according to $\eqref{wtjwzbj}$ and $\eqref{eq2.8}_2$, we set $g=-u_-w_x(0,t)$, $f(x)=-w_x(x,t)$ for any fixed $t\in[0,\infty)$ in Lemma \ref{le-Bessel}. Then, by employing $\eqref{wtjLwq}$, we can obtain the asymptotic behavior of $z$, which completes the proof of Theorem \ref{thm2.2}.

\section{Asymptotics to Superposition of Nonlinear Waves} \label{sec-5}
\subsection{Reformulation of the Problem in the Case of $u_-<0<u_+$}
Referring to the preceding sections, we set
\begin{equation*}
  \Phi_3 (x,t):=\bar{u}_2 (x,t)+\tilde{u}_4 (x,t), \ \ \ \ \Psi_3 (x,t):=\bar{q}_2 (x,t)+\tilde{q}_4 (x,t),
\end{equation*}
as an asymptotic state as $t\rightarrow \infty$, where $(\tilde{u}_4,\tilde{q}_4 )$ and $(\bar{u}_2,\bar{q}_2)$ are given in Lammas \ref{lemxsb} and \ref{lem2.1.1}, respectively. For simplicity, $\bar{u}_2$ and $\tilde{u}_4$ are denoted by $\bar{u}$ and $\tilde{u}$, respectively.
The perturbation
\begin{equation*}
  \begin{cases}
    w (x,t)=u (x,t)-\Phi_3 (x,t)=u (x,t)-\bar{u} (x)-\tilde{u} (x,t),\\[1mm]
    z (x,t)=q (x,t)-\Psi_3 (x,t)=q (x,t)-\bar{q} (x)+\tilde{u}_{x} (x,t),\\
  \end{cases}
\end{equation*}
satisfies the reformulated problem
\begin{equation}\label{fhbrdfc} 
  \begin{cases}
    w_t+(\Phi_3 w)_x+ww_x+z_x=- \bar{u} \tilde{u}_x- \tilde{u}\bar{u}_x,\\[1mm]
    -z_{xx}+z+w_x=-\tilde{u}_{xxx},\\[1mm]
    w(0,t)=0,\\[1mm]
    w(x,0)=w_0(x)=u_0(x)-\bar{u} (x)-\tilde{u} (x,0).
  \end{cases}
\end{equation}

We seek the solutions of \eqref{fhbrdfc} in the set of functions $X_3(0,T)$ defined by
$$X_3(0,T)=\left\{(w,z)|w\in C^0([0,T);H^2)\cap C^1([0,T);H^1),  z\in C^0([0,T);H^3)\cap L^2(0,T;H^3)\right\}.$$
Firstly, we state the global existence and uniform stability result for the reformulated problem \eqref{fhbrdfc}.

\begin{thm}\label{thmfhb} 
Suppose that the boundary condition and far field states satisfy $u_-<0<u_+$, the initial data $w_0\in H^2(\mathbbm{R}_+)$ and the wavelength $\delta=|u_--u_+|$ are sufficiently small. Then there are two positive constants $\varepsilon_1$ and $C=C(\varepsilon_1)$ such that if  $\|w_0\|_2 +\delta\le \varepsilon_1$, the problem \eqref{fhbrdfc} admits a unique solution $(w(x,t),z(x,t))\in X_3(0,+\infty)$ satisfying
\begin{equation*}
\|w(t)\|_2^2+\|z(t)\|_3^2+\int_0^t \ (\|w_x(\tau)\|_1^2+\|z(\tau)\|_3^2 ) \,d{\tau}\le C(\|w_0\|_2^2+\delta^\frac{1}{3} ),
\end{equation*}
and the asymptotic behavior
\begin{equation}\label{fhbdsjxw}
\begin{split}
 &\sup_{x\in \mathbbm{R}_+} |\partial_x^kw(x,t)|\rightarrow 0 \ \  \text{as} \ \ t\rightarrow \infty, \ \ k=0,1, \\[3mm]
 &\sup_{x\in \mathbbm{R}_+} |\partial_x^kz(x,t)|\rightarrow 0 \ \  \text{as} \ \ t\rightarrow \infty, \ \ k=0,1,2.
\end{split}
\end{equation}
\end{thm}

The combination of the local existence and the \emph{a priori} estimates proves Theorem \ref{thmfhb}.

\begin{prop}[Local existence]
Suppose the boundary condition and far field states satisfy $u_-<0<u_+$, the initial data satisfies $w_0 \in H^2(\mathbbm{R}_+)$ and $\|w_0\|_{2}+ \delta\le\varepsilon_1$. Then there are two positive constants $C=C(\varepsilon_1)$ and $T_0=T_0(\varepsilon_1)$ such that the problem \eqref{fhbrdfc} has a unique solution $(w,z)\in X_3(0,T_0)$, which satisfies
\begin{equation*}
\|w(t)\|_2^2+\|z(t)\|_3^2+\int_0^t \  ( \| w_x(\tau) \|_1^2+\| z(\tau) \|_3^2  )  \,d{\tau} \le C(\| w_0 \|_2^2+\delta^\frac{1}{3} ),
\end{equation*}
for $t\in [0,T_0]$.
\end{prop}

\begin{prop}[A priori estimates] \label{propfhbxygj} 
Let $T$ be a positive constant. Suppose that the problem \eqref{fhbrdfc} has a unique solution $(w,z)\in X_3(0,T)$. Then there exists positive constants $\varepsilon_2(\le\varepsilon_1)$ and $C=C(\varepsilon_2)$ such that if $\|w_0\|_2+\delta\le\varepsilon_2$, then we have the estimate
\begin{equation*}
\|w(t)\|_2^2+\|z(t)\|_3^2+\int_0^t \ (\| w_x(\tau)\|_1^2+\| z(\tau)\|_3^2 )  \,d{\tau}
\le C(\|w_0\|_2^2+\delta^\frac{1}{3} ),
\end{equation*}
for $t\in[0,T]$.
\end{prop}

\subsection{\emph{A priori} Estimates}
Under the assumptions of Theorem \ref{thmfhb}, to show the \emph{a priori} estimate in Proposition \ref{propfhbxygj}, we devote ourselves to the estimates on the solution $(w,z)\in X_3(0,T)$ (for some $T>0$) of \eqref{fhbrdfc} under the \emph{a priori} assumption
\begin{equation}\label{fhbxyjs} 
 |w_x(t)|_\infty\le\varepsilon_0,
\end{equation}
where $0<\varepsilon_0\ll 1$. For simplicity, we divide the proof of the \emph{a priori} estimate into several lemmas.

\begin{lem}\label{lemfhbwx}
There are positive constants $\varepsilon_1(\le \varepsilon_0)$ and $C=C(\varepsilon_1)$ such that if $\|w_0\|_2+\delta \le \varepsilon_1$, then
\begin{equation*}
  \begin{aligned}
    \|w(t)\|_1^2+&\int_0^t \left(\|\sqrt{\Phi_{3x}}w(\tau)\|^2+\|\sqrt{\Phi_{3x}}w_x(\tau)\|^2+\|z(\tau)\|_2^2\right) \,\mathrm{d}\tau+\int_0^t w_x^2(0,\tau) \,\mathrm{d}\tau
   \le C(\|w_0\|_1^2+\delta^{\frac{1}{3} })
  \end{aligned}
\end{equation*}
holds for $t\in[0,T]$.
\end{lem}
{\it\bfseries Proof.}
Multiplying \eqref{fhbrdfc}$_1$ by $w$ and \eqref{fhbrdfc}$_2$ by $z$, and adding the two resulting equations up, we obtain
\begin{equation}\label{fhbjbnlgj} 
  \begin{aligned}
    \frac{1}{2}\frac{\mathrm{d}}{\mathrm{d}t}w^2+\frac{1}{2}\Phi_{3x}w^2+z_x^2+z^2
    +\left\{\frac{1}{2} \Phi_3 w^2+\frac{1}{3}w^3-z_xz+zw\right\}_x=- \bar{u} \tilde{u}_xw- \tilde{u}\bar{u}_xw- \tilde{u}_{xxx}z.
  \end{aligned}
\end{equation}

We differentiate \eqref{fhbrdfc}$_1$ with respect to $x$ and multiply it by $w_x$, and multiply \eqref{fhbrdfc}$_1$ by $-z_{xx}$.
Finally adding these two equations up, we obtain
\begin{equation}\label{fhbwx1}
  \begin{aligned}[b]
    \frac{1}{2}\frac{\mathrm{d}}{\mathrm{d}t}&w_x^2+\frac{3}{2} \Phi_{3x} w_x^2+z_{xx}^2+z_x^2+\left\{\frac{1}{2} \Phi_3 w_{x}^2+\frac{1}{2}ww_x^2-z_{x}z+(\bar{u} \tilde{u})_{xx}w \right\}_x\\[3mm]
    &=-\frac{1}{2}w_x^3- \Phi_{3xx}ww_x+\tilde{u}_{xxx}z_{xx}+(\bar{u} \tilde{u})_{xxx}w.
  \end{aligned}
\end{equation}
Adding \eqref{fhbjbnlgj} and \eqref{fhbwx1} up, we get
\begin{equation}\label{fhbx3}
\begin{aligned}[b]
  &\frac{1}{2}\frac{\mathrm{d}}{\mathrm{d}t}(w^2+w_x^2)+\frac{1}{2}\Phi_{3x}w^2
  +\frac{3}{2}\Phi_{3x}w_x^2+z_{xx}^2+2z_x^2+z^2\\[3mm]
  & \ \  +\left\{ \frac{1}{2}\Phi_{3}w^2+\frac{1}{3}w^3+zw+\frac{1}{2} \Phi_3w_x^2+\frac{1}{2}ww_x^2-2z_xz +(\bar{u} \tilde{u})_{xx}w  \right\}_x\\[3mm]
  &=- \bar{u} \tilde{u}_xw- \tilde{u}\bar{u}_xw- \tilde{u}_{xxx}z-\frac{1}{2}w_x^3- \Phi_{3xx}ww_x+\tilde{u}_{xxx}z_{xx}+(\bar{u} \tilde{u})_{xxx}w.
\end{aligned}
\end{equation}
Integrating \eqref{fhbx3} over $\mathbbm{R}_+\times(0,t)$, combining it with $\bar{u}(0)=u_-<0$, $\tilde{u}(0,t)=0$ and $w(0,t)=0$, we have
\begin{equation}\label{fhbx4}
\begin{aligned}[b]
  &\|w(t)\|_1^2+\int_0^t \left(\|\sqrt{\Phi_{3x}}w(\tau)\|^2+\|\sqrt{\Phi_{3x}}w_x(\tau)\|^2+\|z_{xx}(\tau)\|^2+2\|z_x(\tau)\|^2+\|z(\tau)\|^2\right) \,\mathrm{d}\tau-u_-\int_0^t w_x^2(0,t) \,\mathrm{d}\tau\\[3mm]
  &\le C \left(\|w_0\|_1^2+ \int_0^t |z_x(0,t)||z(0,t)| \,\mathrm{d}\tau +\int_0^t\int_{\mathbbm{R}_+}\left(|\bar{u}|\tilde{u}_x|w|+\tilde{u}\bar{u}_x|w|+|\tilde{u}_{xxx}|(|z|+|z_{xx}|)\right)\,\mathrm{d}\tau \right. \\[3mm]
  & \left. \ \ \ \ \ \ \ \  +\int_0^t\int_{\mathbbm{R}_+} (|w_x|^3
 +| \Phi_{3xx}ww_x|+|(\bar{u} \tilde{u})_{xxx}w|)\,\mathrm{d}x \mathrm{d}\tau \right).
\end{aligned}
\end{equation}

Firstly, from \eqref{fhbrdfc}$_1$, we know that $z_x(0,t)=-u_-(w_x(0,t)+\tilde{u}_x(0,t))$ holds.
According to Lemma \ref{lemxsb} $\mathrm{(ii)}$ with $k=1$, for some small $u_-$ corresponding to small $\delta$ satisfying $C u_-^2\le -\frac{u_-}{2}$, we have
\begin{equation}\label{fhbbjgj}
\begin{aligned}[b]
\int_0^t  |z_x(0,\tau)||z(0,\tau)|  \,\mathrm{d}\tau
&\le C \int_0^t   |z_x(0,\tau)|^2 \,\mathrm{d}\tau +\frac{1}{8}\int_0^t    |z(0,\tau)|^2  \,\mathrm{d}\tau\\[3mm]
&\le C u_-^2 \int_0^t ( |w_x(0,\tau)|^2+|\tilde{u}_x(0,\tau)|^2) \,\mathrm{d}\tau +\frac{1}{8} \int_0^t |z(\tau)|_\infty^2 \,\mathrm{d}\tau \\[3mm]
&\le -\frac{u_-}{2}  \int_0^t |w_x(0,\tau)|^2 \,\mathrm{d}\tau +C \delta^\frac{1}{3} (1+t)^{-1} +\frac{1}{8} \int_0^t \|z_x(\tau)\|^2 \,\mathrm{d}\tau +\frac{1}{8} \int_0^t  \|z(\tau)\|^2 \,\mathrm{d}\tau.
\end{aligned}
\end{equation}

Secondly, we estimate the right-hand side of \eqref{fhbx4} as follows:
\begin{equation*}
  \int_{\mathbbm{R}_+} |\bar{u}| \tilde{u}_x|w| \,\mathrm{d}x= \int_{\mathbbm{R}_+} (-\bar{u}) \tilde{u}_x|w| \,\mathrm{d}x
  =\int_{0}^{u_+t} (-\bar{u}) \tilde{u}_x|w| \,\mathrm{d}x +\int_{u_+t} ^\infty (-\bar{u}) \tilde{u}_x|w| \,\mathrm{d}x
  =:I_1+I_2.
\end{equation*}
By virtue of $\bar{u}<0, \tilde{u}>0$, using Lemmas \ref{lem2.1.1} and \ref{lemxsb}, we get
\begin{equation}\label{fhbu1}
  \begin{aligned}[b]
    I_1&\le |w|_\infty |\tilde{u}_x|_\infty \int_{0}^{u_+t} (- \bar{u})  \,\mathrm{d}x \\[3mm]
    &\le C \|w_x(t)\|^{\frac{1}{2} }\|w(t)\|^{\frac{1}{2} }\delta^{\frac{1}{4}}(1+t)^{-\frac{7}{8} }\int_{0}^{u_+t} \frac{\delta}{1+\delta x}  \,\mathrm{d}x \\[3mm]
    &\le \frac{1}{48}  \|w_x(t)\|^2+ C \delta^{\frac{1}{3} } \|w(t)\|^{\frac{2}{3} }\left\{(1+t)^{-\frac{7}{8}}\ln(1+t)\right\}^{\frac{4}{3} }\\[3mm]
    &\le  \frac{1}{48}  \|w_x(t)\|^2+C \delta^{\frac{1}{3} } (1+t)^{-\frac{7}{6}} ( \ln(1+t) )^{\frac{4}{3} },
  \end{aligned}
\end{equation}
and
\begin{equation}\label{fhbu2}
  \begin{aligned}[b]
    I_2&\le |w|_\infty \left\{ -[\bar{u}\tilde{u}]_{u_+t}^\infty +\int_{u_+t}^{\infty}  \bar{u}_x\tilde{u}  \,\mathrm{d}x \right\}\\[3mm]
    &\le C u_+ \|w_x(t)\|^{\frac{1}{2} }\|w(t)\|^{\frac{1}{2} }\int_{u_+t}^{\infty}  \bar{u}_x  \,\mathrm{d}x\\[3mm]
    &\le C \delta \|w_x(t)\|^{\frac{1}{2} }\|w(t)\|^{\frac{1}{2} }\int_{u_+t}^{\infty}  \frac{\delta^2}{(1+ \delta x)^{2}}   \,\mathrm{d}x\\[3mm]
    &\le \frac{1}{48}  \|w_x(t)\|^2+ C \delta^{\frac{8}{3} }\|w(t)\|^{\frac{2}{3} }(1+t)^{-\frac{4}{3} }\\[3mm]
    &\le \frac{1}{48}  \|w_x(t)\|^2+ C \delta^{\frac{1}{3} }(1+t)^{-\frac{4}{3} }.
  \end{aligned}
\end{equation}
In deriving the second inequality of $\eqref{fhbu1}$, we have used from Lemma \ref{lemxsb} $\mathrm{(iv)}$ and $\mathrm{(v)}$ that
 $$|\tilde{u}_x|_\infty\le C(\delta(1+t)^{-\frac{1}{2} })^\frac{1}{4}((1+t)^{-1})^{\frac{3}{4} }\le C \delta^\frac{1}{4}(1+t)^{-\frac{7}{8} }.$$
 And in deriving the last inequality of $\eqref{fhbu1}$ and $\eqref{fhbu2}$, we have used the fact that $\| w(t) \| $ is bounded. In a similar fashion to \eqref{fhbu1} and \eqref{fhbu2}, we can obtain
\begin{equation}\label{fhbuu}
  \begin{aligned}
    \int_{\mathbbm{R}_+} \bar{u}_x \tilde{u}|w| \,\mathrm{d}x
    \le \frac{1}{24}  \|w_x\|^2+C \delta^{\frac{1}{3} } \left\{ (1+t)^{-\frac{4}{3} } +(1+t)^{-\frac{7}{6}}( \ln(1+t) )^{\frac{4}{3} }\right\}.
  \end{aligned}
\end{equation}
Using Cauchy-Schwarz inequality and Lemma \ref{lemxsb} $\mathrm{(v)}$ with $k=3, l=0$, we get
\begin{equation}\label{eq-uxxxzxx}
  \int_{\mathbbm{R}_+} |\tilde{u}_{xxx}(z+z_{xx})|  \,\mathrm{d}x
  \le \frac{1}{8} (\|z_{xx}(t)\|^2+\|z(t)\|^2)+ C\delta^\frac{1}{3} (1+t)^{-\frac{5}{2} },
\end{equation}
and
\begin{equation}\label{eq-uuw}
\begin{aligned}[b]
      \int_{\mathbbm{R}_+} |(\bar{u} \tilde{u})_{xxx}w| \,\mathrm{d}x
      &\le C \int_{\mathbbm{R}_+} (|\bar{u}_{xxx}|\tilde{u}|w|+|\bar{u}_{xx}\tilde{u}_xw|+|\bar{u}_x \tilde{u}_{xx}w|+|\bar{u}\tilde{u}_{xxx}w|) \,\mathrm{d}x\\[3mm]
      &\le C \left|\frac{\bar{u}_{xxx}}{\bar{u}_{x}} \right|_\infty \int_{\mathbbm{R}_+} \tilde{u}\bar{u}_{x}|w| \,\mathrm{d}x
      +C |w|_\infty( |\bar{u}_{xx}|_1|\tilde{u}_x|_\infty +|\bar{u}_{x}|_1|\tilde{u}_{xx}|_\infty +|\bar{u}|_\infty |\tilde{u}_{xxx}|_1 )\\[3mm]
      & \le C \delta |w|_\infty \int_{\mathbbm{R}_+} \tilde{u}\bar{u}_{x} \,\mathrm{d}x+C\delta|w|_\infty (1+t)^{-1}\\[3mm]
      &\le \frac{1}{24} \|w_x(t)\|^2+C \delta^{\frac{1}{3} } \left\{ (1+t)^{-\frac{4}{3} } +(1+t)^{-\frac{7}{6}}( \ln(1+t) )^{\frac{4}{3} }\right\}.
\end{aligned}
\end{equation}
From the Lemmas \eqref{lem2.1.1} and \eqref{lemxsb}, we can get estimates of the remaining terms on the right-hand side of \eqref{fhbx4}:
\begin{equation}\label{eq-wx3}
  \begin{aligned}
          \int_{\mathbbm{R}_+} |w_x|^3 \,\mathrm{d}x
          \le |w_x|_\infty\|w_x(t)\|^2
          \le \varepsilon_0\|w_x(t)\|^2,
  \end{aligned}
\end{equation}
and
\begin{equation}\label{fhbwx3}
  \begin{aligned}[b]
          \int_{\mathbbm{R}_+} |\Phi_{3xx}ww_x| \,\mathrm{d}x
          &\le \int_{\mathbbm{R}_+} (|\tilde{u}_{xx}ww_x|+|\bar{u}_{xx}ww_x|) \,\mathrm{d}x\\[3mm]
          &\le \frac{1}{24} \|w_x(t)\|^2+C|w|_\infty \delta^\frac{1}{3} (1+t)^{-\frac{3}{2} }+C \left| \frac{\bar{u}_{xx}^2}{\bar{u}_{x}}\right|_\infty \|\sqrt{\bar{u}_x }w(t)\|^2.\\
  \end{aligned}
\end{equation}
On the other hand, we note that $w_x^2\le 3(z_{xx}^2+z^2+\tilde{u}_{xxx}^2)$ due to \eqref{fhbrdfc}$_2$, which yields
\begin{equation}\label{fhbwz}
   \|w_x(t)\|^2\le 3(\|z_{xx}(t)\|^2+\|z(t)\|^2+\|\tilde{u}_{xxx}(t)\|^2).
\end{equation}
In the end, substituting \eqref{fhbbjgj}-\eqref{fhbwx3} into \eqref{fhbx4}, and using \eqref{fhbwz}, for some small $\delta$ and $\varepsilon_0(<\frac{1}{8} )$, we can conclude that
\begin{equation*}
  \begin{aligned}
    \|w(t)\|_1^2+\int_0^t \left(\|\sqrt{\Phi_{3x}}w(\tau)\|^2+\|\sqrt{\Phi_{3x}}w_x(\tau)\|^2+\|z(\tau)\|_2^2\right) \,\mathrm{d}\tau-u_-\int_0^t w_x^2(0,\tau) \,\mathrm{d}\tau
   \le C(\|w_0\|_1^2+\delta^\frac{1}{3} ).
  \end{aligned}
\end{equation*}
This completes the proof of Lemma \ref{lemfhbwx}.
$\hfill\Box$

From \eqref{fhbwz} and the results of Lemma \ref{lemfhbwx}, we can easily get the following estimate.

\begin{cor}
Under the assumptions of Lemma \ref{lemfhbwx}, the estimate
\begin{equation*}
  \int_0^t \|w_x(\tau)\|^2 \,\mathrm{d}\tau \le C(\|w_0\|_1^2+\delta^{\frac{1}{3} })
\end{equation*}
holds for $t\in[0,T]$.
\end{cor}

\begin{lem}\label{lemfhbwxx}
There are positive constants $\varepsilon_2(\le \varepsilon_1)$ and $C=C(\varepsilon_2)$
such that if $\|w_0\|_2+\delta \le \varepsilon_2$, then
\begin{equation}\label{fhbwxx}
    \|w_{xx}(t)\|^2+\int_0^t \left(\|\sqrt{\Phi_{3x}}w_{xx}(\tau)\|^2+\|z_{xxx}(\tau)\|^2+\|z_{xx}(\tau)\|^2+w_{xx}^2(0,\tau)\right) \,\mathrm{d}\tau \le C(\|w_0\|_2^2+\delta^{\frac{1}{3} }),
\end{equation}
for $t\in [0,T]$.
\end{lem}
{\it\bfseries Proof.}
We differentiate \eqref{fhbrdfc}$_1$ twice with respect to $x$ and multiply it by $w_{xx}$,
and differentiate \eqref{fhbrdfc}$_2$ with respect to $x$ and multiply it by $-z_{xxx}$.
Then, adding these two equations up, we obtain
\begin{equation}\label{fhbwxx1}
    \begin{aligned}[b]
      \frac{1}{2}\frac{\mathrm{d}}{\mathrm{d}t}w_{xx}^2&+\frac{5}{2}\Phi_{3x}w_{xx}^2+z_{xxx}^2+z_{xx}^2
      +\left\{\frac{1}{2} \Phi_3 w_{xx}^2+\frac{1}{2}ww_{xx}^2-z_{xx}z_x+(\bar{u} \tilde{u})_{xxx}w_x\right\}_x \\[2mm]
      &=- \Phi_{3xxx}ww_{xx}-3 \Phi_{3xx}w_xw_{xx}-\frac{5}{2}w_xw_{xx}^2+\partial_x^4\tilde{u}z_{xxx}+\partial_x^4(\bar{u} \tilde{u})w_x.
    \end{aligned}
\end{equation}
Integrating \eqref{fhbwxx1} over $\mathbbm{R}_+\times(0,t)$, combining $\Phi_3 (0,t)=u_-<0$ and $w(0,t)=0$, we have
\begin{equation}\label{fhbwxxjf}
    \begin{aligned}[b]
      &\|w_{xx}(t)\|^2+\int_0^t \left(\|\sqrt{\Phi_{3x}}w_{xx}(\tau)\|^2+\|z_{xxx}(\tau)\|^2+\|z_{xx}(\tau)\|^2\right) \,\mathrm{d}\tau -u_-\int_0^t w_{xx}^2(0,\tau) \,\mathrm{d}\tau \\[3mm]
      &\le C\left( \|w_0\|_2^2 +\int_0^t |z_{xx}(0,\tau)||z_{x}(0,\tau)| \,\mathrm{d}\tau +\int_0^t  |(\bar{u} \tilde{u})_{xxx}(0,\tau)||w_x(0,\tau)| \,\mathrm{d}\tau \right. \\[3mm]
      &\left. \ \ \ \ +\int_0^t\int_{\mathbbm{R}_+} (|\Phi_{3xxx}ww_{xx}|+ |\Phi_{3xx}w_xw_{xx}|+|w_xw_{xx}^2|+ |\partial_x^4\tilde{u}z_{xxx}|+|\partial_x^4(\bar{u} \tilde{u})w_x|) \,\mathrm{d}x \mathrm{d}\tau   \right).
    \end{aligned}
\end{equation}

Firstly, we estimate the second and the third terms on the right-hand side of \eqref{fhbwxxjf}.
Using the Cauchy-Schwarz inequality, we have
\begin{equation}\label{fhbwxxbj}
    \begin{aligned}[b]
      C \int_0^t |z_{xx}(0,\tau)||z_{x}(0,\tau)| \,\mathrm{d}\tau
      &\le C \int_0^t |z_{xx}(\tau)|_\infty |z_{x}(\tau)|_\infty  \,\mathrm{d}\tau \\[3mm]
      &\le \frac{1}{4} \int_0^t \|z_{xxx}(\tau)\|^2 \,\mathrm{d}\tau +C \int_0^t (\|z_{xx}(\tau)\|^2+\|z_{x}(\tau)\|^2) \,\mathrm{d}\tau \\[3mm]
      &\le \frac{1}{4} \int_0^t \|z_{xxx}(\tau)\|^2 \,\mathrm{d}\tau +C (\|w_0\|_1^2+\delta^\frac{1}{3} ).
    \end{aligned}
\end{equation}
For $(\bar{u} \tilde{u})_{xxx}=\bar{u}_{xxx} \tilde{u}+3 \bar{u}_{xx}\tilde{u}_x+3 \bar{u}_x \tilde{u}_{xx}+\bar{u} \tilde{u}_{xxx}$ and $\tilde{u}(0,t)=0$, combining Lemma \ref{lemfhbwx} and Lemma \ref{lemxsb} $\mathrm{(ii)}$, we get
\begin{equation*}
    \begin{aligned}[b]
      C \int_0^t  |(\bar{u} \tilde{u})_{xxx}(0,\tau)||w_x(0,\tau)| \,\mathrm{d}\tau
      &\le C\delta \int_0^t  (|\tilde{u}_x(0,\tau)|+|\tilde{u}_{xx}(0,\tau)|+|\tilde{u}_{xxx}(0,\tau)|)|w_x(0,\tau)| \,\mathrm{d}\tau\\[3mm]
      &\le \int_0^t w_x^2(0,\tau) \,\mathrm{d}\tau + C \delta^\frac{1}{3}   (1+t)^{-1 }.
    \end{aligned}
\end{equation*}

Next, we estimate the last five terms on the right-hand side of \eqref{fhbwxxjf}. Using Cauchy-Schwarz inequality, we have
\begin{equation*}
\begin{aligned}[b]
  \int_{\mathbbm{R}_+} |\Phi_{3xxx}ww_{xx}| \,\mathrm{d}x
  &\le  \int_{\mathbbm{R}_+}  (|\bar{u}_{xxx}ww_{xx}|+|\tilde{u}_{xxx} ww_{xx}|) \,\mathrm{d}x\\[3mm]
  &\le \frac{1}{36} \|w_{xx}(t)\|^2+ C \left|\frac{\bar{u}_{xxx}^2}{\bar{u}_{x}} \right|_\infty \|\sqrt{\bar{u}_{x}}w(t)\|^2+C \delta^\frac{1}{3} |w(t)|_\infty (1+t)^{-\frac{5}{2} },
\end{aligned}
\end{equation*}
\begin{equation*}
\begin{aligned}
 \int_{\mathbbm{R}_+} |\Phi_{3xx}w_xw_{xx}|  \,\mathrm{d}x  \le \frac{1}{36} \|w_{xx}(t)\|^2+ |\Phi_{3xx}(t)|_\infty \|w_x(t)\|^2\le \frac{1}{36} \|w_{xx}(t)\|^2+C \delta\|w_x(t)\|^2,
\end{aligned}
\end{equation*}
\begin{equation*}
\begin{aligned}
 \int_{\mathbbm{R}_+} |\partial_x^4\tilde{u}z_{xxx}| \,\mathrm{d}x \le \frac{1}{4} \|z_{xxx}(t)\|^2+  C\delta^\frac{1}{3} (1+t)^{-\frac{7}{2} },
\end{aligned}
\end{equation*}
and by utilizing the \emph{a priori} assumption \eqref{fhbxyjs},
\begin{equation*}
\begin{aligned}
  \int_{\mathbbm{R}_+} |w_xw_{xx}^2|  \,\mathrm{d}x \le |w_x(t)|_\infty \|w_{xx}(t)\|^2\le \varepsilon_0\|w_{xx}(t)\|^2.
\end{aligned}
\end{equation*}
By the same method as \eqref{fhbu1} and \eqref{fhbu2}, we can get the last one estimate:
\begin{equation}\label{fhbwxxzh}
\begin{aligned}[b]
   \int_{\mathbbm{R}_+} |\partial_x^4(\bar{u} \tilde{u})w_x|  \,\mathrm{d}x
   &\le C \int_{\mathbbm{R}_+} (|\partial_x^4\bar{u}\tilde{u}|+|\bar{u}_{xxx}\tilde{u}_x|+|\bar{u}_{xx}\tilde{u}_{xx}|
   +|\bar{u}_{x}\tilde{u}_{xxx}|+|\bar{u}\partial_x^4\tilde{u}|)|w_x|  \,\mathrm{d}x \\[3mm]
   &\le C \left|\frac{\partial_x^4\bar{u}}{\bar{u}_{x}} \right|_\infty |w_x|_\infty \int_{\mathbbm{R}_+} \bar{u}_{x}\tilde{u}\,\mathrm{d}x +C \delta|w_x|_\infty  (|\tilde{u}_x|_\infty+|\tilde{u}_{xx}|_\infty+|\tilde{u}_{xxx}|_\infty+|\partial_x^4\tilde{u}|_\infty ) \\[3mm]
   &\le C \delta|w_x|_\infty\left( \int_{0} ^{u_+t} \bar{u}_{x}\tilde{u}\,\mathrm{d}x +\int_{u_+t}^{\infty} \bar{u}_{x}\tilde{u}  \,\mathrm{d}x \right)  +C\delta|w_x|_\infty(1+t)^{-1}\\[3mm]
   &\le \frac{1}{36} \|w_{xx}(t)\|^2 +C \delta^\frac{1}{3} \|w_x(t)\|^\frac{2}{3}  \left\{ (1+t)^{-\frac{4}{3} }+(1+t)^{-\frac{7}{6} }(\ln(1+t))^{\frac{4}{3} } \right\}.
\end{aligned}
\end{equation}

In the end, to estimate $\int_0^t \|w_{xx}(\tau)\|^2 \,\mathrm{d}\tau $, for \eqref{fhbrdfc}$_2$, we note that
\begin{equation}\label{fhbwz2}
  \|w_{xx}(t)\|^2\le 3(\|z_{xxx}(t)\|^2+\|z_x(t)\|^2+\|\partial_x^4\tilde{u}(t)\|^2).
\end{equation}
Substituting \eqref{fhbwxxbj}-\eqref{fhbwxxzh} into \eqref{fhbwxxjf}, combining \eqref{fhbwz2} and Lemma \ref{lemfhbwx}, we have
\begin{equation*}
\begin{aligned}[b]
  &\|w_{xx}(t)\|^2+\int_0^t \left(\|\sqrt{\Phi_{3x}}w_{xx}(\tau)\|^2+\|z_{xxx}(\tau)\|^2+\|z_{xx}(\tau)\|^2 + w_{xx}^2(0,\tau)\right) \,\mathrm{d}\tau\\[3mm]
  \le& \left(\frac{3}{4}+ \varepsilon_0\right) \|z_{xxx}(t)\|^2+C(\|w_0\|_2^2 +\delta^{\frac{1}{3} }).
\end{aligned}
\end{equation*}
That is, for $\varepsilon_0<\frac{1}{4} $, we complete the proof of \eqref{fhbwxx}.
$\hfill\Box$

The equation \eqref{fhbwz2} and Lemma \ref{lemfhbwxx} yields Corollary \ref{corfhbwxx}.

\begin{cor}\label{corfhbwxx}
Under the assumptions of Lemma \ref{lemfhbwxx},
\begin{equation*}
  \int_0^t \|w_{xx}(\tau)\|^2 \,\mathrm{d}\tau \le C (\|w_0\|_2^2 +\delta^{\frac{1}{3}}), \quad \forall t\in[0,T].
\end{equation*}
\end{cor}

\begin{lem}\label{lemfhbz}
Under the assumptions of Lemma \ref{lemfhbwxx},
\begin{equation}\label{fhbzH3}
    \|z(t)\|_3^2\le C(\|w_0\|_2^2+\delta^{\frac{1}{3}}), \quad \forall t\in[0,T].
\end{equation}
\end{lem}
{\it\bfseries Proof.}
Firstly, from \eqref{fhbrdfc}$_1$, we have the equation at the boundary $x=0$,
\begin{equation*}
    z_x(0,t)=-u_-(w_x(0,t)+\tilde{u}_x(0,t)),
\end{equation*}
which plays an essential role in estimating boundary terms. From \eqref{fhbrdfc}$_2$, it holds
\begin{equation}\label{fhbzpf}
   z_{xx}^2+z^2+2z_x^2=w_x^2+\tilde{u}_{xxx}^2+2(z_xz)_x+2w_x\tilde{u}_{xxx}.
\end{equation}
Integrating \eqref{fhbzpf} over $\mathbbm{R}_+$, by Cauchy-Schwarz inequality, we get
\begin{equation}\label{fhbzxx}
  \begin{aligned}[b]
     \|z_{xx}(t)\|^2+2\|z_x(t)\|^2+\|z(t)\|^2
    &\le 2\|w_x(t)\|^2+2\|\tilde{u}_{xxx}(t)\|^2+2|z_x(0,t)||z(0,t)| \\[2mm]
    &\le C(\|w_0\|_2^2+\delta^{\frac{1}{3} }) +C|z_x(0,t)|^2+\frac{1}{2} |z(t)|_\infty^2\\[2mm]
    &\le C(\|w_0\|_2^2+\delta^{\frac{1}{3} }) +C(|w_x(t)|_\infty^2 +|\tilde{u}_x(t)|_\infty^2)+\frac{1}{2}  \|z_x(t)\|^2+\frac{1}{2}  \|z(t)\|^2 \\[2mm]
    &\le  C(\|w_0\|_2^2+\delta^{\frac{1}{3} }) +\frac{1}{2}  \|z_x(t)\|^2+\frac{1}{2}  \|z(t)\|^2.
  \end{aligned}
\end{equation}
Differentiating \eqref{fhbrdfc}$_2$ with respect to $x$ and integrating the resulting equation over $\mathbbm{R}_+$,
combining \eqref{fhbzxx}, Lemma \ref{lemxsb} and Lemma \ref{lemfhbwxx}, we get
\begin{equation}\label{fhbzxxx}
  \begin{aligned}
    \|z_{xxx}(t)\|^2
    \le C(\|z_x(t)\|^2+\|w_{xx}(t)\|^2+\|\partial_x^4\tilde{u}(t)\|^2)
    \le C(\|w_0\|_2^2+\delta^{\frac{1}{3} }).
  \end{aligned}
\end{equation}
The combination of \eqref{fhbzxx} and \eqref{fhbzxxx} completes the proof of Lemma \ref{lemfhbz}.
$\hfill\Box$

\subsection{Asymptotic Behavior toward the Superposition of Nonlinear Waves}
Once the \emph{a priori} estimates is established, by combining the local existence, the global existence of unique solution of \eqref{fhbrdfc} and its asymptotic behavior are easily obtained. That is, the global in time solution
\begin{equation*}
\begin{cases}
	w\in C^0([0,\infty);H^2), \quad w_x \in L^2(0,\infty;H^1),\\[1mm]
	z\in C^0([0,\infty);H^3)\cap L^2(0,\infty;H^3).
\end{cases}
\end{equation*}
Then, the \emph{a priori} estimates again assert that
\begin{equation}
\begin{cases}\label{fhbxygjjg}
	\sup\limits_{t\ge0} (\|w(t)\|_2^2+\|z(t)\|_3^2)< \infty,\\[1mm]
	\int_0^t (\|w_x\|_1^2+\|z\|_3^2) \,\mathrm{d}\tau < \infty .
\end{cases}
\end{equation}

To complete the proof of Theorem \ref{thmfhb}, we need to show that
\begin{equation}\label{fhbdsjxw1}
	\int_0^\infty \left|\frac{\mathrm{d}}{\mathrm{d}t}\|w_x\|^2 \right| \,\mathrm{d}t < \infty,
\end{equation}
 it follows
\begin{equation}\label{fhbwxqy0}
	\|w_x\|\rightarrow 0,\ \ ( t \rightarrow 0).
\end{equation}
Using the Sobolev inequality, we can obtain the desired asymptotic behavior in Theorem \ref{thmfhb}
\begin{equation}\label{fhbLwq}
\begin{cases}
\sup\limits_{x\in \mathbbm{R}_+}|w(x,t)|\le \sqrt{2}\|w\|^{\frac{1}{2} }\|w_x\|^{\frac{1}{2} } \rightarrow 0, \ \ (t \rightarrow 0),\\[4mm]
\sup\limits_{x\in \mathbbm{R}_+}|w_x(x,t)|\le \sqrt{2}\|w_x\|^{\frac{1}{2} }\|w_{xx}\|^{\frac{1}{2} } \rightarrow 0, \ \ (t \rightarrow 0).
\end{cases}
\end{equation}
The combination of \eqref{fhbwxqy0} and \eqref{fhbLwq} can completes the proof of Theorem \ref{thmfhb}.

The proof of \eqref{fhbdsjxw1} can be easily obtained. In fact, using the similar estimatea as $\eqref{eq-wx3}$, $\eqref{fhbwx3}$, $\eqref{eq-uxxxzxx}$ and $\eqref{eq-uuw}$, combining \eqref{fhbxygjjg}, we can get from $\eqref{fhbwx1}$ that
\begin{equation*}
 \begin{aligned}[b]
 	\int_0^\infty \left|\frac{\mathrm{d}}{\mathrm{d}t}\|w_x\|^2 \right| \,\mathrm{d}t
 	&\le C \int_0^\infty |z_x(0,t)z(0,t)| \,\mathrm{d}t + C \int_0^\infty\int_{\mathbbm{R}_+} (|w_x^3|+|\Phi_{3xx}ww_x|+|\tilde{u}_{xxx}z_{xx}|+|(\bar{u} \tilde{u})_{xxx}w|)\,\mathrm{d}x \mathrm{d}t	\\[3mm]
 	&<\infty.
 \end{aligned}
\end{equation*}

Thus, we get the asymptotic behavior of $w$ and we finish the proof of \eqref{fhbdsjxw}$_1$ in Theorem \ref{thmfhb}. In the end, using the Lemma \ref{le-Bessel}, and setting $g=-u_-(w_x(0,t)+\tilde{u}_x(0,t))$, $f(x)=-w_x(x,t)-\tilde{u}_{xxx}(x,t)$ for any fixed $t\in[0,\infty)$, we can obtain the asymptotic behavior of $z$, which completes the proof of Theorem \ref{thmfhb}.

\vspace{6mm}

\noindent {\bf Acknowledgements:}
The research was supported by
Guangdong Basic and Applied Basic Research Foundation $\#$2020B1515310015, $\#$2021A1515010367,
the National Natural Science Foundation of China $\#$11771150, $\#$11831003.

\end{document}